\documentclass[10pt, times, twocolumn]{article}
\usepackage[top=2cm, bottom=2cm, left=2cm, right=2cm]{geometry}
\usepackage{graphics}
\usepackage{amssymb}
\usepackage{amsfonts,amsmath,mathtools,amsthm}
\usepackage{wrapfig}
\usepackage{mathrsfs}
\usepackage{xcolor}
\usepackage{mathtools}
\usepackage{enumitem}
\usepackage{stmaryrd}
\usepackage[numbers]{natbib}
\providecommand{\keywords}[1]{\emph{Keywords}: #1}

\setlist{nolistsep}

\let\div\undefined
\DeclareMathOperator{\div}{div}								
\DeclareMathOperator{\diag}{diag}								
\DeclareMathOperator{\deviatoric}{dev}							
\newcommand{\tp}{^T}										
\newcommand{\diff}[2]{\frac{\mathrm{d} #1}{\mathrm{d} #2}} 			
\newcommand{\pdiff}[2]{\frac{\partial #1}{\partial #2}} 				
\newcommand{\pdiffarg}[2]{\partial_{ #1}{ #2}} 						
\newcommand{\nablaarg}[2]{\nabla_{ #1}{ #2}} 						
\DeclarePairedDelimiter\norm{\lVert}{\rVert} 						
\newcommand{\var}{\delta}									
\newcommand{\bdot}{:}										
\newcommand{\duality}[4]{\!\left< #1, #2 \right>_{#3\times#4}}			
\newcommand{\bduality}[2]{\!\left< #1, #2 \right>}					
\newcommand{\dual}{^{*}}										
\newcommand{\conj}{^{*}}										
\newcommand{\sym}[1]{{#1}\txtsup{(s)}}								

\newcommand{\rbr}[1]{\left( #1 \right)} 							
\newcommand{\sbr}[1]{\left[ #1 \right]} 							
\newcommand{\cbr}[1]{\left\{ #1 \right\}} 							
\newcommand{\fixed}[2]{\left. #1 \right|_{#2}} 						
\newcommand{\abs}[1]{\lvert #1 \rvert} 				
\newcommand{\txtsub}[1]{_{\text{#1}}} 							
\newcommand{\txtsup}[1]{^{\text{#1}}} 							
\newcommand{\elem}{_e} 										
\newcommand{\qp}{_d} 										
\newcommand{\qppnt}{d} 										
\newcommand{\gp}{_g} 										
\newcommand{\interp}{^h} 									
\newcommand{\pnt}{\bbx}										
\newcommand{\force}{\bq}									
\newcommand{\stiff}{\bK}										
\newcommand{\internal}{\txtsup{int}}								
\newcommand{\external}{\txtsup{ext}}							

\newcommand{\ba}{\boldsymbol{a}}

\newcommand{\be}{\boldsymbol{e}}
\newcommand{\bbf}{\boldsymbol{f}}
\newcommand{\bg}{\boldsymbol{g}}

\newcommand{\bl}{\boldsymbol{l}}

\newcommand{\bn}{\boldsymbol{n}}
\newcommand{\bq}{\boldsymbol{q}}
\newcommand{\br}{\boldsymbol{r}}
\newcommand{\bs}{\boldsymbol{s}}

\newcommand{\bu}{\boldsymbol{u}}
\newcommand{\bv}{\boldsymbol{v}}
\newcommand{\bbx}{\boldsymbol{x}}
\newcommand{\bw}{\boldsymbol{w}}
\newcommand{\by}{\boldsymbol{y}}

\newcommand{\balpha}{\boldsymbol{\alpha}}
\newcommand{\bbbeta}{\boldsymbol{\beta}}
\newcommand{\bbeta}{\boldsymbol{\eta}}
\newcommand{\beps}{\boldsymbol{\varepsilon}}

\newcommand{\blambda}{\boldsymbol{\lambda}}

\newcommand{\bsigma}{\boldsymbol{\sigma}}
\newcommand{\btau}{\boldsymbol{\tau}}

\newcommand{\bphi}{\boldsymbol{\varphi}}
\newcommand{\bxi}{\boldsymbol{\xi}}

\renewcommand{\phi}{\varphi}
\newcommand{\bzero}{\boldsymbol{0}}

\newcommand{\bB}{\boldsymbol{B}}
\newcommand{\bC}{\boldsymbol{C}}

\newcommand{\bE}{\boldsymbol{E}}
\newcommand{\bF}{\boldsymbol{F}}
\newcommand{\bG}{\boldsymbol{G}}
\newcommand{\bH}{\boldsymbol{H}}
\newcommand{\bI}{\boldsymbol{I}}
\newcommand{\bJ}{\boldsymbol{J}}
\newcommand{\bK}{\boldsymbol{K}}
\newcommand{\bL}{\boldsymbol{L}}

\newcommand{\bN}{\boldsymbol{N}}
\newcommand{\bP}{\boldsymbol{P}}
\newcommand{\bQ}{\boldsymbol{Q}}
\newcommand{\bR}{\boldsymbol{R}}
\newcommand{\bS}{\boldsymbol{S}}
\newcommand{\bT}{\boldsymbol{T}}

\newcommand{\bV}{\boldsymbol{V}}
\newcommand{\bW}{\boldsymbol{W}}

\newcommand{\bPi}{\boldsymbol{\Pi}}

\newcommand{\cA}{\mathcal{A}}

\newcommand{\cE}{\mathcal{E}}

\newcommand{\cK}{\mathcal{K}}

\newcommand{\cL}{\mathcal{L}}
\newcommand{\cR}{\mathcal{R}}

\newcommand{\fC}{\mathbb{C}}

\newcommand{\dmn}{\Omega}									
\newcommand{\bdy}{\Gamma}									
\newcommand{\side}{\gamma}									
\newcommand{\node}{V}										
\newcommand{\intbr}[1]{\int_{\dmn} #1  \text{d}\dmn}					
\newcommand{\intbdybr}[1]{\int_{\bdy} #1 \, \text{d}\bdy} 				
\newcommand{\intelembr}[1]{\int_{\dmn_e} #1  \text{d}\dmn} 			
\newcommand{\displ}{\bu}										
\newcommand{\stn}{\beps}									
\newcommand{\sts}{\bsigma}									
\newcommand{\pmlt}{\Delta\lambda}								
\newcommand{\gf}{\cL}										
\newcommand{\hw}{\cE}										
\newcommand{\hpr}{\Pi}										
\newcommand{\displprm}{\ba}   								
\newcommand{\stnprm}{\be}									
\newcommand{\stsprm}{\bbbeta}								
\newcommand{\pmltprm}{\Delta\bl}								
\newcommand{\displmtx}{\bN}   								
\newcommand{\graddisplmtx}{\bB}   								
\newcommand{\stnmtx}{\bE}									
\newcommand{\stsmtx}{\bP}									
\newcommand{\pmltmtx}{\bL}									
\newcommand{\avgstsmtx}{\overline{\stsmtx}}						
\newcommand{\compmtx}{\bQ}									
\newcommand{\stsstnmtx}{\bG}									
\newcommand{\cmplmtx}{\bH}									
\newcommand{\yfvct}{\bF}										
\newcommand{\flowmtx}{\bR}									
\newcommand{\forces}{\bq}									
\newcommand{\res}{\br}										
\newcommand{\elflex}{\bH}									
\newcommand{\elstiff}{\bK}									
\newcommand{\act}{\txtsub{act}}								
\newcommand{\nn}{N\txtsub{N}}								

\newcommand{\proj}{\bPi}										
\newcommand{\lift}{\bV}										
\newcommand{\fltdisplprm}{\overline{\ba}}	   						
\newcommand{\frdisplprm}{\blambda}   							
\newcommand{\ie}{\Phi}										
\newcommand{\iemin}{\hat{\Phi}}								
\newcommand{\helmholtz}{\varphi}								
\newcommand{\gibbs}{\varphi\conj}								
\newcommand{\diss}{D}										
\newcommand{\yf}{F}										
\newcommand{\intvarstn}{\mathscr{I}}							
\newcommand{\intvarsts}{\mathscr{F}}							
\newcommand{\eps}{\varepsilon}

\newcommand{\p}{^\text{p}}
\newcommand{\e}{^\text{e}}
\newcommand{\ep}{^\text{ep}}
\newcommand{\trial}{^\text{tr}}

\newcommand{\n}{_{n}}
\newcommand{\npo}{_{n+1}}
\newcommand{\kin}{\txtsup{k}}
\newcommand{\iso}{\txtsup{i}}

\newcommand{\y}{_{\text{y}_0}}

\theoremstyle{remark}
\newtheorem*{remark}{Remark}

\title{An overview of mixed finite elements \\ for the analysis of inelastic bidimensional structures}
\date{July 9, 2018}
\author{
	Nicola A.~Nodargi\thanks{E-mail: nodargi@ing.uniroma2.it} \\ 
	Department of Civil Engineering and Computer Science\\
	University of Rome Tor Vergata\\
	via del Politecnico 1, 00133 Rome, Italy}

\begin{document}

\maketitle

\begin{abstract}
As inelastic structures are ubiquitous in ma\-ny engineering fields, a central task in computational mechanics is to develop accurate, robust and efficient tools for their analysis. Motivated by the poor performances exhibited by standard displacement-based finite element formulations, attention is here focused on the use of mixed methods as approximation technique, within the small strain framework, for the mechanical problem of inelastic bidimensional structures. Despite a great flexibility characterizes mixed element formulations, several theoretical and numerical aspects have to be carefully taken into account in the design of a high-performance element. The present work aims at providing the basis for methodological analysis and comparison in such aspects, within the unified mathematical setting supplied by generalized standard material model and with special interest towards elastoplastic media. A critical review of the state-of-the-art computational methods is delivered in regard to variational formulations, selection of interpolation spaces, numerical solution strategies and numerical stability. Though those arguments are interrelated, a topic-oriented presentation is resorted to, for the very rich available literature to be properly examined. Finally, the performances of several significant mixed finite element formulations are investigated in numerical simulations. \\[3ex]
\keywords{Mixed finite element; Material nonlinearity; Plasticity; Element state determination; Numerical stability}
\end{abstract}

\section{Introduction}
Design and safety assessment procedures in many engineering fields, such as structural or geotechnical engineering, and industrial applications e.g.~in robotics, automotive, aerospace and biomedical industries, require to predict the behavior of inelastic structures. Thus, a central task in computational mechanics is to conceive computational tools for their analysis, capable to combine accuracy, robustness and efficiency together with appropriate computational demand.

A numerical method for computational inelasticity is hinged on: (i) a \emph{state update algorithm} to solve the constitutive equations of the material under investigation, and (ii) an \emph{approximation technique} of the structural problem capable of capturing distinctive features of inelastic behavior, such as strain and stress concentrations. The first research developments in the former aspect trace back to the mid-nineteen-fifties, with the classical expositions of the infinitesimal theory of plasticity by Hill \cite{Hill_1950} and Koiter \cite{Koiter_1960} and the radial return method for numerical integration of the $J_2$-plasticity evolution equations by Wilkins~\cite{Wilkins_1964}. Since then, computational plasticity has been receiving continuous attention and a comprehensive review of the subject can be found in~\cite{Scalet_Auricchio_ARCME_2017}. Here, attention is focused on structural approximation techniques, among which a prominent role is taken by the finite element method (FEM). The origin of FEM dates back to the mid-nineteen-fifties as well, its start being possibly recognized in the work by Turner, Clough, Martin and Topp~\cite{Turner_Topp_JAS_1956}. Other major contributors during that pioneering age were Melosh, who showed the link between compatible formulations and the minimum potential energy principle~\cite{Melosh_AIAA_1963}, Irons, who proposed the isoparametric formulation \cite{Irons_AIAA_1966, Zienkiewicz_Scott_1969}, and Wilson, who developed the first FEM computer program~\cite{Wilson_PhD_1963} (for an extended historical background on the compatible FEM, the reader is referred to~\cite{Felippa_2011}). 

The displacement-based or compatible FEM starts from the assumption of a kinematic model, exactly satisfying the compatibility condition, to prescribe the equilibrium condition in a weak form. Originally proposed for linear elastic problems, its natural extension to materially nonlinear context just requires the enforcement of the nonlinear constitutive equations, instead of the linear ones, at quadrature points. Unfortunately, compatible formulations for inelastic structures are unable to meet essential computational needs as accuracy, robustness and limited computational cost (for instance, see \cite{Taylor_Auricchio_CM_2003}). In fact, starting from generally adopted polynomial interpolations of the displacement field, the compatible strain (i.e.~determined by differentiation) turns out to be a very poor approximation of the actual strain spatial distribution, possibly highly nonlinear due to inelastic effects. As a consequence, inaccurate stresses are derived from the inelastic constitutive relationships, irrespective of weak equilibrium enforcement. The overall quality of the solution can only be recovered adopting computationally expensive fine meshes, thus reducing the efficiency of the approach. 

Alternative to compatible FEM, in the mid-nineteen-sixties Pian proposed another type of finite elements, the derivation of the stiffness matrix being obtained by an assumed stress distribution and a minimum complementary principle~\cite{Pian_AIAA_1964}. The term hybrid element was coined for such a formulation (for a discussion on the evolution of hybrid elements, see~\cite{Pian_FEAD_1995, Pian_IJNME_2000}). As prescribing the interpolation of more than one independent field, the Pian hybrid element can be regarded as a first instance of mixed finite element formulations. During the nine-seventies, the latter were diffusely investigated as a promising strategy to tackle problems characterized by some physical constraint, such as the incompressible or nearly incompressible behavior of rubber-like media, or the shear constraint in plate problems (e.g., see \cite{Bathe_1996}). In the words of Felippa~\cite{Felippa_2011}, the main advantages of mixed formulations are: (i) \emph{relaxed continuity requirements}, achieved by approximating as primary unknowns not only the displacement, but also the strain and/or the stress fields, (ii) \emph{better stress solution}, computed without differentiation of the displacement field, and (iii) \emph{better displacement solution}, also for coarse and irregular meshes. When considering inelastic structures, mixed methods offer the possibility to overcome the complications arising with displacement-based formulations. Indeed, the introduction of independent stress and/or strain interpolations gives the flexibility to overcome the strict constraint imposed on the displacement approximation by the compatibility condition. Hence, a high-performance finite element formulation that combines accuracy (even for coarse meshes) together with robustness and limited computational cost is potentially conceivable.

The present work aims at supplying an insight into the application of mixed finite element methods to the analysis of inelastic structures, which can provide guidance both to the development of novel formulations and to the use of existing ones. In addition to the review purpose, emphasis is devoted to the construction of a unified mathematical setting simplifying a methodological analysis of the subject and allowing for a direct comparison. To this end, the material nonlinearity is restricted to the class of generalized standard materials in the sense of Halphen and Nguyen \cite{Halphen_Nguyen_JMec_1975}, with special consideration for hardening associative plastic behavior. In particular, the material results to be equipped with an (incremental) energetic structure~\cite{Miehe_Lambrecht_JMPS_2002, Petryk_CRM_2003, Mielke_JMA_2004, Mosler_CMAME_2010, Nodargi_Bisegna_FSI_2014, Nodargi_Bisegna_IJNME_2015}, which is exploited for a synthetical exposition. Though many of the results here presented are also valid in the threedimensional case and might be extended to the finite strain framework, the discussion is referred to bidimensional structures under small strain assumption. In the following, theoretical and numerical aspects to be carefully taken into account in the design of mixed finite element are examined. Specifically:
\begin{itemize}
	\item {an overview of the most investigated \emph{mixed variational formulations} is provided, along with their discretized version;}
	\item {popular choices of the \emph{interpolation spaces} for the unknown fields are reviewed, with the customary assumption of element-supported interpolation of unknown fields additional to displacements;}
	\item {classical and recently addressed \emph{numerical solution strategies} of the resulting discrete structural problem are analyzed, focusing on the properties of the existing algorithms and highlighting relevant strengths and weaknesses;}
	\item {the \emph{numerical stability}, i.e.~the fulfillment of the ellipticity requirement and the inf--sup condition of Babu\v{s}ka and Brezzi \cite{Bathe_1996, Boffi_Brezzi_Fortin_2013}, is explored for the formulations under consideration.}	
\end{itemize}
Though in the design of a mixed finite element all the aforementioned aspects are interrelated (for instance, selection of interpolation spaces and numerical stability or variational formulation and numerical solution strategy), the present exposition is intended to be topic-oriented, in such a way that the very rich available literature can be properly examined. Finally, the performances of several significant mixed finite element formulations are investigated in numerical simulations.

The paper is organized as follows. In Section~\ref{s:constitutive_model}, the constitutive model of generalized standard material is introduced. Mixed variational formulations and relevant discretizations are addressed in Section~\ref{s:gsm_variational}. Section~\ref{s:interpolation} is devoted to the discussion on the selection of the interpolation spaces for the unknowns. Numerical solution strategies of the mixed discrete problem are presented in Section~\ref{s:numerical_strategy}. In Section~\ref{s:stability}, the numerical stability is investigated. Numerical simulations aiming at a comparison of several mixed finite elements available in the literature are presented in Section~\ref{s:numerical_simulations}. Conclusions and perspectives are outlined in Section~\ref{s:conclusions}.

\section{Constitutive model}\label{s:constitutive_model}
In this section, some results about the theory of generalized standard materials are presented, with the aim to derive a sufficiently general framework for covering a significant class of inelastic materials, also including hardening associative elastoplasticity. 

The model of generalized standard material, first introduced in the seminal work of Halphen and Nguyen \cite{Halphen_Nguyen_JMec_1975}, is based on the two fundamental concepts of internal variables and dissipation potential. The internal variables are additional variables, with respect to strain and temperature, needed to completely define the state of the system and to describe dissipative mechanisms typical of the material under investigation \cite{Germain_Suquet_JAM_1983}. On the other hand, the dissipation potential is a convex scalar function which governs the evolution of the internal variables according to the assumption of normal dissipation~\cite{Halphen_Nguyen_JMec_1975}. Initially proposed as a generalization of classical plasticity theory (in that context the normal dissipation hypothesis leads to associative flow law, expressed by the principle of maximum dissipation~\cite{Eve_Rockafellar_QAM_1990, Reddy_Martin_CMAME_1991, Han_DayaReddy_1999}), the model of generalized standard material has been then successfully applied and extended to complex material behaviors (e.g., see \cite{Michel_Suquet_IJSS_2003} and references therein). In particular, it has been shown that, in a time-discrete setting, the evolution of a generalized standard material follows the minimizing path of an incremental energy, given by the sum of free energy and of dissipation potential, with respect to the internal variables \cite{Miehe_Lambrecht_JMPS_2002, Petryk_CRM_2003, Mielke_JMA_2004, Mosler_CMAME_2010, Nodargi_Bisegna_FSI_2014, Nodargi_Bisegna_IJNME_2015}. The introduction of the material incremental energy allows a unified treatment of the variational formulations to be discussed in Section~\ref{s:gsm_variational}.

\subsection{Generalized standard materials}  \label{ss:evolution_equations}
In the framework of small deformation theory and under the assumption of isothermal conditions, the constitutive behavior of a generalized standard material is defined in terms of the infinitesimal strain tensor $\stn$ and of a generalized vector of internal variables $\intvarstn$ \cite{Halphen_Nguyen_JMec_1975}. The energy storage in the deforming material is described by a strictly convex Helmholtz free energy $\helmholtz$, depending on the strain $\stn$ and on the internal variables $\intvarstn$. Accordingly, by standard thermodynamic arguments \cite{Armero_ECM_2004}, the constitutive equations for the stress $\sts$ and the generalized vector of internal forces $\intvarsts$ conjugated to internal variables $\intvarstn$ turn out to be:
\begin{align}
	\begin{aligned}
		&\sts\!\rbr{\stn, \intvarstn} = \nablaarg{\stn}{\helmholtz\!\rbr{\stn, \intvarstn}}, \\[1ex]
		&\intvarsts\!\rbr{\stn, \intvarstn} = -\nablaarg{\intvarstn}{\helmholtz\!\rbr{\stn, \intvarstn}},
	\end{aligned}
\label{eq:const_law}
\end{align}
in which $\nabla$ denotes gradient operator. The flux in time $\dot{\intvarstn}$ of the internal variables~$\intvarstn$ is related to the internal forces $\intvarsts$ by a dissipation potential $\diss$, assumed to be a convex function of the flux $\dot{\intvarstn}$ of the internal variables:
\begin{equation}
	\intvarsts \in \pdiffarg{\dot{\intvarstn}} \diss(\dot{\intvarstn}),
\label{eq:evolution_law}	
\end{equation}
where $\partial$ denotes the subdifferential operator of convex functions~\cite{Rockafeller_1970}. By comparing equations~\eqref{eq:const_law} and~\eqref{eq:evolution_law}, the following differential equation, usually referred to as Biot's equation of standard dissipative system (for instance, see \cite{Miehe_Lambrecht_JMPS_2002} and references therein), is derived:
\begin{equation}
	\nabla_{\!\intvarstn}{\helmholtz\!\rbr{\stn, \intvarstn}} + \pdiffarg{\dot{\intvarstn}}{\diss(\dot{\intvarstn})} \ni \bzero.
\label{eq:gsm_ev}
\end{equation}

Equations~\eqref{eq:const_law} and \eqref{eq:gsm_ev}, complemented with initial conditions on the internal variables~$\intvarstn$, completely determine the material evolution for given Helmholtz free energy $\helmholtz$ and dissipation potential $\diss$. In particular, for rate-independent material behavior, the latter is positively homogeneous of degree one (e.g., see~\cite{Armero_ECM_2004}). Consequently it has a cone-like graph and is not differentiable at the origin. It is remarked that the dissipation potential $\diss$ is customarily introduced as the support function of the elastic domain $\cK$ of the forces $\intvarsts$ (for instance, see~\cite{Simo_Hughes_1998}): 
\begin{equation}
	\diss(\dot{\intvarstn}) = \sup_{\intvarsts \in \cK} \{\intvarsts \cdot \dot{\intvarstn}\},
\label{eq:diss}
\end{equation} 
and its subdifferential set at the origin coincides with the elastic domain $\cK$ itself.

\subsection{Incremental energy minimization}\label{ss:time_discretization}
Apart from classical evolution equations for which closed-form solutions are available (e.~g., see \cite{Simo_Hughes_1998}), Biot's equation \eqref{eq:gsm_ev} has generally to be solved by numerical integration. Yet many approaches have been deeply investigated to date (for a comprehensive review, see~\cite{Scalet_Auricchio_ARCME_2017}), a common approximation strategy consists in a backward Euler integration scheme, which is particularly attractive because of its stability also for large time steps (for instance, see \cite{Borja_2013}). Accordingly, the internal variables $\intvarstn$ are assumed to vary linearly in each time step $\sbr{t\n, t\npo}$ of a suitable partition of the time span of interest~$\sbr{0, T}$. The result is an incremental form of Biot's equation, governing the response at time $t\npo$ for given state at time $t\n$:
\begin{equation}
	\partial_{\Delta\intvarstn} \helmholtz\!\rbr{\stn, \intvarstn\n + \Delta\intvarstn}  + 
	\partial_{\Delta\intvarstn} \diss\!\rbr{\Delta\intvarstn}  \ni 
	\bzero,
\label{eq:Biot_equation_discrete}
\end{equation}
where subscript~$n$ [resp., $n+1$] denotes evaluation at time $t = t\n$ [resp., $t = t\npo$] and $\Delta$ stands for the increment within~$\sbr{t\n, t\npo}$.

According to equation~\eqref{eq:Biot_equation_discrete}, the material evolution is proven to take place along the minimizing path of the convex function
\begin{equation}
	\iemin\!\rbr{\stn, \Delta\intvarstn} = \helmholtz\!\rbr{\stn, \intvarstn\n + \Delta\intvarstn} + \diss\!\rbr{\Delta\intvarstn}
\label{eq:incremental_energy}	
\end{equation}
with respect to the increment of internal variables $\Delta\intvarstn$ \cite{Miehe_Lambrecht_JMPS_2002}. This argument motivates the definition of the material incremental energy as:
\begin{equation}
	\ie\!\rbr{\stn} = \inf_{\Delta\intvarstn} \cbr{\iemin\!\rbr{\stn, \Delta\intvarstn}}.	
\label{eq:incremental_energy_min}	
\end{equation}
Using equation \eqref{eq:Biot_equation_discrete}, it is also a simple matter to check that infimum \eqref{eq:incremental_energy_min} is attained at $\Delta\intvarstn = \bzero$ if and only if the trial stress, i.e. the stress computed assuming no internal variable evolution, belongs to the subdifferential set of the dissipation potential at the origin. Excluding that trivial case, a direct computation shows that, in any point of differentiability of the dissipation potential, the incremental energy $\ie$ plays the role of stress incremental potential:
\begin{align}
	\begin{aligned}
		\diff{\ie}{\beps}\!\rbr{\stn}
			&= \Big\{\nablaarg{\stn}{\iemin}\!\rbr{\stn, \Delta\intvarstn} \\ 
			&\hspace{0.5cm} + \nablaarg{\stn}{\Delta\intvarstn}\!\rbr{\stn} \nablaarg{\Delta\intvarstn}{\iemin}\!\rbr{\stn, \Delta\intvarstn}\Big\}\Big|_{{\nablaarg{\Delta\intvarstn}{\iemin} = \bzero}} \\[1ex]
			&= \fixed{\nablaarg{\stn}{\helmholtz}\!\rbr{\stn, \Delta\intvarstn}}
			{{\nablaarg{\Delta\intvarstn}{\iemin} = \bzero}} \\[1ex]
			&= \sts\!{}.
	\end{aligned}
\label{eq:stress_potential}
\end{align}
It is remarked that $\text{d}\ie / \text{d}\stn$ denotes the total derivative of the incremental energy $\ie$ with respect to the strain~$\stn$, i.e. it also takes into account the contribution stemming from the variation of the increment of internal variables $\Delta\intvarstn$ due to the variation of the strain $\stn$.

\begin{remark}\label{r:complementary_ie}
Beside the incremental energy, it will be useful to introduce its partial Legendre conjugate with respect to the strain~$\stn$, i.e.~the incremental complementary energy:
\begin{equation}
	\ie\conj\!\rbr{\sts} = \sup_{\stn}\cbr{\sts \cdot \stn - \ie\!\rbr{\stn}}.
\label{eq:incremental_compl_energy}
\end{equation}
By exploiting the definition of incremental energy via equations~\eqref{eq:incremental_energy} and~\eqref{eq:incremental_energy_min} and by exchanging $\inf$ and $\sup$ operators, it is a simple matter to check that: 
\begin{equation}
	\ie\conj\!\rbr{\sts} = - \inf_{\Delta\intvarstn}\cbr{-\helmholtz\conj\!\rbr{\sts, \intvarstn\n + \Delta\intvarstn} + \diss\!\rbr{\Delta\intvarstn}},
\label{eq:incremental_compl_energy_2}
\end{equation}
in which $\helmholtz\conj$ is the Gibbs free energy, i.e.~the partial Legendre transformation of the Helmholtz free energy $\helmholtz$ with respect to the strain $\stn$:
\begin{equation}
	\gibbs\!\rbr{\sts, \intvarstn\n + \Delta\intvarstn} = \sup_{\stn}\cbr{\sts \cdot \stn - \helmholtz\!\rbr{\stn, \intvarstn\n + \Delta\intvarstn}}.
\label{eq:Gibbs}
\end{equation}
In the context of generalized standard materials, the following classical result of convex analysis~\cite{Rockafeller_1970}:
\begin{equation}
	\sts = \pdiffarg{\stn}{\ie\!\rbr{\stn}} \iff \stn = \pdiffarg{\sts}{\ie\conj\!\rbr{\sts}},
\label{eq:characterization}
\end{equation}
ensures the incremental complementary energy~$\ie\conj$ to define an inverse constitutive law, which maps the stress $\sts$ into the strain~$\stn$ (the increment of the generalized vector of internal variables~$\Delta\intvarstn$ being derived as the minimum point in equation~\eqref{eq:incremental_compl_energy_2}). From a computational point of view, the inversion of the constitutive law can be regarded as a problem equivalent to the material state update, provided the Helmholtz free energy $\helmholtz$ is replaced by the opposite of the Gibbs free energy $\gibbs$ (the latter can be computed analytically via equation~\eqref{eq:Gibbs} in case of linear elastic constitutive law). In passing, it is observed that plane stress condition is naturally imposed within this setting. Conversely, for considering plane strain condition, it is sufficient to introduce the constraint of vanishing out-of-plane strain components in the computation of Gibbs free energy $\helmholtz\conj$ from equation~\eqref{eq:Gibbs}.
\end{remark}

\subsection{Associative elastoplastic materials}\label{s:elastoplastic_materials}
Originally proposed to be a generalization of the classical plasticity theory, the model of generalized standard material naturally covers hardening elastoplastic behavior. In this case, the internal variables~$\intvarstn$ are required to capture the two phenomena of plastic yielding and material hardening. In regard to plastic yielding, an additive decomposition of the infinitesimal strain $\stn$ into elastic and plastic parts, $\stn\e$ and $\stn\p$ respectively, is introduced \cite{Armero_ECM_2004}:
\begin{equation}
	\stn = \stn\e + \stn\p.
\label{eq:add_decomp}
\end{equation}
Accordingly, the plastic strain $\stn\p$ is assumed as an internal variable describing the yielding mechanism, whereas equation~\eqref{eq:add_decomp} is thought of as definition of the elastic strain $\stn\e$. On the other hand, material hardening is modeled selecting as internal variables strain-like kinematic and isotropic hardening variables, $\balpha\kin$ and $\alpha\iso$, respectively. Correspondingly to such choices of internal variables $\intvarstn$, the conjugated forces $\intvarsts$ are the stress $\sts$ and the stress-like kinematic [resp. isotropic] hardening variable $\bq\kin$ [resp.~$q\iso$]. The Helmholtz free energy is usually expressed by~\cite{Armero_ECM_2004}:
\begin{equation}
	\helmholtz = \helmholtz\e + \helmholtz\iso + \helmholtz\kin,
\label{eq:elastic_energy_bis}	
\end{equation}
where $\helmholtz\e$ is the elastic energy and $\helmholtz\kin$ [resp.~$\helmholtz\iso$] is the kinematic [resp.~isotropic] hardening potential. In case of linear elastic response with linear hardening, such potentials result in:
\begin{align}
	\begin{aligned}
		&\helmholtz\e\!\rbr{\stn, \stn\p\n+\Delta\stn\p} = \frac{1}{2} \rbr{\stn - \stn\p\n - \Delta\stn\p} \cdot \fC\e \rbr{\stn - \stn\p\n - \Delta\stn\p}, \\[0.5ex]
		&\helmholtz\kin\!\rbr{\alpha\kin\n+\Delta\alpha\kin} = \frac{1}{2} k\kin \norm{\alpha\kin\n+\Delta\alpha\kin}^2, \\[0.5ex]
		&\helmholtz\iso\!\rbr{\alpha\iso\n+\Delta\alpha\iso} = \frac{1}{2} k\iso \rbr{\alpha\iso\n+\Delta\alpha\iso}^2,
	\end{aligned}
\label{eq:elastic_energy}
\end{align}
where $\fC\e$ is the elastic tensor, $k\kin$ [resp.~$k\iso$] is the kinematic [resp.~isotropic] hardening modulus and $\norm{\bullet}$ denotes euclidean vector norm. Consequently, the stress~$\sts$ and the forces~$\intvarsts$ work-conjugated to the internal variables~$\intvarstn$ are derived by means of the constitutive law~\eqref{eq:const_law}:
\begin{align}
	\begin{aligned}
		&\sts = \fC\e \rbr{\stn - \stn\p\n - \Delta\stn\p}, \\[1ex]
		&\bq\kin = - k\kin \rbr{\balpha\kin\n + \Delta\balpha\kin}, \\[1ex]
		&q\iso = - k\iso \rbr{\alpha\iso\n + \Delta\alpha\iso}.
	\end{aligned}
\end{align}
In passing, it is observed that the stress $\sts$ depends on the elastic strain~$\stn\e$. 

In the context of elastoplastic materials, the stress~$\sts$ is constrained to belong to an admissible set $\cK$, the elastic domain, defined as the zero-sublevel set of a convex yield function~$\yf$:
\begin{equation}
	\cK := \{ \bsigma : \yf (\sts, \bq\kin, q\iso)  \leq 0 \} \,.
\end{equation}
In particular, the effect of stress-like kinematic [resp. isotropic] hardening variable $\bq\kin$ [resp.~$q\iso$] is to determine a shift [resp.~a scaling] of the yield domain $\cK$. Within this setting, the definition of dissipation potential \eqref{eq:diss} is regarded as the statement of the principle of maximum dissipation and results in the normality rule between plastic strain rate and yield surface during the plastic flow (associative behavior) \cite{Eve_Rockafellar_QAM_1990, Reddy_Martin_CMAME_1991, Han_DayaReddy_1999}.

\begin{remark}
Considering for simplicity the case of perfect plasticity, the dissipation potential introduced in equation~\eqref{eq:diss} turns out to be:
\begin{equation}
	\diss\!\rbr{\Delta\stn\p} = \sup_{\sts \in \cK} \cbr{\sts \!\cdot\! \Delta\stn\p} = \sup_{\sts} \cbr{\sts \!\cdot\! \Delta\stn\p \!-\! I_{\cK}\!\rbr{\sts}}\!,
\end{equation}
where~$I_{\cK}$ denotes the indicator function of the elastic domain~$\cK$:
\begin{equation}
	I_{\cK}\!\rbr{\sts} = \begin{dcases}
		0, \quad &\sts \in \cK, \\
		+\infty, \quad &\text{otherwise}
	\end{dcases}.
\label{eq:indicator_fun}
\end{equation}
By solving the sup problem with respect to the stress~$\sts$, the following variational inclusion is derived:
\begin{equation}
	\bzero \in \Delta\stn\p - N_{\cK}\!\rbr{\bsigma},
\label{eq:normality_law}
\end{equation}
with $N_{\cK}\!\rbr{\sts} = \partial_{\sts}I_{\cK}\!\rbr{\sts}$ as the normal cone of the yield locus~$\cK$ at the stress point~$\bsigma$. As $\cK$ is the zero-sublevel set of the convex yield function $\yf$, the normal cone~$N_{\cK}\!\rbr{\sts}$ can be represented by (e.g., see~\cite{Han_DayaReddy_1999}):
\begin{equation}
	N_{\cK}\!\rbr{\bsigma} =
	\begin{dcases}
		\varnothing, \quad &\yf\!\rbr{\bsigma} > 0, \\[1ex]
		\cbr{\Delta\lambda \, \pdiffarg{\bsigma}{\yf\!\rbr{\sts}}, \,\, \Delta\lambda \geq 0}, \quad &\yf\!\rbr{\bsigma} = 0, \\[1ex] 
		\bzero, \quad &\yf\!\rbr{\bsigma} < 0,
	\end{dcases}
\end{equation}
or, equivalently:
\begin{multline}
	N_{\cK}\!\rbr{\bsigma} = \cbr{ \Delta\lambda \, \pdiffarg{\bsigma}{\yf}\!\rbr{\sts} \,\bdot\, \right. \\[0.5ex] \left. \Delta\lambda \geq 0, \,\, \yf\!\rbr{\bsigma} \leq 0, \,\,\Delta\lambda \, \yf\!\rbr{\bsigma} = 0 }.
\label{eq:normal_cone}
\end{multline}
Accordingly, the variational inclusion~\eqref{eq:normality_law} is rephrased in:
\begin{align}
	\begin{aligned}
		&\Delta\stn\p = \Delta\lambda \, \pdiffarg{\bsigma}{\yf}\!\rbr{\sts}, \\[0.5ex] 
		&\Delta\lambda \geq 0, \,\, \yf\!\rbr{\bsigma} \leq 0, \,\,\Delta\lambda \, \yf\!\rbr{\bsigma} = 0,
	\end{aligned}
\end{align}
responding to the customary format of the normality law between increment of plastic strain and yield surface, complemented with plastic admissibility and consistency conditions. Analogous conclusions hold when hardening is considered. For a comprehensive treatment of the subject, for example see~\cite{Han_DayaReddy_1999}.

\end{remark}

\section{Mixed variational formulations}\label{s:gsm_variational}
Within the context of inelastic structures constituted of generalized standard materials, mixed variational formulations can be systematically derived from a generalization of the classical Hu--Washizu (HW) functional \cite{Washizu_1982}, involving as independent variables the fluxes of internal variables in addition to displacement, strain and stress fields.

In the time-discrete framework introduced in Section~\ref{ss:time_discretization}, let $\dmn \in \mathbb{R}^{N\txtsub{dim}}$, with $N\txtsub{dim} = 2$ as the space dimension, be the domain occupied by the body at time~$t\n$. In the small deformation regime, the \emph{generalized HW functional}, governing its mechanical evolution within the time step $\sbr{t\n, t\npo}$, is introduced by:
\begin{align}
	\begin{aligned}
		&\gf : V \times E \times I \times S \rightarrow \mathbb{R}, \\[0.5ex]
		&\gf = \intbr{\cbr{\iemin\!\rbr{\stn, \intvarstn\n + \Delta\intvarstn} - \sts \cdot \rbr{\stn - \sym\nabla{\displ}}}}.
\label{eq:gf}
	\end{aligned}
\end{align}
Here $\displ \in V$, $\stn \in E$ and $\sts \in S$ respectively denote current displacement, strain and stress fields, $\intvarstn\n \in I$ is the vector of internal variables at time~$t\n$ and $\Delta\intvarstn \in I$ is the increment of the vector of internal variables in the finite time step $t\npo-t\n$. Natural choices of the functional spaces $V$, $E$, $I$ and $S$ consist in:
\begin{align}
	\begin{aligned}
		&V = \sbr{H^1_{\bdy_{\displ}}\!\rbr{\dmn}}^{N\txtsub{dim}}, \\[1ex]
		&E = S = \sbr{L^2\!\rbr{\dmn}}^{N\txtsub{dim}\times N\txtsub{dim}}\txtsub{Sym}, \\[1ex]
		&I = \sbr{L^2\!\rbr{\dmn}}^{N_{I}},
	\end{aligned}
\end{align}
in which $L^2$ is the space of square integrable functions, $H^1$ is the space of square integrable functions along with their first weak derivatives, $H^1_{\bdy_{\displ}}$ is the subspace of $H^1$-functions vanishing (in the trace sense) over the constrained part~$\bdy_{\displ}$ of the boundary~$\bdy$, and $N_{I}$ is the number of required internal variables. It is recalled, equation~\eqref{eq:incremental_energy}, that $\iemin$ is a convex function given by the sum of the Helmholtz free energy $\helmholtz$ and the dissipation potential $\diss$. Hereafter the stronger hypothesis of strict convexity of $\iemin$, i.e.~the presence of some hardening mechanism in the material model, will be assumed. Moreover $\sym\nabla$ stands for the symmetric gradient operator.

Upon noticing that the first variation of $\gf$ reads:
\begin{align}
	\begin{aligned}
		\var\gf = &-\intbr{\div{\sts} \cdot \var\displ \,} + \intbdybr{\sts \bn \cdot \var\displ \,} \\[1ex]
		&+\intbr{\cbr{\nablaarg{\stn}{\iemin} \! \rbr{\stn, \intvarstn\n + \Delta\intvarstn} - \sts} \cdot \var\stn \,} \\[1ex]
		&+\intbr{\cbr{\nablaarg{\Delta\intvarstn}{\iemin}\!\rbr{\stn, \intvarstn\n + \Delta\intvarstn}} \cdot \var\Delta\intvarstn \,},	\\[1ex]
		&-\intbr{\rbr{\stn - \sym\nabla{\displ}} \cdot \var\sts \,}
	\end{aligned}
\end{align}
its stationary conditions with respect to $\displ$, $\stn$, $\Delta\intvarstn$ and $\sts$, can be interpreted as equilibrium,
constitutive, material evolution and compatibility conditions. Hence the solution of the mechanical problem at hand turns out to be the (unique) stationary point of the functional~$\gf$. Similar functionals have been investigated for elastoplasticity in~\cite{Comi_Perego_CMAME_1995, Capsoni_Corradi_CMAME_1997}, by considering the internal variables (i.e.~plastic strain and strain-like hardening variables) as primary unknowns to be independently interpolated.

Several variational statements can be derived from the generalized HW functional by a priori enforcing its stationary condition with respect to some of the fields it depends on  (e.g., see \cite{Bathe_CS_2001}). Though from a theoretical perspective any possible formulation is interesting \emph{per se}, some general considerations are in order. For the formulation to be directly implemented in a standard displacement-driven finite element computer program, it is convenient to have the displacement field as an independent variable. In addition, because of its essential role in engineering applications and in order to avoid accuracy losses due to differentiation of interpolated quantities, an explicit interpolation of the stress field is desirable as well. 

One main choice concerns the possibility of assuming the strain or the flux of internal variables as independent quantities, with the underlying idea that interpolating both would be cumbersome and computational costly. In fact, a reasonable claim is that capturing accurately the strain field spatial distribution, with its possible highly nonlinear distribution, makes unnecessary the internal variable field interpolation. That leads to the HW principle (Section~\ref{ss:HW}), which can be derived by minimization of the generalized HW functional with respect to the fluxes of internal variables, and hence takes displacement, stress and strain fields as primary variables. Examples of its application in elastoplasticity can be found in~\cite{Weissman_Jamjian_IJNME_1993, Comi_Perego_CMAME_1995, Capsoni_Corradi_CMAME_1997, Kasper_Taylor_C&S_2000, Moharrami_Cocchetti_CS_2015, Nodargi_Bisegna_C&S_2017}. In the same line of having a rich and flexible strain approximation, a modification of the HW principle, in which the strain field is represented as the sum of a compatible part and of an enhanced part, yields the enhanced strain (ES) formulations (Section~\ref{ss:enhanced_strain}), explored in elastoplasticity in~\cite{Simo_Hughes_JAM_1986, Simo_Rifai_IJNME_1990, Piltner_Taylor_IJNME_1995, Piltner_Taylor_IJNME_1999, Piltner_CM_2000}. 

As an alternative point of view, it can be argued that possible highly nonlinear spatial distributions of the strain field would only stem from similar patterns of the flux of internal variables. Accordingly, proper interpolation of the latter might be sufficient for the design of high-performance finite elements. That results in the return map functional (Section~\ref{ss:return_map}), obtained by minimization of the generalized HW functional with respect to the strain field, and depending on displacement, stress and flux of internal variables. A remarkable feature of such a formulation is the weak imposition of the material evolution equations at element level (instead of quadrature-point level), thus recalling the notion of \emph{constitutive law of the finite element} tracing back to Maier's work \cite{Maier_Meccanica_1968, Maier_Meccanica_1969}. For elastoplasticity, a deeply explored strategy amounts at the adoption of the complementary mixed (CM) functional, initially proposed in~\cite{Simo_Taylor_CMAME_1989}, and after addressed in~\cite{Pinsky_CMAME_1987, Bilotta_Casciaro_CMAME_2007, Krabbenhoft_Wriggers_IJNME_2007, Mendes_Castro_FEAD_2009, Bilotta_Garcea_CS_2012, Bilotta_Garcea_FEAD_2011} (for an overview on complementarity problems in structural engineering, the reader is referred to~\cite{Bolzon_ARCME_2017}). That can be regarded as a reduced version of the return map functional, in which only the plastic multiplier field is interpolated among all internal variables, thus implying a weak element flow law.   

The most simple structure of a mixed formulation responds to the Hallinger--Reissner (HR) principle (Section~\ref{ss:HR}), which follows from taking the stationarity of the generalized HW functional with respect to both strain field and fluxes of internal variables, thus depending on displacement and stress fields only. Based on the clue that an accurate description of the stress field overcomes the need of an explicit approximation of strain and flux of internal variables, several recent applications have dealt with elastoplasticity~\cite{Schroder_Miehe_CMAME_1997, Contrafatto_Ventura_IJNME_2004, Leonetti_Aristodemo_FEAD_2015, Bilotta_Leonetti_FEAD_2016, Schroder_Starke_CMAME_2017}.

In the following, a detailed discussion on such formulations is presented. Methodological analysis and comparisons are fostered by the  unified mathematical setting resulting from the framework of generalized standard material structures. Nevertheless, special attention is devoted to the particular case of hardening associative elastoplastic media. 

\begin{remark}\label{r:hybrid}
Some slight modifications of the generalized HW functional in equation~\eqref{eq:gf} might be considered. A noteworthy example concerns the restriction of such functional to the space of self-equilibrated stresses~$S\txtsub{div} = \cbr{\sts \in S \bdot \, \div\sts = \bzero}$. Upon integrating by parts, a generalized hybrid HW functional is derived (e.g., see~\cite{Pian_AIAA_1964, Madeo_Casciaro_FEAD_2012, Nodargi_Bisegna_C&S_2017}):
\begin{align}
	\begin{aligned}
		&\gf\txtsub{hyb} : V \times E \times I \times S\txtsub{div} \rightarrow \mathbb{R}, \\[0.5ex]
		&\gf\txtsub{hyb} = 
		\intbr{\cbr{\iemin\!\rbr{\stn, \intvarstn\n + \Delta\intvarstn} - \sts \cdot \stn}} +
		\intbdybr{\sts\bn \cdot \displ},
	\end{aligned}
\end{align}
involving the displacement field only over the domain boundary~$\bdy$. Accordingly, in a discrete setting, functional~$\gf\txtsub{hyb}$ requires interpolating the displacement field just over the element boundary. On the other hand, a self-equilibrated stress approximation has to be generated (in Section~\ref{sss:Airy}, a systematic approach resorting to Airy's function approach is reviewed). By the same technique here presented, a hybrid version can be derived for all the functionals to be discussed.
\end{remark}

\subsection{Hu--Washizu functional} \label{ss:HW}
The \emph{classical HW functional} is obtained by taking the stationarity of the functional~$\gf$ with respect to the increment of internal variable~$\Delta\intvarstn$:
\begin{align}
	\begin{aligned}
		&\hw : V \times E \times S \rightarrow \mathbb{R}, \\[0.5ex]
		&\hw =
		\intbr{\cbr{\ie \! \rbr{\stn} - \sts \cdot \rbr{\stn - \sym\nabla{\displ}} }},
	\end{aligned}
\label{eq:HW}
\end{align}
where the definition of the incremental energy~\eqref{eq:incremental_energy_min} has been exploited. To derive a finite element approximation of problem~\eqref{eq:HW} on a given mesh~$\dmn = \cup_{e=1}^{N\elem} \dmn\elem$, a sequence of interpolation spaces, parametrized by a mesh parameter $h \rightarrow 0$, is introduced for the unknown fields:
\begin{align}
	\begin{aligned}
		&V\interp = \cbr{\displ\interp \in V \bdot \displ\interp = \displmtx \displprm,\, \displprm \in \mathbb{R}^{N_{\displ}}}, \\[1ex]
		&E\interp = \cbr{\stn\interp \in E \bdot \fixed{\stn\interp}{\dmn\elem} = \stnmtx\elem \stnprm\elem,\, \stnprm\elem \in \mathbb{R}^{N_{\stn}}}, \\[1ex]
		&S\interp = \cbr{\sts\interp \in S \bdot \fixed{\sts\interp}{\dmn\elem} = \stsmtx\elem \stsprm\elem,\, \stsprm\elem \in \mathbb{R}^{N_{\sts}}},
	\end{aligned}
\label{eq:HW_interp_spaces}
\end{align}
where $\displprm$ collects nodal degrees of freedom (DOFs), $\stnprm\elem$ and $\stsprm\elem$ are the interpolation parameters of strain and stress at element level, respectively, and $\displmtx$, $\stnmtx\elem$, $\stsmtx\elem$ are the relevant interpolation matrices (customary Voigt notation is adopted for interpolated strains~$\sts\interp$ and strains $\stn\interp$). Hence, the following mixed finite element functional is derived:
\begin{align}
	\begin{aligned}
		&\hw\interp : V\interp \times E\interp \times S\interp \rightarrow \mathbb{R}, \\[0.5ex]
		&\hw\interp = 
			\sum_{e=1}^{N\elem} \cbr{
			\intelembr{\ie \! \rbr{\stnmtx\elem\stnprm\elem}} - 
			\stsprm\elem\txtsup{T}\stsstnmtx\elem\txtsup{T}\stnprm\elem +
			\stsprm\elem\txtsup{T}\compmtx\elem\displprm\elem
			},
	\end{aligned}
\label{eq:HW_interp_fun}
\end{align}
in which $\displprm\elem$ gathers the nodal DOFs of the element~$e$, $\graddisplmtx$ is the discrete symmetric gradient operator, i.e.~such that $\sym\nabla\displ\interp = \graddisplmtx\displprm$, and $\graddisplmtx\elem$ collects its columns pertaining to the element~$e$, $\compmtx\elem$ is the element compatibility matrix and~$\stsstnmtx\elem$ is the element stress-strain operator. In particular, the latter are defined by:
\begin{equation}
	\compmtx\elem = \intelembr{\stsmtx\elem\txtsup{T}\graddisplmtx\elem}, \quad
	\stsstnmtx\elem = \intelembr{\stnmtx\elem\txtsup{T}\stsmtx\elem}.
\label{eq:comp_stsstn_operator}
\end{equation}
The finite element solution is then obtained by taking the stationarity of the functional~$\hw\interp$. To this end, because the interpolation of strain and stress is inter-element discontinuous, the stationarity with respect to $\stnprm\elem$ and~$\stsprm\elem$ can be imposed at element level and reads:
\begin{align}
	\begin{aligned}
		\bzero &=\intbr{\stnmtx\elem\txtsup{T}\partial_{\stn}\ie \! \rbr{\stnmtx\elem\stnprm\elem}} - \stsstnmtx\elem\stsprm\elem, \\[1ex]
		\bzero &= -\stsstnmtx\elem\txtsup{T}\stnprm\elem + \compmtx\elem\displprm\elem.
	\end{aligned}
\label{eq:HW_fun_stat}
\end{align}
whereas the vector of internal nodal forces results to be:
\begin{equation}
	\force\internal\elem = \compmtx\elem\txtsup{T}\stsprm\elem,
\label{eq:HW_q_int}
\end{equation}
to be assembled at structural level for imposing global equilibrium as in standard displacement-based finite elements.

\subsection{Enhanced strain functional} \label{ss:enhanced_strain}
The \emph{ES functional} can be derived from the HW functional~\eqref{eq:HW} by considering strain fields of the form \cite{Simo_Rifai_IJNME_1990}:
\begin{equation}
	\stn = \sym\nabla\displ + \tilde\stn,
\label{eq:enhanced_stn_def}
\end{equation}
where $\tilde\stn$ is referred to as the enhanced part of the strain field. Accordingly, the ES functional results to be:
\begin{align}
	\begin{aligned}
		&\tilde{\cE} : V \times \tilde E \times S \rightarrow \mathbb{R}, \\[0.5ex]
		&\tilde{\cE} = 
		\intbr{ \cbr{ \ie\!\rbr{\sym\nabla\displ + \tilde\stn} - \sts \cdot \tilde\stn } },
	\end{aligned}
\label{eq:enhanced_stn}
\end{align}
with~$\tilde E = \sbr{L^2\!\rbr{\dmn}}^{N\txtsub{dim} \times N\txtsub{dim}}\txtsub{Sym}$ as the ambient space of the enhanced part of the strain field. Moving on to the discrete formulation, let the interpolation spaces~$V\interp$ and~$S\interp$ be defined as in equations~\eqref{eq:HW_interp_spaces}$\txtsub{1,3}$ and let the interpolated enhanced strain space~$\tilde E\interp$ be given by:
\begin{equation}
	\tilde{E}\interp = \cbr{\tilde\stn\interp \in \tilde{E} \bdot \fixed{\tilde\stn\interp}{\dmn\elem} = \tilde\stnmtx\elem \tilde\stnprm\elem,\, \tilde\stnprm\elem \in \mathbb{R}^{N_{\tilde\stn}}},
\label{eq:assumed_strain_interp_space}		
\end{equation}
where $\tilde\stnprm\elem$ is the vector collecting the element enhanced strain interpolation parameters and~$\tilde\stnmtx\elem$ is the relevant element interpolation matrix. Accordingly, the discrete ES functional turns out to be:
\begin{align}
	\begin{aligned}
		&\tilde{\cE}\interp : V\interp \times \tilde E\interp \times S\interp \rightarrow \mathbb{R}, \\[0.5ex]
		&\tilde{\cE}\interp  = 
		\sum_{e=1}^{N\elem} \cbr{ \intelembr{
		\ie(\graddisplmtx\elem\displprm\elem  + \tilde\stnmtx\elem \tilde\stnprm\elem) } -
		\stsprm\elem\txtsup{T}\tilde\stsstnmtx\elem\txtsup{T}\stnprm\elem },
	\end{aligned}
\label{eq:enhanced_fun}
\end{align}
in which the element stress-enhanced-strain field operator has been introduced by:
\begin{equation}
	\tilde\stsstnmtx\elem = \intelembr{\tilde\stnmtx\elem\txtsup{T}\stsmtx\elem}.
\label{eq:tilde_bG}
\end{equation}
Proceeding as in the previous section, because enhanced strain and stress interpolations are discontinuous across element boundaries, the stationarity conditions of~$\tilde \cE\interp$ with respect to~$\tilde\stnprm\elem$ and~$\stsprm\elem$ can be imposed at element level, yielding:
\begin{align}
	\begin{aligned}
		\bzero &= \intelembr{\tilde\stnmtx\elem\txtsup{T} \partial_{\stn}\ie(\graddisplmtx\elem\displprm\elem + \tilde\stnmtx\elem \tilde\stnprm\elem)} -
		\tilde\stsstnmtx\elem\stsprm\elem, \\[1ex]
		\bzero &= \tilde\stsstnmtx\elem\txtsup{T} \tilde\stnprm\elem,
	\end{aligned}
\label{eq:enhanced_fun_stat}
\end{align}
with the vector of element nodal internal forces:
\begin{equation}
	\force\internal\elem = \intelembr{\graddisplmtx\elem\txtsup{T} \partial_{\stn}\ie(\graddisplmtx\elem\displprm\elem + \tilde\stnmtx\elem \tilde\stnprm\elem)},
\label{eq:enhanced_q_int}
\end{equation}
being evaluated in the same fashion of standard compatible formulations, provided the enhanced part of the strain is accounted for. 

It is observed that the compatibility condition~\eqref{eq:enhanced_fun_stat}$\txtsub{2}$, which weakly constraints the enhanced part of the strain field to be vanishing, prescribes a compatibility requirement on strain and stress approximations. Indeed, two limit situations can be considered. In the original approach proposed in~\cite{Simo_Rifai_IJNME_1990}, $L^{2}$-orthogonal interpolations of enhanced strain and stress fields are assumed. Hence, equation~\eqref{eq:enhanced_fun_stat}$\txtsub{2}$ is identically satisfied and the second term in equation~\eqref{eq:enhanced_fun_stat}$\txtsub{1}$ vanishes. Accordingly, the stress parameters~$\stsprm\elem$ are effectively eliminated from the finite element formulation and the derivation of a stress recovery method becomes a central issue. Oppositely, in case the matrix~$\tilde\stsstnmtx\elem\txtsup{T}$ results to be full-column rank, equation~\eqref{eq:enhanced_fun_stat}$\txtsub{2}$ implies that $\tilde\stnprm\elem = \bzero$ and the present ES formulation boils down to a standard displacement-based formulation.

\begin{remark}\label{r:assumed_stn}
As noted in~\cite{Simo_Rifai_IJNME_1990}, the ES formulation includes, as a particular case, the method of incompatible modes originally proposed in~\cite{Wilson_Ghaboussi_1973}. In fact, the latter is characterized by the following enhanced strain interpolation space:
\begin{equation}
	\tilde{E}\interp = \cbr{\tilde\stn\interp \in \tilde{E} \bdot \fixed{\tilde\stn\interp}{\dmn\elem} = \graddisplmtx\elem\txtsup{i}\displprm\elem\txtsup{i},\, \displprm\elem\txtsup{i} \in \mathbb{R}^{N_{\tilde\stn}}},
\end{equation}
where the columns of the matrix~$\graddisplmtx\elem\txtsup{i}$ constitute a system of element-supported incompatible modes and the parameters $\displprm\elem\txtsup{i}$ denote the relevant element-level DOFs. Such description of the enhanced strain field will be exploited~in Section~\ref{ss:enh_stn_interp} for the derivation of a suitable interpolation of the enhanced strain field.
\end{remark}

\begin{remark}
It is noteworthy to compare ES formulations with assumed strain methods, also referred to as B-bar methods. Assumed strain methods can be introduced as particular instances of the HW formulation, provided the strain interpolation space is selected to be~\cite{Simo_Hughes_JAM_1986}:
\begin{equation}
	E\interp = \cbr{\stn\interp \in E \bdot \stn\interp = \bar\graddisplmtx\displprm,\,\,\, \displprm \in \mathbb{R}^{N_{\displ}}}	
\end{equation} 
with $\bar\graddisplmtx$ as assumed strain interpolation matrix and~$\displprm$ the nodal DOFs. As a consequence, the interpolated assumed strain functional, derived from equation~\eqref{eq:HW}, reads:
\begin{align}
	\begin{aligned}
		&\bar{\cE}\interp \bdot V\interp \times S\interp \rightarrow \mathbb{R}, \\[0.5ex]
		&\bar{\cE}\interp = 
		\sum_{e=1}^{N\elem} \cbr{\intelembr{
			\ie\!\rbr{\bar\graddisplmtx\elem\displprm\elem} } -
			\stsprm\elem\txtsup{T}\!\rbr{\bar\compmtx\elem - \compmtx\elem}\!\displprm\elem },
	\end{aligned}
\label{eq:assumed_strain_fun}
\end{align}
where the element compatibility matrix~$\compmtx\elem$ is defined as in equation~\eqref{eq:comp_stsstn_operator}$\txtsub{1}$ and its B-bar counterpart, i.e.~the element assumed-strain compatibility matrix, is introduced by:
\begin{equation}
	\bar\compmtx\elem = \intelembr{\stsmtx\elem\txtsup{T}\bar\graddisplmtx\elem}.
\end{equation}
The stationarity condition with respect to the stress parameters~$\stsprm\elem$ is enforced at element level and yields:
\begin{equation}
	\bzero = \rbr{\bar\compmtx\elem - \compmtx\elem} \displprm\elem, 
\label{eq:assumed_strain_stat}
\end{equation}
whence the vector of element nodal internal forces is given by:
\begin{equation}
	\force\internal\elem = \intelembr{\bar\graddisplmtx\elem\txtsup{T} \partial_{\stn}\ie\!\rbr{\bar\graddisplmtx\elem\displprm\elem}} - 
		\rbr{\bar\compmtx\elem - \compmtx\elem}\txtsup{T}\!\stsprm\elem.
\label{eq:assumed_strain_q_int}
\end{equation}
Similarly to the ES formulation, the compatibility condition~\eqref{eq:assumed_strain_stat} can be interpreted as a constraint on the stress interpolation space for a given selection of the assumed strain space. In particular, the strictest constraint possible, i.e.~prescribing that $\bar\compmtx\elem - \compmtx\elem = \bzero$, yields a variational consistent version of the B-bar procedure, as originally introduced in the context of nearly-incompressible media in~\cite{Hughes_IJNME_1980}. In such a case the stress parameters~$\stsprm\elem$ are eliminated from the finite element formulation and a stress recovery method is needed~\cite{Simo_Hughes_JAM_1986}. In Section~\ref{sss:HW_identical_stsstn} it will be shown that a variationally consistent stress recovery strategy is naturally furnished by the interpretation of B-bar methods with eliminated stresses as instances of a HW formulation with identical stress and strain interpolations.

In particular, for linear material problems, the ES formulation with eliminated stresses can be interpreted as a B-bar method under the assumption~\cite{Simo_Rifai_IJNME_1990}:
\begin{equation}
	\bar\graddisplmtx = \graddisplmtx - \tilde\stnmtx \rbr{\intbr{\tilde\stnmtx\txtsup{T}\fC\txtsup{e}\tilde\stnmtx}}^{-1} \!\rbr{\intbr{\tilde\stnmtx\txtsup{T}\fC\txtsup{e}\graddisplmtx}},
\end{equation}
where it is recalled that $\fC\txtsup{e}$ denotes the linear elastic material stiffness tensor.
\end{remark}
 
\begin{remark}
An alternative version of the ES formulation, proposed in~\cite{Piltner_Taylor_IJNME_1995, Piltner_Taylor_IJNME_1999, Piltner_CM_2000}, is based on the following modification of the HW functional:
\begin{align}
	\begin{aligned}
		&\tilde{\cE} \bdot V \times E \times \tilde E \times S \rightarrow \mathbb{R}, \\[0.5ex]
		&\tilde{\cE} = 
			\intbr{ \cbr{ \ie\!\rbr{\stn} - \sts \cdot \rbr{\stn - \sym\nabla\displ - \tilde\stn} } },
	\end{aligned}
\label{eq:assumed_strain_fun_alt}
\end{align}
and differs from the one in equation~\eqref{eq:enhanced_stn} for the weak enforcement of condition~\eqref{eq:enhanced_stn_def}. \\
Upon introducing interpolation spaces as in equations \eqref{eq:HW_interp_spaces} and~\eqref{eq:assumed_strain_interp_space}, the discrete version of functional~$\tilde{\cE}$ follows as:
\begin{align}
	\begin{aligned}
		&\tilde{\cE} \bdot V\interp \times E\interp \times \tilde E\interp \times S\interp \rightarrow \mathbb{R}, \\[0.5ex]
		&\tilde{\cE}\interp = 
		\sum_{e=1}^{N\elem} \cbr{
			\intelembr{\ie\!\rbr{\stnmtx\elem\stnprm\elem}} - \stsprm\elem\txtsup{T}\, (
			\stsstnmtx\elem\txtsup{T}\stnprm\elem - 			
			\compmtx\elem\displprm\elem -
			\tilde\stsstnmtx\elem\txtsup{T}\tilde\stnprm\elem)}.
	\end{aligned}
\end{align}
Its stationarity conditions with respect to~$\stnprm\elem$, $\stsprm\elem$ and~$\tilde\stnprm\elem$ is enforced at element level and read:
\begin{align}
	\begin{aligned}
		\bzero &= \intelembr{\stnmtx\elem\txtsup{T}\partial_{\stn}\ie\!\rbr{\stnmtx\elem\stnprm\elem}} - \stsstnmtx\elem\stsprm\elem, \\[0.5ex]
		\bzero &= \stsstnmtx\elem\txtsup{T}\stnprm\elem - \compmtx\elem\displprm\elem - \tilde\stsstnmtx\elem\txtsup{T}\tilde\stnprm\elem, \\[1ex]
		\bzero &= \tilde\stsstnmtx\elem\stsprm\elem,
	\end{aligned}
\label{eq:enhanced_fun_stat_alt}	
\end{align}
the vector of element nodal internal forces resulting in:
\begin{equation}
	\force\internal\elem = \compmtx\elem\txtsup{T}\stsprm\elem.
\label{eq:enhanced_fun_q_int_alt}	
\end{equation}
As for the ES formulation~\eqref{eq:enhanced_stn}, equation \eqref{eq:enhanced_fun_stat_alt}$\txtsub{3}$ turns out to be a constraint on the stress interpolation space for a prescribed enhanced strain space. However, an $L^2$-orthogonal selection would now cancel out the operator~$\tilde\stsstnmtx\elem$ from the discrete formulation, thus reducing such approach to classical HW functional discussed in Section~\ref{ss:HW}.
\end{remark}

\begin{remark}
For a direct comparison of the two ES formulations in equations~\eqref{eq:assumed_strain_fun} and~\eqref{eq:assumed_strain_fun_alt}, it is required eliminating the strain field~$\stn$ from~$\tilde{\cE}$.  Accordingly, a priori enforcing the infimum of~$\tilde{\cE}$ with respect to~$\stn$, the following functional is obtained:
\begin{align}
	\begin{aligned}
		&\tilde\hpr \bdot V \times \tilde E \times S  \rightarrow \mathbb{R}, \\[0.5ex]
		&\tilde\hpr = 
			 \intbr{ \cbr{ -\ie\conj\!\rbr{\sts} + \sts \cdot \rbr{\sym\nabla\displ + \tilde\stn} } },
	\end{aligned}
\end{align}
with $\ie\conj$ as the complementary incremental energy, defined in equation~\eqref{eq:incremental_compl_energy}. Accordingly, $\tilde\hpr$ can be regarded as an enhanced strain version of the HR functional to be discussed in Section~\ref{ss:HR}.
\end{remark}

\subsection{Return map functional} \label{ss:return_map}
In~\cite{Simo_Taylor_CMAME_1989}, Simo and coauthors proposed a CM finite element formulation for elastoplasticity, based on an explicit interpolation of the plastic multiplier in addition to displacements and stresses. Such formulation is here regarded as a reduced version of a slightly more general one, derived from the so-called \emph{return map functional}, in which the overall plastic strain is assumed as primary variable. In particular, by taking the infimum of the generalized HW functional~$\gf$ with respect to the total strain $\beps$, the return map functional results in:
\begin{align}
	\begin{aligned}
		&\cR \bdot V \times S \times I \rightarrow \mathbb{R}, \\[0.5ex]
		&\cR =
		\intbr{\cbr{-\gibbs \! \rbr{\sts, \intvarstn\n + \Delta\intvarstn} + \diss \! \rbr{\Delta\intvarstn} + \bsigma \cdot \sym\nabla{\bu} }},
	\end{aligned}
\end{align}
where the definitions of function~$\iemin$, equation~\eqref{eq:incremental_energy}, and of Gibbs free energy, equation~\eqref{eq:Gibbs}, have been used. 

Aiming at the derivation of the \emph{CM functional} explored in~\cite{Simo_Taylor_CMAME_1989}, perfect plasticity is considered for simplicity. Accordingly, the only internal variable is the plastic strain and the functional~$\cR$ reduces to:
\begin{align}
	\begin{aligned}
		&\cR =
		\intbr{\Big\{- \gibbs\!\rbr{\bsigma} + \diss \! \rbr{\Delta\beps\p}  \\[0.5ex]
		&\hspace{2cm} + \bsigma \cdot \rbr{\sym\nabla{\bu} - \beps\p\n - \Delta\beps\p} \Big\}}.
	\end{aligned}
\label{eq:rm}	
\end{align}
Let the increment of plastic strain be expressed in the form~$\Delta\beps\p = \Delta\lambda \bn$, where $\Delta\lambda \geq 0$ is the plastic multiplier and $\bn$, such that $\norm{\bn} = 1$, is the plastic flow direction. The return map functional can be then represented by:
\begin{align}
	\begin{aligned}
		& \hat{\cR} \bdot V \times S \times L \times I \rightarrow \mathbb{R}, \\[1ex]
		&\hat{\cR}= 
		\intbr{\Big\{- \gibbs\!\rbr{\bsigma} +\diss \! \rbr{\Delta\lambda \bn} \\[0.5ex] 
			& \hspace{2cm}+ \bsigma \cdot \rbr{\sym\nabla{\bu} - \beps\p\n - \Delta\lambda \bn} \Big\}\, }, \\[1ex]
		&\text{subject to } \Delta\lambda \geq 0, \quad \norm{\bn} = 1,
	\end{aligned}
\end{align}
where~$L = L^2\!\rbr{\dmn}$ is the ambient space of the plastic multiplier field. Because of the first-degree positive homogeneity of the dissipation potential, the unitary constraint on the plastic flow direction $\bn$ amounts at a rescaling of the plastic multiplier $\Delta\lambda$. Hence, a free minimization of $\hat{\cR}$ with respect to $\bn$ can be enforced, yielding the CM functional:
\begin{align}
	\begin{aligned}
		&\check{\cR} \bdot V \times S \times L \rightarrow \mathbb{R}, \\[1ex]
		&\check{\cR}=
		\intbr{ \cbr{- \gibbs\!\rbr{\bsigma} + \bsigma \cdot \rbr{\sym\nabla{\bu}-\beps\p\n} - \Delta\lambda I_{\cK}\!\rbr{\bsigma}}}, \\[1ex]
		&\text{subject to } \Delta\lambda \geq 0.
	\end{aligned}
\end{align}
It is recalled that~$I_{\cK}$ denotes the indicator function of the elastic domain~$\cK$, equation~\eqref{eq:indicator_fun}, thus implying the pointwise imposition of the plastic admissibility condition $\yf\!\rbr{\sts}\leq 0$. The functional~$\check{\cR}$ can be reformulated in terms of the yield function~$\yf$, instead of the elastic domain~$\cK$. To this end, the stationary condition of $\check{\cR}$ with respect to the stress $\bsigma$ is considered:
\begin{multline}
	\bzero \in \delta_{\bsigma}{\check{\cR}} = \intbr{ \Big\{- \pdiffarg{\bsigma}\gibbs\!\rbr{\bsigma} + \rbr{\sym\nabla{\bu}-\beps\p\n} \\[0.5ex] 
	-\Delta\lambda N_{\cK}\!\rbr{\bsigma}\Big\} \cdot \var\bsigma \,},
\label{eq:comp_fun_stat_sts}
\end{multline}
resulting in a variational inclusion with~$N_{\cK}$ as the normal cone of the yield locus~$\cK$. On the basis of the normal cone representation in equation~\eqref{eq:normal_cone}, the CM functional~$\check{\cR}$ can be equivalently stated by~\cite{Simo_Taylor_CMAME_1989}:
\begin{align}
	\begin{aligned}
		&\check{\cR} = 
		\intbr{\cbr{- \gibbs\!\rbr{\bsigma} + \bsigma \cdot \rbr{\sym\nabla{\bu} - \beps\p\n } - \Delta\lambda \yf\!\rbr{\bsigma} }}, \\[1ex]
		&\text{subject to} \quad \Delta\lambda \geq 0, \quad
		\yf\!\rbr{\sts} \leq 0, \quad
		\Delta\lambda \, \yf\!\rbr{\sts} = 0.
	\end{aligned}
\end{align}

For the discretization of functional~$\check{\cR}$, displacement and stress interpolation spaces, $V\interp$ and~$S\interp$ respectively, are given as in equations~\eqref{eq:HW_interp_spaces}$\txtsub{1,3}$. Moreover, the following plastic multiplier interpolation space is introduced:
\begin{equation}
	L\interp = \cbr{\pmlt\interp \in L \bdot \fixed{\pmlt\interp}{\dmn\elem} = \pmltmtx\elem \pmltprm\elem,\,\,\, \pmltprm\elem \in \mathbb{R}^{N_{\pmlt}}}, \\[1ex]
\end{equation}
where~$\pmltprm\elem$ is the vector collecting element plastic multiplier interpolation parameters and~$\pmltmtx\elem$ is the matrix of the relevant interpolation functions, assumed to satisfy~$\pmltmtx\elem \geq \bzero$ in~$\dmn\elem$. Assuming the Gibbs free energy of a linear elastic material, the discrete CM functional results to be:
\begin{align}
	\begin{aligned}
		&\check{\cR}\interp \bdot V\interp \times S\interp \times L\interp \rightarrow \mathbb{R}, \\[1ex]
		&\check{\cR}\interp = 
		\sum_{e=1}^{N\elem} \cbr{-
			\frac{1}{2} \stsprm\elem\txtsup{T} \cmplmtx\elem \stsprm\elem +  
			\stsprm\elem\txtsup{T} \sbr{ \compmtx\elem\displprm\elem - \rbr{\stnprm\p\n}\elem } - 
			\pmltprm\elem\txtsup{T} \yfvct\!\rbr{\stsprm\elem}  }, \\[1ex]
		&\text{subject to } \pmltprm\elem \geq \bzero	
	\end{aligned}
\end{align}
where $\compmtx\elem$ is the element compatibility matrix introduced in equation~\eqref{eq:comp_stsstn_operator}$\txtsub{1}$, $\cmplmtx\elem\e$ denotes the elastic compliance matrix:
\begin{equation}
	\cmplmtx\elem\e = \intelembr{ \stsmtx\elem\txtsup{T} \rbr{\fC\e}^{-1} \stsmtx\elem },
\label{eq:return_map_H}
\end{equation}
and the following positions hold:
\begin{align}
	\begin{aligned}
		&\yfvct\!\rbr{\stsprm\elem} = \intelembr{ \pmltmtx\elem\txtsup{T} \yf\!\rbr{\stsmtx\elem \stsprm\elem} }, \\[0.5ex]
		&\rbr{\stnprm\p\n}\elem = \intelembr{\stsmtx\elem\txtsup{T}\beps\p\n},
	\end{aligned}
\label{eq:return_map_position}
\end{align}
respectively with the meaning of element yield function and element plastic strain at previous time step. Exploiting the element-supported interpolations of stress and plastic multiplier fields, the stationarity condition with respect to~$\stsprm\elem$ and the Kuhn-Tucker conditions with respect to~$\pmltprm\elem$ are imposed at element level:
\begin{align}
	\begin{aligned}
		& \bzero = \compmtx\elem\displprm\elem - \rbr{\stnprm\p\n}\elem - \flowmtx\!\rbr{\stsprm\elem}\pmltprm\elem - \cmplmtx\elem\e\stsprm\elem, \\[1ex]
		& \pmltmtx\elem\pmltprm\elem \geq \bzero, \quad
		\yfvct\!\rbr{\stsprm\elem} \leq \bzero, \quad
		\pmltprm\elem\txtsup{T} \yfvct\!\rbr{\stsprm\elem} = 0,
	\end{aligned}
\label{eq:return_map_stat}
\end{align}
with the element plastic flow operator defined as the derivative of the element yield function~$\yfvct$ with respect to the stress parameters~$\stsprm\elem$:
\begin{equation}
	\flowmtx\!\rbr{\stsprm\elem} = \pdiff{\yfvct}{{\stsprm\elem}}\!\rbr{\stsprm\elem} = \intelembr{\stsmtx\elem\txtsup{T} \partial_{\sts}\yf\!\rbr{\stsmtx\elem \stsprm\elem} \pmltmtx\elem }.
\end{equation}
The vector of element nodal internal forces to be assembled at structural level for the imposition of structural equilibrium follows as the stationarity condition of~$\check\cR\interp$ with respect to the element DOFs~$\displprm\elem$: 
\begin{equation}
	\force\internal\elem = \compmtx\elem\txtsup{T}\stsprm\elem.
\label{eq:return_map_q_int}
\end{equation}
Distinctive feature of such a formulation, which has been deeply investigated also in recent contributions~\cite{Bilotta_Casciaro_CMAME_2007, Bilotta_Garcea_FEAD_2011, Bilotta_Garcea_CS_2012}, is that the flow law is weakly imposed at element level, equation~\eqref{eq:return_map_stat}$\txtsub{1}$. On the other hand, plastic admissibility and consistency conditions, equations~\eqref{eq:return_map_stat}$\txtsub{2}$, can be either weakly or strongly imposed, depending on the selected interpolation of the plastic multiplier field, as it will be discussed in Section~\ref{ss:plastic_mlt}. Deep similarities can be here recognized with the problem of non-smooth multi-surface plasticity at material point level (for a review, see~\cite{Karaoulanis_ARCME_2013}), upon interpreting each component of the element yield function~$\yfvct$ as parametrizing one surface of the multi-surface yield locus. In Section~\ref{ss:element_return_mapping}, such point of view will be exploited for the exposition of an element return mapping procedure resembling the integration scheme customarily adopted for the elastoplastic material state update.

\subsection{Hallinger--Reissner functional}\label{ss:HR}
The \emph{classical HR functional} can be derived by imposing the stationarity of the HW functional~$\hw$ with respect to the strain $\stn$ (or equivalently, by imposing the stationarity of the return map functional $\cR$ with respect to the increment of internal variable~$\Delta\intvarstn$):
\begin{align}
	\begin{aligned}
		&\hpr \bdot V \times S \rightarrow \mathbb{R}, \\[0.5ex]
		&\hpr = 
		\intbr{ \cbr{- \ie\conj\!\rbr{\sts} + \sts \cdot \sym\nabla{\displ}}},
	\end{aligned}
\label{eq:hpr}
\end{align}
where $\ie\conj$ is the incremental complementary energy defined in equation~\eqref{eq:incremental_compl_energy}. 

Introducing interpolation spaces~$V\interp$ and~$S\interp$ as in equations~\eqref{eq:HW_interp_spaces}$\txtsub{1,3}$, the discrete HR functional becomes:
\begin{align}
	\begin{aligned}
		&\hpr\interp \bdot V\interp \times S\interp \rightarrow \mathbb{R}, \\[0.5ex]
		&\hpr\interp = 
		\sum_{e=1}^{N\elem}\cbr{ \intelembr{ - \ie\conj\!\rbr{\stsmtx\elem\stsprm\elem} } + 
			\stsprm\elem\txtsup{T}\compmtx\elem\displprm\elem },
	\end{aligned}
\label{eq:hpr_discrete}
\end{align}
with $\compmtx\elem$ as the element compatibility matrix in equation~\eqref{eq:comp_stsstn_operator}$\txtsub{1}$. The stationary condition with respect to the stress parameters~$\stsprm\elem$ is imposed at element level, yielding:
\begin{equation}
	\bzero = -\intelembr{\stsmtx\elem\txtsup{T}\partial_{\sts}\ie\conj\!\rbr{\stsmtx\elem\stsprm\elem}} + \compmtx\elem\displprm\elem.
\label{eq:HR_stat_2}
\end{equation}
On the other hand, the stationarity with respect to the element nodal DOFs~$\displprm\elem$ gives the definition of the vector of element nodal internal forces:
\begin{equation}
	\force\internal\elem = \compmtx\elem\txtsup{T} \stsprm\elem,
\label{eq:HR_q_int}
\end{equation}
to be assembled at structural level in order to impose global equilibrium equations. 

It is observed that the present HR formulation differs from the CM functional discussed in Section~\ref{ss:return_map} because of the local imposition of the flow law underlying the definition of the incremental complementary energy~$\ie\conj$.

\begin{remark}\label{r:HR_Schroder_Starke}
A recent contribution tackling the elastoplastic structural problem by means of the HR functional is~\cite{Schroder_Starke_CMAME_2017}. Therein an equivalent reformulation of the compatibility condition~\eqref{eq:HR_stat_2} is presented. In fact, recalling the inverse constitutive law~\eqref{eq:characterization}, the following relationship holds true:
\begin{equation}
	\partial_{\sts}\ie\conj\!\rbr{\stsmtx\elem\stsprm\elem} = \stn =  \rbr{\fC\e}^{-1}\stsmtx\elem\stsprm\elem + \stn\p\n + \Delta\stn\p\!\rbr{\stsmtx\elem\stsprm\elem},
\end{equation}
where the increment of plastic strain~$\Delta\stn\p$ is intended as a known function of the stress~$\sts$, i.e.~obtained by solving the sup problem in the definition of the incremental complementary energy~$\ie\conj$, equation~\eqref{eq:incremental_compl_energy}. It is remarked that such function is well-defined only in case of hardening behavior (for example, see~\cite{Nodargi_Bisegna_C&S_2017}). Accordingly, the compatibility condition~\eqref{eq:HR_stat_2} can be rewritten as:
 \begin{equation}
	\bzero = \compmtx\elem\displprm\elem - \rbr{\stnprm\p\n}\elem - \Delta\stnprm\elem\p\!\rbr{\stsprm\elem} - \cmplmtx\elem\e\stsprm\elem,
\label{eq:HR_stat}
\end{equation}
with the discrete increment of plastic strain~$\Delta\stnprm\elem\p$ depending on the stress parameters~$\stsprm\elem$ by:
\begin{equation}
	\Delta\stnprm\elem\p\!\rbr{\stsprm\elem} = \intelembr{\stsmtx\elem\txtsup{T} \Delta\stn\p\!\rbr{\stsmtx\elem\stsprm\elem}},
\label{eq:HR_plstn_prm}
\end{equation}
whereas the elastic compliance matrix $\cmplmtx\elem\e$ and the discrete plastic strain at previous time step $\rbr{\stnprm\n\p}\elem$  have been introduced in equations~\eqref{eq:return_map_H} and~\eqref{eq:return_map_position}$\txtsub{2}$, respectively.
\end{remark}

\section{Interpolation of unknown fields}\label{s:interpolation}
The variational formulations investigated in Section~\ref{s:gsm_variational} require the selection of suitable interpolation spaces for the relevant unknown fields. Here attention is restricted to the design of quadrilateral elements, which usually exhibit superior performances than triangular elements~\cite{Bathe_1996, Zienkiewicz_2013} (a mathematical treatment of the issues raised by finite element approximation on quadrilateral meshes can be found in~\cite{Arnold_Gastaldi_CNME_2001, Arnold_Falk_MC_2002}). 

The starting point for a displacement interpolation is of course represented by standard Lagrangian approximation. Nevertheless, with the aim of conceiving a high-performance formulation, several potential improvements have been sought for in the literature. Among those, the introduction of corner rotational DOFs, initially explored by Allman \cite{Allman_IJNME_1988}, has been recognized as an attractive strategy to pursue intermediate accuracy between linear and quadratic elements with translations only~\cite{Yunus_Cook_IJNME_1989, Aminpour_IJNME_1992, Ibrahimbegovic_FEAD_1990,  Pimpinelli_FEAD_2004, Choo_Lee_FEAD_2006, Choi_Lee_CMAME_2006, RezaieePajand_Karkon_EJMS_2013, Rebiai_Belounarb_Measurement_2014}. On the other hand, drilling DOFs are associated to zero-energy spurious modes that need to be suppressed, via a modification of underlying variational formulation \cite{Hughes_Brezzi_CMAME_1989, Ibrahimbegovic_Wilson_IJNME_1990, Nodargi_Bisegna_IJNME_2016} or an ad-hoc treatment of the relevant interpolation \cite{Madeo_Casciaro_FEAD_2012}. Alternatively to the introduction of drilling rotations, improved performances can be achieved by using biquadratic or serendipity quadrilaterals, which in turn require a larger number of DOFs~\cite{Bergmann_Mukherjee_IJNME_1990, Darilmaz_Kumbasar_CS_2006, Cen_Zhou_CMAME_2011, Madeo_Zucco_FEAD_2014, Nodargi_Bisegna_C&S_2017}.

The first example of independent stress interpolation traces back to the assumed stress rectangular element proposed by Pian~\cite{Pian_AIAA_1964}. In \cite{Pian_Sumihara_IJNME_1984}, Pian and Sumihara introduced a filter technique of the stress field approximation via incompatible displacement modes, which can be exploited for deriving stress interpolation spaces of arbitrary dimension. As a result, the popular $5$-$\beta$ stress interpolation was proposed. Its applications, or slight modifications of its, are adopted in~\cite{Simo_Taylor_CMAME_1989, Simo_Rifai_IJNME_1990, Bilotta_Garcea_FEAD_2011, Bilotta_Garcea_CS_2012, Krabbenhoft_Wriggers_IJNME_2007, Xie_Zhou_CNME_2008}. In the line of hybrid element design, several authors have also considered the assumption of a self-equilibrated stress approximation, which allows to interpolate the displacement field just over the element boundary to impose boundary equilibrium \cite{Pian_AIAA_1964, Xie_Zhou_CNME_2008, Cen_Zhou_CMAME_2011, Madeo_Casciaro_FEAD_2012, RezaieePajand_Karkon_EJMS_2013, Madeo_Zucco_FEAD_2014, Nodargi_Bisegna_C&S_2017}. 

In HW formulations, emphasis has been placed on the central role of strain interpolation to reproduce highly nonlinear effects, such as strain concentration, that are typical of inelastic structures. In this regard, two main strategies have been explored: to equally interpolate strain and stress fields or, alternatively, to assume a discontinuous strain interpolation within each finite element. The former approach, suggested by linear elastic applications and followed in~\cite{Taylor_IJNME_2000, Wisniewski_Turska_C&S_2009, Nodargi_Bisegna_IJNME_2016}, has the major advantage to result in a closed form numerical solution procedure at element level (Section~\ref{sss:HW_identical_stsstn}). By contrast, though requiring an iterative element state determination scheme, the latter approach appears more promising in accurately capturing strain nonlinear distributions stemming from material nonlinearity, as shown in the context of beam problems (for instance, see \cite{Spacone_Taucer_EESD_1996, Taylor_Auricchio_CM_2003}). The extension to plane problems has been investigated in~\cite{Nodargi_Bisegna_C&S_2017}. 

A strain enhancement is considered in ES formulations for deriving a rich and flexible strain approximation. Following the path originating from the incompatible displacement method introduced in~\cite{Wilson_Ghaboussi_1973}, an enhanced strain interpolation space can be derived from selection of suitable incompatible displacement modes~\cite{Taylor_Wilson_IJNME_1976, Ibrahimbegovic_FEAD_1990, Hueck_Wriggers_CNME_1994, Piltner_Taylor_IJNME_1995, Piltner_Taylor_IJNME_1999, Piltner_CM_2000}. Differently, Simo and Rifai proposed to assume an independent interpolation of the enhanced strain over the parent element and then to exploit the usual push-forward transformation between the parent element and the physical one to derive a frame-invariant approximation~\cite{Simo_Rifai_IJNME_1990}. Actually, the two approaches are related to each other, for the existence of a displacement field corresponding to the assumed enhanced strains~\cite{Simo_Rifai_IJNME_1990}. 

Finally, when dealing with CM formulations, a discretization of the plastic multiplier field is needed. Several choices have been considered in~\cite{Simo_Taylor_CMAME_1989, Bilotta_Garcea_FEAD_2011, Bilotta_Casciaro_CMAME_2007}, each corresponding to a specific weak imposition of the flow law at element level.

In this section, some of the aforementioned discrete spaces are reviewed. Such exposition does not aim to be exhaustive or comprehensive, due to the noteworthy research effort in that direction. On the other hand, it is explicitly remarked that in the design of a mixed finite element, numerical stability requirements represent a strict compatibility constraint on the selection of interpolation of unknown fields (Section~\ref{s:stability}). Because interpolations are defined at element level, whenever no confusion may arise, subscript~$e$ referring to an individual element is henceforth suppressed in order to avoid cumbersome notation.

\subsection{Displacement}
A brief review of the most commonly adopted interpolations for the displacement field is here presented. Let $\pnt_i = \cbr{x_i, y_i}$ denote the coordinates of the typical node $\node_i$ in the element reference frame, and let $\displ_i = \cbr{u_i, v_i}$ denote the corresponding displacement vector. It is assumed that $i = 1, \dots, \nn$, with $\nn$ as the number of element nodes. It is useful to introduce the two vectors collecting nodal coordinates and nodal displacements:
\begin{equation}
	\hat{\pnt} = \cbr{\pnt_1; \dots{}; \pnt_{\nn}}, \quad
	\displprm = \cbr{\displ_1; \dots{}; \displ_{\nn}},
\end{equation}
in which the semicolon symbol stands for column stacking. Following the isoparametric concept in the finite element method and parametrizing the reference element by $\bxi = \rbr{\xi, \eta} \in \sbr{-1,1}\times \sbr{-1,1}$, the physical element geometry and the element displacement field are approximated by the same shape functions: 
\begin{equation}
	\pnt\!\rbr{\bxi} = \displmtx\!\rbr{\bxi}\hat{\pnt}, \quad
	\displ\!\rbr{\bxi} = \displmtx\!\rbr{\bxi}\displprm,
\end{equation}
where the displacement interpolation matrix~$\displmtx$ has the following structure:
\begin{equation}
	\displmtx = \sbr{\,N_1 \bI_2 \mid \dots{} \mid N_{\nn} \bI_2\,},
\end{equation}
with $\bI_2$ denoting the identity matrix of dimension~$2$ and the separator symbol denoting row concatenation. Adopting such notation, bilinear and quadratic serendipity isoparametric shape functions are discussed in the following.

\subsubsection{Bilinear isoparametric interpolation}\label{sss:bilinear}
In the design of a four-node quadrilateral element ($\nn = 4$), bilinear isoparametric interpolation is customarily adopted for the displacement field approximation (for example, see~\cite{Pian_Sumihara_IJNME_1984, Simo_Taylor_CMAME_1989, Simo_Rifai_IJNME_1990, Piltner_Taylor_IJNME_1995, Bilotta_Garcea_FEAD_2011, Bilotta_Garcea_CS_2012, Schroder_Starke_CMAME_2017, Piltner_CM_2000}). The relevant shape functions are given by:
\begin{align}
	\begin{aligned}
		&N_1\!\rbr{\bxi} &= \rbr{1-\xi}\rbr{1-\eta}/4, \\[1ex]
		&N_2\!\rbr{\bxi} &= \rbr{1+\xi}\rbr{1-\eta}/4, \\[1ex]
		&N_3\!\rbr{\bxi} &= \rbr{1+\xi}\rbr{1+\eta}/4, \\[1ex]
		&N_4\!\rbr{\bxi} &= \rbr{1-\xi}\rbr{1+\eta}/4.
	\end{aligned}
\end{align}
For future convenience, the corresponding Jacobian matrix is explicitly derived (for instance, see~\cite{Schroder_Starke_CMAME_2017}):
\begin{equation}
	\bJ\!\rbr{\bxi} = \pdiff{\pnt}{\bxi}\!\rbr{\bxi} = \sbr{\begin{matrix}
		a_1 + a_2 \eta & b_1 + b_2 \eta \\[1ex]
		a_3 + a_2 \xi & b_3 + b_2 \xi
	\end{matrix}},
\end{equation}
where the following positions have been introduced:
\begin{align}
	\begin{aligned}
		a_1 &= \rbr{-x_1+x_2+x_3-x_4}/4, \\[1ex]
		a_2 &= \rbr{x_1-x_2+x_3-x_4}/4, \\[1ex]
		a_3 &= \rbr{-x_1-x_2+x_3+x_4}/4, \\[1ex]
		b_1 &= \rbr{-y_1+y_2+y_3-y_4}/4, \\[1ex]
		b_2 &= \rbr{y_1-y_2+y_3-y_4}/4, \\[1ex]
		b_3 &= \rbr{-y_1-y_2+y_3+y_4}/4.
	\end{aligned}
\label{eq:position_J}
\end{align}

\subsubsection{Serendipity isoparametric interpolation}\label{sss:serendipity}
Serendipity isoparametric interpolation is a valuable choice for achieving quadratic convergence of the formulation, yet avoiding the internal DOFs needed in biquadratic interpolation. Assuming $\nn = 8$ as number of element nodes, and labelling corner nodes as 1--4 and mid-side nodes as 5--8, the corresponding shape functions result to be~\cite{Zienkiewicz_2013}:
\begin{align}
	\begin{aligned}
		N_1\!\rbr{\bxi} &= \rbr{1-\xi}\rbr{1-\eta}\rbr{-1-\xi-\eta}/4, \\[0.5ex]
		N_2\!\rbr{\bxi} &= \rbr{1+\xi}\rbr{1-\eta}\rbr{-1+\xi-\eta}/4, \\[0.5ex]		
		N_3\!\rbr{\bxi} &= \rbr{1+\xi}\rbr{1+\eta}\rbr{-1+\xi+\eta}/4, \\[0.5ex]
		N_4\!\rbr{\bxi} &= \rbr{1-\xi}\rbr{1+\eta}\rbr{-1-\xi+\eta}/4, \\[0.5ex]	
		N_5\!\rbr{\bxi} &= \rbr{1-\xi^2}\rbr{1-\eta}/2, \\[0.5ex]
		N_6\!\rbr{\bxi} &= \rbr{1+\xi}\rbr{1-\eta^2}/2, \\[0.5ex]
		N_7\!\rbr{\bxi} &= \rbr{1-\xi^2}\rbr{1+\eta}/2, \\[0.5ex]
		N_8\!\rbr{\bxi} &= \rbr{1-\xi}\rbr{1-\eta^2}/2.
	\end{aligned}
\label{eq:serendipity}
\end{align}
The restriction of such functions to the element boundary amounts at a quadratic Lagrangian interpolation of the side displacements:
\begin{align}
	\begin{aligned}
		\side_i \rbr{\zeta} &= N_1\!\rbr{\zeta} \bbx_i + N_2\!\rbr{\zeta} \bbx_j + N_3\!\rbr{\zeta} \bbx_k, \quad \zeta \in \sbr{-1,1}, \\[1ex]
		\fixed{\displ}{\side_i}\!\rbr{\zeta} &= N_1\!\rbr{\zeta} \displ_i + N_2\!\rbr{\zeta} \displ_j + N_3\!\rbr{\zeta} \displ_k, \quad \zeta \in \sbr{-1,1},
	\end{aligned}
\label{eq:displ_mtx_1}
\end{align}
where $\side_i$ is the element side joining nodes $\node_i$ and $\node_j$, $\node_k$ denotes the mid-side node and the Lagrangian shape functions are:
\begin{align}
	\begin{aligned}
		&N_1\!\rbr{\zeta} = -\zeta \rbr{1-\zeta}/2, \\[0.5ex]
		&N_2\!\rbr{\zeta} = \zeta \rbr{1+\zeta}/2, \\[0.5ex]
		&N_3\!\rbr{\zeta} =  \rbr{1+\zeta}/2.
	\end{aligned}
\label{eq:Lagrange}
\end{align}
In~\cite{Nodargi_Bisegna_C&S_2017}, the serendipity isoparametric interpolation~\eqref{eq:serendipity} is used as geometric map, whereas the element boundary displacement interpolation~\eqref{eq:Lagrange} is adopted in conjunction with a self-equilibrated stress approximation.

\subsection{Stress}
The quadrilateral proposed by Pian and Sumihara in~\cite{Pian_Sumihara_IJNME_1984} represents a landmark formulation of an assumed stress finite element. In particular, high popularity has been gained by the so-called $5$-$\beta$ stress interpolation and the related filter technique useful to eliminate undesired stress modes. An alternative, systematic strategy to derive a stress interpolation which also guarantees strong fulfillment of internal equilibrium equations, consists in resorting to Airy's function approach ~\cite{RezaieePajand_Karkon_EJMS_2013, Madeo_Zucco_FEAD_2014, Nodargi_Bisegna_C&S_2017}.  In this section, some aspects of the stress approximation are considered. In the following, the ordering $\sts = \cbr{\sigma_{x}; \sigma_{y}; \tau_{xy}}$ is adopted to represent the stress field in Voigt notation.

\subsubsection{Pian-Sumihara stress interpolation}\label{sss:Pian_Sumihara}
Referring to a general quadirilateral, the following complete linear approximation in the natural coordinates~$\bxi$ is assumed as starting point in~\cite{Pian_Sumihara_IJNME_1984}:
\begin{equation}
	\sts\!\rbr{\bxi} = \stsmtx\!\rbr{\bxi}\!\stsprm, \,\, 
		\stsmtx\!\rbr{\bxi} = \sbr{\begin{matrix}  
			1 & \xi & \eta & 0 & 0 & 0 & 0 & 0 & 0 \\[1.5ex]
			0 & 0 & 0 & 1 & \xi & \eta & 0 & 0 & 0  \\[1.5ex]
			0 & 0 & 0 & 0 & 0 & 0 & 1 & \xi & \eta 
		\end{matrix}}\!.
\label{eq:Pian_sts}
\end{equation}
In order to derive the Pian-Sumihara $5$-$\beta$ stress interpolation, the Pian filter technique is used. That exploits auxiliary displacement fields $\displ_{\lambda}$, additional to the compatible displacements characterizing the formulation under consideration, to impose the following further equilibrium constraints:
\begin{equation}
	\bzero = \intbr{\sts \cdot \sym\nabla\displ_{\lambda}\,}.
\end{equation}
Such procedure yields an improvement of the equilibrium conditions, provided that (i) the auxiliary displacement fields are linearly independent from the displacement interpolation, and (ii) the auxiliary strain field has zero mean over the integration domain, in order to not modify the representation of homogeneous stresses~\cite{Pian_Sumihara_IJNME_1984, Bilotta_Casciaro_CMAME_2002}. In particular, by selecting the auxiliary displacement field:
\begin{align}
	\begin{aligned}
		&\displ_{\lambda}\!\rbr{\bxi} = \displmtx_{\lambda}\!\rbr{\bxi} \displprm_{\lambda}, \\[1ex]	
		&\displmtx_{\lambda} = \sbr{\,N^{\lambda}_1 \bI_2 \mid N^{\lambda}_{2} \bI_2\,}, \\[1ex]
		&N_1^{\lambda}\!\rbr{\bxi} = 1-\xi^2, \,\,\, N_2^{\lambda}\!\rbr{\bxi} = 1-\eta^2,
	\end{aligned}
\end{align}
with $\displprm_{\lambda}$ as free parameters, the descending stress interpolation reduces to~\cite{Pian_Sumihara_IJNME_1984}:
\begin{equation}
	\sts\!\rbr{\bxi} = \stsmtx\!\rbr{\bxi}\stsprm, \,\, 
		\stsmtx\!\rbr{\bxi} = \sbr{\begin{matrix}  
			1 & 0 & 0 & a_1^2\eta & a_3^2\xi \\[1.5ex]
			0 & 1 & 0 & b_1^2\eta & b_3^2\xi \\[1.5ex]
			0 & 0 & 1 & a_1b_1\eta &a_3b_3\xi  
		\end{matrix}}\!,
\label{eq:Pian_sts_5b}
\end{equation}
with $a_1, \dots{}, a_3$, $b_1, \dots{}, b_3$ given in equation~\eqref{eq:position_J}. 

\begin{remark}\label{r:sts_dim}
An important point to be investigated in the design of a mixed finite element formulation is the dimension of the stress interpolation space~$N_{\sts}$. In fact, the latter determines the strictness of the element compatibility constraint, and the larger it is, the stiffer the element is expected to be. Furthermore, the larger the number of stress interpolation parameters, the higher the computational cost turns out to be. In many cases, such arguments might motivate the choice to limit the dimension of stress interpolation space as much as possible, coherently with the need to suppress all deformation modes, thus resulting in an isostatic element. Then the Pian filter technique represents a useful tool for conceiving a stress interpolation space with desired dimension, starting from a complete interpolation up to a prescribed degree. Instances of such strategy are e.g.~employed in~\cite{Madeo_Casciaro_FEAD_2012, Nodargi_Bisegna_C&S_2017}. 
\end{remark}

\subsubsection{Self-equilibrated stress interpolation by Airy's function approach}\label{sss:Airy}
A systematic strategy for generating stress interpolation consists in using Airy's function approach. Two main advantages can be highlighted. First, self-equili\-brated stress modes are automatically derived, thus allowing interpolation of the displacement field only over the element boundary (hybrid elements)~\cite{Madeo_Casciaro_FEAD_2012, Nodargi_Bisegna_C&S_2017}. Second, stress modes whose corresponding elastic strain is compatible can be selected at no further expense, thus improving element performances in materially linear problems~\cite{Piltner_Taylor_IJNME_1995, Madeo_Zucco_FEAD_2014}. On the other hand, because the internal equilibrium condition is required to hold over the physical element, the stress interpolation has to be expressed in physical coordinates and some caution is necessary to guarantee frame-invariance. A possible strategy is described in~\cite{Madeo_Casciaro_FEAD_2012}, amounting at the introduction of a local element reference frame. In the following the physical point coordinates in such a frame are still denoted by $\pnt = \cbr{x,y}$.

It is easily shown that the homogeneous modes:
\begin{equation}
	\cbr{1; 0; 0}, \quad \cbr{0;1; 0}, \quad \cbr{0;0;1},
\end{equation}
satisfy internal equilibrium and compatibility in linear problems. On the other hand, the only stress modes of degree $n\geq1$ with the same properties are (for instance, see \cite{Cen_Zhou_CMAME_2011}):
\begin{align}
	\begin{aligned}
		&\cbr{y^n; 0; 0}, \\[1ex] 
		&\cbr{0;x^n; 0}, \\[1ex] 
		&\cbr{0;n x^{n-1}y; -x^n}, \\[1ex]  
		&\cbr{n xy^{n-1}; 0; -y^n}.
	\end{aligned}
\end{align}	
Accordingly, a complete interpolation of stress field up to degree $n$ comprises $3+4n$ stress modes \cite{Nodargi_Bisegna_C&S_2017}. As discussed in Remark~\ref{r:sts_dim}, the Pian filter technique can be then adopted to eliminate overabundant stress modes. Such strategy is for instance exploited in~\cite{Cen_Zhou_CMAME_2011, Madeo_Casciaro_FEAD_2012, Madeo_Zucco_FEAD_2014, Nodargi_Bisegna_C&S_2017}.

\subsection{Strain}
Strain interpolation is required in HW formulations discussed in Section~\ref{ss:HW}. Seeking for an approximation that is at least suitable for linear elastic applications, a common choice consists in an interpolation identical to the one selected for the stress field. On the other hand, an accurate description of the possibly highly nonlinear spatial distribution of the strain field, as arising in nonlinear applications, can be addressed considering independent strain and stress approximations. In particular, a piecewise-constant interpolation is discussed. In the following, ordering $\beps  = \cbr{\eps_{x}; \eps_{y}; \eps_{xy}}$ is adopted for representing the strain field in Voigt notation.

\subsubsection{Identical strain-stress interpolation}\label{sss:identical_stn_sts}
Though stress and strain approximations do not need to satisfy any condition~\cite{Piltner_Taylor_IJNME_1995}, the selection of identical interpolations appears as a reasonable choice in order to ensure high accuracy at least in linear elastic problems (e.g.~see~\cite{Weissman_Jamjian_IJNME_1993, Piltner_CM_2000, Nodargi_Bisegna_IJNME_2016}). Accordingly, the numbers of stress and strain parameters coincide, i.e.~$N_{\sts} = N_{\stn}$, and the strain interpolation matrix is chosen to be:\
\begin{equation}
	\stnmtx = \stsmtx.
\label{eq:identical_stn_sts}
\end{equation}
As a consequence, the stress-strain operator~$\stsstnmtx$, defined in equation~\eqref{eq:comp_stsstn_operator}$\txtsub{2}$, and involved  in both the element compatibility and constitutive relationships, equations \eqref{eq:HW_fun_stat} or~\eqref{eq:enhanced_fun_stat_alt}$\txtsub{1,2}$, is a square symmetric matrix. In Section~\ref{sss:HW_identical_stsstn}, it will be shown that a closed-form procedure is available for the state element determination within a HW formulation with identical strain-stress interpolation, under the assumption that such a matrix is invertible, i.e.~that the stress interpolation matrix~$\stsmtx$ has full-column rank. Moreover, the equivalence with the assumed strain formulation with eliminated stresses discussed in Remark~\ref{r:assumed_stn} is highlighted.

\subsubsection{Piecewise-constant strain interpolation}\label{sss:piecewise_stn}
A piecewise-constant strain interpolation is adopted wi\-thin each finite element in~\cite{Nodargi_Bisegna_C&S_2017}. Such a choice is motivated by the possible highly nonlinear spatial distribution of the strain field, as inherited from the spatial behaviour of internal variables (for instance, plastic strain in an elastoplastic medium). That distribution might be better captured by a discontinuous, piecewise-constant, rather than, e.g., by a smooth, a priori fixed, polynomial interpolation. Similar choices are explored for beam elements in \cite{Spacone_Taucer_EESD_1996, Taylor_Auricchio_CM_2003} and for membrane elements, referring to the plastic multiplier, in \cite{Bilotta_Casciaro_CMAME_2007, Bilotta_Garcea_FEAD_2011}.

For introducing the strain interpolation, a partition of the typical element into $N\qp$ subdomains $\dmn\qp$ is considered by the image through the reference map of an analogous partition of the parent element (a particular choice might be to subdivide the parent element into $N\qp$ equal-sized square subdomains~\cite{Nodargi_Bisegna_C&S_2017}). The strain approximation is assumed to be constant over each subdomain $\dmn\qp$ with value $\stnprm\qp = \cbr{\rbr{e_x}\qp; \rbr{e_y}\qp; \rbr{e_{xy}}\qp}$, whence the dimension of the strain interpolation space results $N_\varepsilon = 3N\qp$, and the relevant interpolation parameters can be collected into the $N_\varepsilon \times 1$ vector:
\begin{equation}
	\stnprm = \cbr{\stnprm_1; \dots{}; \stnprm_{N\qp}}.
\end{equation}
The strain interpolation matrix can be then expressed by:
\begin{equation}
	\stnmtx = \sbr{\, \chi_1\bI_3 \mid \dots{} \mid \chi_{N\qp} \bI_3 \,},
\label{eq:stn_interp}
\end{equation}
in which $\chi\qp$ is the characteristic function of subdomain $\dmn\qp$ and $\bI_3$ is the identity matrix of dimension $3$. Accordingly, the stationary conditions~\eqref{eq:HW_fun_stat} gets simplified into:
\begin{align}
	\begin{aligned}
		\bzero &= {\pdiffarg{\stn}{\ie} \! \rbr{\stnprm\qp} -\avgstsmtx\qp\stsprm} \,, \\[0.5ex]
		\bzero &= - \textstyle\sum_{\qppnt=1}^{N\qp}{\avgstsmtx\qp^{\text{T}} \stnprm\qp \abs{\dmn\qp}} + \compmtx \displprm,
	\end{aligned}
\label{eq:HW_fun_stat_piecewise_stn}
\end{align}
with $\abs{\dmn\qp}$ as the area of the element subdomain $\dmn\qp$ and $\avgstsmtx\qp$ as the averaged stress interpolation matrix $\stsmtx$ over~$\dmn\qp$:
\begin{equation}
	\avgstsmtx\qp = \frac{1}{\abs{\dmn\qp}} \int_{\dmn\qp} \stsmtx \text{d}\dmn.
\end{equation}
\begin{remark}
It is worth noticing that the piecewise-con\-stant strain interpolation reduces constitutive equations \eqref{eq:HW_fun_stat_piecewise_stn}$\txtsub{1}$ to material state update problems over each element subdomain. Unfortunately, no advantage can be directly derived from such a simple structure, because these equations are coupled to each other by compatibility conditions~\eqref{eq:HW_fun_stat_piecewise_stn}$\txtsub{2}$. Conversely, potential convergence difficulties might be expected for the solution of a large-dimension system of coupled nonlinear and nonsmooth equations. That observation motivates the iterative procedure proposed in~\cite{Nodargi_Bisegna_C&S_2017} and reviewed in Section~\ref{ss:element_iem}, amounting at solving the material state update problems in an uncoupled fashion.
\end{remark}

\begin{remark}
As a possible development of the present setting, a piecewise-linear strain interpolation might be considered within each finite element. Though in that case the constitutive equations pertaining to a single element subdomain are coupled to each other, the iterative procedure proposed in~\cite{Nodargi_Bisegna_C&S_2017} might be exploited to solve at least the subdomain equations in an uncoupled format.
\end{remark}

\subsection{Enhanced strain}\label{ss:enh_stn_interp}
An interpolation of the enhanced strain field is required in the ES formulation presented in Section~\ref{ss:enhanced_strain}. A customary strategy for its design consists in the derivation from a prescribed incompatible displacement field, \emph{de facto} following the method of incompatible modes originally proposed in~\cite{Wilson_Ghaboussi_1973}. Accordingly, for the incompatible displacement field given by:
\begin{equation}
	\displ\txtsup{i}\!\rbr{\bxi} = \displmtx\txtsup{i}\!\rbr{\bxi} \displprm\txtsup{i},
\end{equation}
the enhanced strain interpolation follows as:
\begin{align}
	\begin{aligned}
		&\tilde\stn\!\rbr{\bxi} = \graddisplmtx\txtsup{i}\!\rbr{\bxi} \displprm\txtsup{i}, \\[0.5ex]
		&\graddisplmtx\txtsup{i}\!\rbr{\bxi} = \bJ^{-\text{T}}\!\rbr{\bxi} \, \sym\nabla_{\bxi} \displmtx\txtsup{i}\!\rbr{\bxi} =
			\frac{1}{J\!\rbr{\bxi}} \,\bJ\conj\!\rbr{\bxi} \sym\nabla_{\bxi} \displmtx\txtsup{i}\!\rbr{\bxi}
	\end{aligned}
\label{eq:B_i}
\end{align}
where $\bJ$ is the gradient of the isoparametric map~$\pnt=\pnt\!\rbr{\bxi}$, $J = \det\bJ$ and $\bJ\conj$ is the cofactor matrix of $\bJ$.

An alternative approach for the selection of the enhanced strain field interpolation is proposed in~\cite{Simo_Rifai_IJNME_1990}. As starting point, an arbitrary interpolation matrix~$\tilde\stnmtx_{\bxi}$ is defined in isoparametric space. Then, interpreting the isoparametric map $\pnt = \pnt\!\rbr{\bxi}$ as a deformation, the interpolation matrix in physical space~$\tilde\stnmtx$ is derived by a push-forward of the interpolation matrix in isoparametric space~$\tilde\stnmtx_{\bxi}$:
\begin{align}
	\begin{aligned}
		&\tilde\stn\!\rbr{\bxi} = \tilde\stnmtx\!\rbr{\bxi} \tilde\stnprm, \\[0.5ex]
		&\tilde\stnmtx\!\rbr{\bxi} = \bJ^{-\text{T}}\!\rbr{\bxi} \tilde\stnmtx_{\bxi}\!\rbr{\bxi} \bJ^{-1}\!\rbr{\bxi}.
	\end{aligned}
\label{eq:tilde_E}
\end{align}
Such transformation guarantees the existence of an incompatible displacement field in isoparametric space~$\displ\txtsup{i}$ such that, via equation~\eqref{eq:B_i}, the same enhanced strain interpolation matrix is derived in physical space, i.e. $\graddisplmtx\txtsup{i} = \tilde\stnmtx$ \cite{Simo_Rifai_IJNME_1990}.

\begin{remark}
The original choice of the incompatible displacement functions~$\displ\txtsup{i}$ traces back to the work by Wilson and his coauthors~\cite{Wilson_Ghaboussi_1973}, where it is assumed that:
\begin{align}
	\begin{aligned}
		&\displmtx\txtsup{i}\!\rbr{\bxi} = \sbr{\,N\txtsup{i}_1 \bI_2 \mid N\txtsup{i}_{2} \bI_2\,}, \\[1ex]
		&N_1\txtsup{i}\!\rbr{\bxi} = 1-\xi^2, \,\,\, N_2\txtsup{i}\!\rbr{\bxi} = 1-\eta^2.
	\end{aligned}
\end{align}
However in~\cite{Taylor_Wilson_IJNME_1976} it is observed that the constant strain patch test cannot be satisfied when using the exact gradient of such incompatible shape functions. Accordingly the following approximate gradient operator is used~\cite{Taylor_Wilson_IJNME_1976}:
\begin{equation}
	\graddisplmtx\txtsup{i}\!\rbr{\bxi} = \frac{J\!\rbr{\bzero}}{J\!\rbr{\bxi}} \,\bJ\conj\!\rbr{\bzero} \sym\nabla_{\bxi} \displmtx\txtsup{i}\!\rbr{\bxi}.
\end{equation}
and, correspondingly, equation~$\eqref{eq:tilde_E}$ is modified in~\cite{Simo_Rifai_IJNME_1990}:
\begin{equation}
	\tilde\stnmtx\!\rbr{\bxi} = \frac{J\!\rbr{\bzero}}{J\!\rbr{\bxi}} \bJ^{-\text{T}}\!\rbr{\bzero} \tilde\stnmtx_{\bxi}\!\rbr{\bxi} \bJ^{-1}\!\rbr{\bzero},
\end{equation}
Enhanced strain interpolations that can be set within such framework are adopted e.g.~in~\cite{Piltner_Taylor_IJNME_1995, Piltner_Taylor_IJNME_1999, Piltner_CM_2000}.
\end{remark}
\begin{remark}
As discussed in Section~\ref{ss:enhanced_strain}, the assumption of $L^2$-orthogonality between stress and enhanced strain interpolations is a sufficient condition for eliminating the stress parameters~$\stsprm$ from the ES formulation. In addition, a sufficient condition for the constant patch test satisfaction is the enhanced strain interpolation matrix~$\tilde\stnmtx$ to have vanishing average over the physical element~\cite{Simo_Rifai_IJNME_1990}. Both conditions can be enforced after representations~\eqref{eq:B_i} or~\eqref{eq:tilde_E} are derived.
\end{remark}

\subsection{Plastic multiplier} \label{ss:plastic_mlt}
As a peculiar feature of the CM formulation discussed in Section~\ref{ss:return_map}, the plastic flow law is weakly enforced at element level. Conversely pointwise or element imposition of plastic admissibility and plastic consistency conditions can be tuned on the basis of the interpolation adopted for the plastic multiplier field. In this section, particular choices explored in the literature are detailed.

\subsubsection{Pointwise or quadrature-point Dirac-delta interpolation}\label{sss:plastic_mlt_pointwise}
In the original formulation of the CM formulation, proposed in~\cite{Simo_Taylor_CMAME_1989}, a pointwise interpolation of the plastic multiplier field is considered, whence conditions~\eqref{eq:return_map_stat} can be written as:
\begin{align}
	\begin{aligned}
		&\bzero = \compmtx\displprm - \stnprm\p\n - 
			\intbr{\stsmtx\txtsup{T} \partial_{\sts}\yf\!\rbr{\stsmtx \stsprm} \Delta\lambda\, } - 
			\cmplmtx\e\stsprm, \\[1ex]
		&\Delta\lambda \geq 0, \quad
		\yf\!\rbr{\sts} \leq 0, \quad
		\Delta\lambda \yf\!\rbr{\sts} = 0.
	\end{aligned}
\label{eq:plastic_admissibility_Simo}	
\end{align}
In practice, numerical quadrature is employed for evaluating integrals and the plastic admissibility condition is pointwise enforced only at quadrature points. Accordingly, denoting by $\cbr{\pnt_1, \cdots{}, \pnt_{N\gp}}$ the set of element quadrature points and $\cbr{w_1, \cdots{}, w_{N\gp}}$ the set of corresponding quadrature weights, it follows that:
\begin{align}
	\begin{aligned}
		&\bzero = \compmtx\displprm - \stnprm\p\n - 
			\textstyle\sum_{g=1}^{N\gp}{ \flowmtx\gp\!\rbr{\stsprm} \Delta \lambda\gp w\gp} - 
			\cmplmtx\e\stsprm, \\[1ex]
		&\Delta \lambda\gp \geq 0, \quad
			\yf\gp\!\rbr{\stsprm} \leq 0, \quad
			\Delta \lambda\gp \yf\gp\!\rbr{\stsprm} = 0,
	\end{aligned}
\label{eq:plastic_admissibility_Simo_2}	
\end{align}
where:
\begin{align}
	\begin{aligned}
		&\yf\gp\!\rbr{\stsprm} = \yf\!\rbr{\stsmtx\gp \stsprm}, \\[1ex]
		&\flowmtx\gp\!\rbr{\stsprm} = \pdiff{\yf\gp}{\stsprm}\!\rbr{\stsprm} = \stsmtx\gp\txtsup{T}\,\partial_{\sts}\yf\!\rbr{\stsmtx\gp \stsprm},
	\end{aligned}
\label{eq:plastic_admissibility_notation}
\end{align}
are respectively introduced for yield function and plastic flow operator at quadrature point~$\pnt\gp$. It is remarked that plastic admissibility and plastic consistency conditions are pointwise enforced at quadrature points, whereas an element flow law is involved.

\begin{remark}
An alternative choice investigated in~\cite{Simo_Taylor_CMAME_1989} amounts at a quadrature-point Dirac-delta interpolation for the interpolation of the plastic multiplier, i.e.~to the following interpolation matrix:
\begin{equation}
	\pmltmtx\!\rbr{\pnt} = \sbr{\delta_1\!\rbr{\pnt}, \dots{}, \delta_{N\gp}\!\rbr{\pnt}},
\end{equation}
where $\delta_g\!\rbr{\pnt} = \delta\!\rbr{\pnt-\pnt\gp}$ and $\delta$ denotes the Dirac delta measure (centered at the origin). As a consequence, conditions~\eqref{eq:return_map_stat} result to be:
\begin{align}
	\begin{aligned}
		&\bzero = \compmtx\displprm - \stnprm\p\n - 
			\textstyle\sum_{g=1}^{N\gp}{ \flowmtx\gp\!\rbr{\stsprm} \Delta l\gp} - 
			\cmplmtx\e\stsprm, \\[1ex]
		&\Delta l\gp \geq 0, \quad
			\yf\gp\!\rbr{\stsprm} \leq 0,\quad
			\Delta l\gp \yf\gp\!\rbr{\stsprm} = 0,
	\end{aligned}
\label{eq:plastic_admissibility_Simo_App}	
\end{align}
Equations~\eqref{eq:plastic_admissibility_Simo_App} only differ from conditions~\eqref{eq:plastic_admissibility_notation} for the presence of the quadrature weights~$w_g$ in the plastic flow term of compatibility relationship. In particular, in the case that a $2 \times 2$ quadrature is used over the element ($w_g=1$), the two discretizations coincide.
\end{remark}

\subsubsection{Piecewise-constant interpolation}
Another strategy, explored in~\cite{Bilotta_Garcea_FEAD_2011}, consists in a piecewise constant interpolation on the element subdomains $\dmn\qp$, $d=1, \dots{}, N\qp$, defined as in Section~\ref{sss:piecewise_stn}. In such a case, the plastic multiplier interpolation matrix is:
\begin{equation}
	\pmltmtx\!\rbr{\pnt} = \sbr{\chi_1\!\rbr{\pnt}, \dots{}, \chi_{N_d}\!\rbr{\pnt}},
\end{equation}
where it is recalled that $\chi\qp$ denotes the characteristic function associated to the element subdomain~$\dmn\qp$. Correspondingly, conditions~\eqref{eq:return_map_stat} boil down to:
\begin{align}
	\begin{aligned}
		&\bzero = \compmtx\displprm - \stnprm\p\n -
			\textstyle\sum_{d=1}^{N\qp}\overline{\flowmtx}\qp\!\rbr{\stsprm} \Delta l\qp \, \abs{\dmn\qp} -
			\cmplmtx\e\stsprm,\\[1ex]
		&\Delta l\qp \geq 0, \quad
			\overline{\yf}\qp\!\rbr{\stsprm} \leq 0,\quad
			\Delta l\qp \, \overline{\yf}\qp\!\rbr{\stsprm} \abs{\dmn\qp} = 0,
	\end{aligned}
\label{eq:plastic_admissibility_Bilotta}		
\end{align}
in which:
\begin{align}
	\begin{aligned}
		&\overline{\yf}\qp\!\rbr{\stsprm} = \frac{1}{\abs{\dmn\qp}}\int_{\dmn\qp}{\yf\!\rbr{\stsmtx\stsprm} \text{d}\dmn}, \\[1ex]
		&\overline{\flowmtx}\qp\!\rbr{\stsprm} = \pdiff{\overline{\yf}\qp}{\stsprm}\!\rbr{\stsprm} = \frac{1}{\abs{\dmn\qp}} {\int_{\dmn\qp}{\stsmtx\txtsup{T}\partial_{\sts}\yf\!\rbr{\stsmtx\stsprm}}\text{d}\dmn},
	\end{aligned}
\end{align}
respectively denote the mean values of the element yield function and of the element plastic flow operator over the element subdomain~$\dmn\qp$. Neither the plastic admissibility nor the plastic consistency conditions are imposed pointwise.

\begin{remark}\label{r:pltmlt_constant}
As a particular instance of the piecewise-constant interpolation, constant interpolation of the plastic multiplier can be considered over the finite element \cite{Bilotta_Casciaro_CMAME_2007}. As a consequence, i.e. assuming $N\qp = 1$, the plastic multiplier interpolation matrix results to be~$\pmltmtx\!\rbr{\pnt} = \sbr{1}$, whence the system of equation~\eqref{eq:plastic_admissibility_Bilotta} reduces to:
\begin{align}
	\begin{aligned}
		&\bzero = \compmtx\displprm - \stnprm\p\n -
			\overline{\flowmtx}\!\rbr{\stsprm} \Delta l \, \abs{\dmn} -
			\cmplmtx\e\stsprm,\\[1ex]
		&\Delta l \geq 0, \quad
			\overline{\yf}\!\rbr{\stsprm} \leq 0, \quad
			\Delta l \, \overline{\yf}\!\rbr{\stsprm} \abs{\dmn} = 0, 
	\end{aligned}
\label{eq:plastic_admissibility_Bilotta_2}		
\end{align}
where it is intended that $\dmn = \dmn_1$, $\Delta l = \Delta l_1$, $\overline{\yf}\!\rbr{\stsprm} = \overline{\yf}_1\!\rbr{\stsprm}$ and $\overline{\flowmtx}\!\rbr{\stsprm} = \overline{\flowmtx}_1\!\rbr{\stsprm}$. A significant difference can be here highlighted in comparison to the other discussed plastic multiplier approximations. In fact, because of the one-parameter interpolation, the element yield function loses its vector form and reduces to an element-avaraged counterpart of the (material) yield function. That allows for a great simplification of the numerical solution strategy of CM functional discrete problem, as discussed in Remark~\ref{r:constant_pmlt}.
\end{remark}

\section{Numerical solution strategy}\label{s:numerical_strategy}
In this section, attention is devoted to numerical schemes employed for the solution of the discrete structural problems resulting from the mixed methods in Sections~\ref{s:gsm_variational} and~\ref{s:interpolation}. 

Apart for intrinsic differences due to the different underlying variational formulations, the common feature of continuous displacement interpolation across element boundaries and element-supported interpolation of additional fields (stress, strain, enhanced strain, plastic multiplier), implies the feasibility of a common displacement driven architecture for the solution strategy. In particular, two levels of solution can be recognized. At element level, additional fields parameters can be statically condensed out for given nodal DOFs and element nodal forces can be consequently computed. Such procedure is referred to as element state determination. Conversely, at structural level, global equilibrium equations, which involve the nodal forces (obtained assembling element contributions) as a function of nodal DOFs, are enforced. Yet the separation of the two levels of solution will be alluded to for the sake of clarity in the following developments, it should not be overemphasized, as there will emerge numerical schemes that tend to blur element and structural levels. In this regard, the key point consists in choosing whether the element state determination has to be fulfilled up to prescribed tolerance or element residuals relevant to the additional fields are allowed to be assembled at structural level. 

In general, it is desirable that a numerical algorithm is \emph{accurate}, \emph{robust} and \emph{efficient}. Several numerical solution strategies have been proposed in the literature for the solution of the element state determination problems under consideration. Specifically, the following approaches can be distinguished: \\[-1ex]
\begin{itemize}
	\item[1.]{\emph{Closed-form solutions}};
	\item[2.]{\emph{Newton's method}};
	\item[3.]{\emph{Element return mapping}};
	\item[4.]{\emph{Mathematical programming}};
	\item[5.]{\emph{Nodal-force-based algorithm}}.\\[-1ex]
\end{itemize}

Closed-form solutions are generally out of reach when considering mixed methods for inelastic structures. In fact, the only instance deals with a HW formulation with identical stress-strain interpolation~\cite{Nodargi_Bisegna_IJNME_2016}, which is shown to be equivalent to a B-bar formulation with eliminated stresses, and hence substantially consists in a displacement-based formulation.

Newton's method represents of course the mainstream option for the solution of nonlinear systems of equations, resulting especially attractive for its local asymptotic quadratic rate of convergence. Indeed, except for mathematical programming strategy, Newton's method is involved, at least in disguise, in all the numerical approaches here listed. Nevertheless, its direct application results in a robust procedure just for ES formulations~\cite{Simo_Hughes_JAM_1986, Simo_Rifai_IJNME_1990, Piltner_Taylor_IJNME_1995, Piltner_Taylor_IJNME_1999, Piltner_CM_2000} and HR formulations~\cite{Schroder_Miehe_CMAME_1997, Schroder_Starke_CMAME_2017}.

Contrarily, referring to CM formulations, the element compatibility equation along with element plastic admissibility and consistency conditions amount to a problem analogous to non-smooth multi-surface plasticity at material point level. In particular, a direct application of Newton's method is precluded because no reliable prediction can be made on the set of active constraints.

In~\cite{Simo_Taylor_CMAME_1989, Bilotta_Casciaro_CMAME_2007, Bilotta_Garcea_FEAD_2011}, an element return mapping strategy is adopted, resembling the classical elastic predictor-plastic corrector scheme employed in elastoplasticity. Accordingly, in the plastic corrector stage, Newton's iterations are performed as a nested loop within the iterations for determining the set of active constraints.

An alternative approach resorts to a rephrasing of the discrete problem in a closest-point-projection format and to the adoption of numerical algorithms from mathematical programming~\cite{Krabbenhoft_Wriggers_IJNME_2007, Mendes_Castro_FEAD_2009, Bilotta_Garcea_FEAD_2011, Bilotta_Garcea_CS_2012}. In fact, if during the nineteen-eighties simplex method and its derivatives were the only possibilities for the solution of general mathematical programs, suffering from a significant increase in computational effort for increasing problem size, many appealing new methods of optimization are now available and represent a competitive option~\cite{Krabbenhoft_Wriggers_IJNME_2007}.

In case of a HW formulation, the coupled structure of quite many nonlinear and nonsmooth constitutive equations makes the element state determination problem a difficult task to be accomplished through a direct application of Newton's method. That difficulty is even more pronounced if a strain interpolation space of large dimension is involved, as desirable for a more accurate description of the strain field at element level.  A possible remedy is offered by the iterative procedure, here referred to as nodal-force-based algorithm, presented in~\cite{Nodargi_Bisegna_C&S_2017} and inspired by the techniques discussed in~\cite{Taylor_Auricchio_CM_2003, Spacone_Taucer_EESD_1996, Neuenhofer_Filippou_JSE_1997, Saritas_Soydas_IJNLM_2012, Nodargi_Bisegna_PAMM_2015}. Relying on a piecewise constant strain interpolation, that procedure allows for independent solution of as many material state update problems as the number of element subdomains, thus reducing the computational cost which stems from coupled nonlinear constitutive equations, and mitigating convergence difficulties of Newton's method.

In the following, the listed approaches are analyzed and discussed adopting the notation introduced in Sections~\ref{s:gsm_variational} and~\ref{s:interpolation}.

\subsection{Closed-form solutions}
A closed-form solution of the element state determination problem, i.e.~a solution not requiring any element iteration, is generally unfeasible for mixed formulation concerning inelastic structures. However, when considering a HW formulation with identical strain-stress interpolation, the method reduces to a B-bar formulation and such a closed-form procedure can be devised. 

\subsubsection{Hu--Washizu formulation with identical stress-strain interpolation}\label{sss:HW_identical_stsstn}
The procedure developed in \cite{Nodargi_Bisegna_IJNME_2016} is here reviewed. As for identical strain-stress interpolation, the stress-strain operator~$\stsstnmtx$ introduced in equation~\eqref{eq:comp_stsstn_operator}$\txtsub{2}$ results to be a square symmetric matrix of dimension~$N_{\sts} = N_{\stn}$ (see Section~\ref{sss:identical_stn_sts}). Assuming it to be invertible, i.e.~the stress interpolation matrix~$\stsmtx$ to have full-column rank, the element compatibility condition~\eqref{eq:HW_fun_stat}$\txtsub{2}$ can be solved for the strain parameters~$\stnprm$, yielding:
\begin{equation}
	\stnprm = \stsstnmtx^{-\text{T}}\compmtx\displprm,
\end{equation}
whence the stress parameters~$\stsprm$ follow from the element constitutive equations~\eqref{eq:HW_fun_stat}$\txtsub{1}$:
\begin{equation}
	\stsprm = \stsstnmtx^{-1} \intbr{\stnmtx\txtsup{T}\partial_{\stn}\ie \! \rbr{\stnmtx\stsstnmtx^{-\text{T}}\compmtx\displprm}}.
\label{eq:HW_identical_sts}			
\end{equation}
Finally, the element internal forces~$\force\internal$ can be derived from the element equilibrium equations~\eqref{eq:HW_q_int}:
\begin{equation}
	\force\internal = \compmtx\txtsup{T} \stsstnmtx^{-1} \intbr{\stnmtx\txtsup{T}\partial_{\stn}\ie \! \rbr{\stnmtx\stsstnmtx^{-\text{T}}\compmtx\displprm}}, 
\label{eq:HW_identical_force}	
\end{equation}
and the element consistent stiffness matrix, to be used at structural level for the imposition of the global equilibrium, follows from their differentiation with respect to the element nodal DOFs~$\displprm$:
\begin{align}
	\begin{aligned}
		\stiff = 
		\compmtx\txtsup{T} \stsstnmtx^{-1} \sbr{\intbr{\stnmtx\txtsup{T}\partial^2_{\stn\stn}\ie \! \rbr{\stnmtx\stsstnmtx^{-\text{T}}\compmtx\displprm}\stnmtx\,}} \stsstnmtx^{-\text{T}}\compmtx.
	\end{aligned}
\label{eq:HW_identical_stiff}		
\end{align}

\begin{remark}
Considering the assumed strain formulation discussed in Remark~\ref{r:assumed_stn} under the assumption of stress elimination ($\compmtx=\bar\compmtx$), the element internal forces~$\force\internal$ and the element stiffness matrix~$\stiff$ follow from equation~\eqref{eq:assumed_strain_q_int}:
\begin{align}
	\begin{aligned}
		&\force\internal = \intbr{\bar\graddisplmtx\txtsup{T} \partial_{\stn}\ie\!\rbr{\bar\graddisplmtx\displprm}}, \\[0.5ex]
		&\stiff = \intbr{\bar\graddisplmtx\txtsup{T} \partial^2_{\stn\stn}\ie\!\rbr{\bar\graddisplmtx\displprm}\bar\graddisplmtx\,}.
	\end{aligned}
\label{eq:assumed_stn_no_sts}
\end{align}
By comparison of expressions~\eqref{eq:HW_identical_force} and~\eqref{eq:HW_identical_stiff} with equations~\eqref{eq:assumed_stn_no_sts}, it can be concluded that HW formulations with identical stress-strain interpolation are equivalent to assumed strain formulations with eliminated stresses provided the B-bar matrix is given by:
\begin{equation}
	\bar\graddisplmtx = \stnmtx\stsstnmtx^{-\text{T}}\compmtx.
\label{eq:HW_identical_Bbar}
\end{equation}
A crucial advantage of the HW setting is that the stress parameters~$\stsprm$ can be directly computed from equation \eqref{eq:HW_identical_sts}, with no need of ad-hoc stress recovery procedures, as required in assumed strain formulations. In turn, from the point of view of assumed strain methods, the selection of the B-bar matrix in equation~\eqref{eq:HW_identical_Bbar} ensures the formulation to be variationally consistent in the sense of~\cite{Simo_Hughes_JAM_1986}. In fact, it is easily checked that $\bar\compmtx = \compmtx$.
\end{remark}

\subsection{Newton's method}
Newton's method represents the main option for the solution of the element state determination problem, offering local asymptotic quadratic rate of convergence. However, its direct application results in a robust procedure only when considering ES and HR formulations, as detailed in this section.

\subsubsection{Enhanced strain formulation}\label{sss:Newton_enh}
The solution procedure of element state determination equations in ES formulations, as proposed in~\cite{Simo_Rifai_IJNME_1990}, is here presented. In particular, referring to the discussion in Section~\ref{ss:enhanced_strain}, it is assumed that $L^2$-orthogonal interpolation of enhanced strain and stress fields has been selected, whence equation~\eqref{eq:enhanced_fun_stat}$\txtsub{2}$ is identically satisfied and the stress parameters~$\stsprm$ are not being explicitly involved in equation~\eqref{eq:enhanced_fun_stat}$\txtsub{1}$.
The latter equation is solved by means of Newton's method for the unknown enhanced strain parameters~$\tilde\stnprm$. Starting from the residual formulation:
\begin{equation}
	\bzero = \res_{\tilde\stn} = \intbr{\tilde\stnmtx\txtsup{T} \partial_{\stn}\ie(\graddisplmtx\displprm + \tilde\stnmtx \tilde\stnprm)\,},
\label{eq:enhanced_fun_stat_no_sts}
\end{equation}
the solution approximation at $j$-th iteration is given by the following linearized equation:
\begin{equation}
	\res_{\tilde\stn}^{j+1} \approx \res_{\tilde\stn}^{j} + 
		\stiff_{\tilde\stn}\!\rbr{\tilde\stnprm^j} \rbr{\tilde\stnprm^{j+1}-\tilde\stnprm^j} = \bzero,
\end{equation}
in which the enhanced strain stiffness matrix is introduced by:
\begin{equation}
	\stiff_{\tilde\stn}\!\rbr{\tilde\stnprm} = \intbr{\tilde\stnmtx\txtsup{T} \partial^{2}_{\stn\stn}\ie(\graddisplmtx\displprm + \tilde\stnmtx \tilde\stnprm) \tilde\stnmtx \,}.
\label{eq:enh_K_eps}
\end{equation}
When the convergence is reached within a fixed tolerance, the element nodal internal forces~$\force\internal$ are computed from equation~\eqref{eq:enhanced_q_int} and standard assembling procedure yields structural equilibrium residuum. The element consistent stiffness matrix to be used in its linearization reads:
\begin{equation}
	\stiff = 
		\stiff_{\displ} + \stiff_{\tilde\stn,\displ}\txtsup{T} \, \pdiff{\tilde\stnprm}{\displprm},
\end{equation}
where $\stiff_{\displ}$ and~$\stiff_{\tilde\stn,\displ}$ respectively denote the compatible strain stiffness matrix and the coupling enhanced-compatible strain stiffness matrix, defined by:
\begin{align}
	\begin{aligned}
		&\stiff_{\displ}\!\rbr{\tilde\stnprm} = \intbr{\graddisplmtx\txtsup{T} \partial^{2}_{\stn\stn}\ie(\graddisplmtx\displprm + \tilde\stnmtx \tilde\stnprm) \graddisplmtx\,}, \\[1ex]
		&\stiff_{\tilde\stn\displ}\!\rbr{\tilde\stnprm} = \intbr{\tilde\stnmtx\txtsup{T} \partial^{2}_{\stn\stn}\ie(\graddisplmtx\displprm + \tilde\stnmtx \tilde\stnprm) \graddisplmtx\,}.
	\end{aligned}
\label{eq:enh_K_u_K_eps_u}	
\end{align}
For the computation of the derivative of the enhanced strain parameters~$\tilde\stnprm$, equation~\eqref{eq:enhanced_fun_stat_no_sts} is differentiated:
\begin{equation}
	 \pdiff{\tilde\stnprm}{\displprm} = - \stiff_{\tilde\stn}^{-1} \stiff_{\tilde\stn,\displ}
\end{equation}
whence it is obtained that:
\begin{equation}
	\stiff = \stiff_{\displ} - \stiff_{\tilde\stn,\displ}\txtsup{T} \stiff_{\tilde\stn}^{-1} \stiff_{\tilde\stn,\displ}.
\end{equation}

\begin{remark}
As stresses do not explicitly enter the element equations~\eqref{eq:enhanced_fun_stat_no_sts}, their recovery at element level is required. A variationally consistent procedure valid in case of linearly elastic material and based on a least-square formulation is discussed in~\cite{Simo_Rifai_IJNME_1990}.
\end{remark}

\begin{remark}
In~\cite{Piltner_Taylor_IJNME_1995, Piltner_Taylor_IJNME_1999, Piltner_CM_2000} Newton's method of solution is applied to the ES formulation~\eqref{eq:enhanced_fun_stat_alt}. In particular, under the assumption of equal stress and strain interpolations, equation~\eqref{eq:identical_stn_sts}, the compatibility equation~\eqref{eq:enhanced_fun_stat_alt}$\txtsub{2}$ can be solved for the strain parameters~$\stnprm$, whence equations~\eqref{eq:enhanced_fun_stat_alt}$\txtsub{1,3}$ constitute a nonlinear system in the unknowns stress parameters~$\stsprm$ and enhanced strain parameters~$\tilde\stnprm$.
\end{remark}

\subsubsection{Hallinger--Reissner formulation}\label{sss:Newton_HR}
Referring to the HR formulation discussed in Section~\ref{ss:HR}, and considering in particular an elastoplastic structural problem (see Remark~\ref{r:HR_Schroder_Starke}), two numerical solution strategies hinged on Newton's method are proposed in \cite{Schroder_Miehe_CMAME_1997, Schroder_Starke_CMAME_2017}. Specifically, the two approaches share the use of Newton's method at structural level, but differ on the choice to solve exactly or not, at each structural iteration, the element equations~\eqref{eq:HR_stat}. In this section, such schemes are briefly reviewed and compared to each other.

In~\cite{Schroder_Miehe_CMAME_1997}, an exact solution, i.e.~up to prescribed tolerance, of element equations is pursued by local Newton's iterations. To this end, after recasting equations~\eqref{eq:HR_stat} in residual format:
\begin{equation}
	\bzero = \res_{\sts} = \compmtx\displprm - \stnprm\p\n - \Delta\stnprm\p\!\rbr{\stsprm} - \cmplmtx\e\stsprm,
\label{eq:HR_stat_res}
\end{equation}
linearization with respect to the unknown stress parameters~$\stsprm$ is performed:
\begin{equation}
	\res_{\sts}^{j+1} \approx \res_{\sts}^{j} - \cmplmtx\ep\!\rbr{\stsprm^j}\rbr{\stsprm^{j+1}-\stsprm^j} = \bzero.
\end{equation}
Here the elastoplastic compliance matrix is introduced by:
\begin{equation}
	\cmplmtx\ep\!\rbr{\stsprm} = \cmplmtx\e + \pdiff{\Delta\stnprm\p}{\stsprm}\!\rbr{\stsprm}, 
\end{equation}
the derivative of the discrete increment of plastic strain $\Delta\stnprm\p$ with respect to the stress parameters~$\stsprm$ following from equation~\eqref{eq:HR_plstn_prm}:
\begin{equation}
	\pdiff{\Delta\stnprm\p}{\stsprm}\!\rbr{\stsprm} = \intbr{\stsmtx\txtsup{T} \, \pdiff{\Delta\stn\p}{\sts}\!\rbr{\stsmtx\stsprm} \stsmtx}.
\end{equation}
After convergence is reached, the element internal nodal forces~$\force\internal$ are computed by means of the element equilibrium equations~\eqref{eq:HR_q_int}, and standard assembling procedure is carried out for imposing structural equilibrium. In particular, observing that from equation~\eqref{eq:HR_stat_res} the derivative of the stress parameters~$\stsprm$ with respect to the DOFs~$\displprm$ is given by:
\begin{equation}
	\pdiff{\stsprm}{\displprm} = \rbr{\cmplmtx\ep}^{-1}\compmtx,
\end{equation}
the element consistent tangent stiffness matrix can be derived from the element equilibrium equations~\eqref{eq:HR_q_int}:
\begin{equation}
	\stiff = \compmtx\txtsup{T} \,  \rbr{\cmplmtx\ep}^{-1}\compmtx.
\label{eq:HR_el_stiff}
\end{equation}

Though the idea that structural equilibrium has to be imposed on a configuration which already satisfies element level equations is customarily accepted in computational mechanics, one main concern can be raised. In fact, significant computational cost is potentially required for achieving convergence at local level, even if a poor structural approximation is still under investigation~\cite{Krabbenhoft_Wriggers_IJNME_2007}. Such consideration motivates the approach recently proposed in~\cite{Schroder_Starke_CMAME_2017}, which amounts at the elimination of local Newton's iterations. Accordingly, the nonlinear system constituted by structural equilibrium equations and element compatibility conditions~\eqref{eq:HR_stat} is tackled in the following residual form:
\begin{align}
	\begin{aligned}
		\bzero &= \res_{\displ} = \cA_{e=1}^{N\elem} \rbr{\compmtx\elem\txtsup{T} \stsprm\elem}, \\[1ex]
		\bzero &= \rbr{\res_{\sts}}\elem = \compmtx\elem\displprm\elem - \rbr{\stnprm\p\n}\elem - \Delta\stnprm\p\elem\!\rbr{\stsprm\elem} - \cmplmtx\elem\e\stsprm\elem,
	\end{aligned}
\end{align}
where $\cA$ stands for the standard assembling operator of finite element method and the subscript $\rbr{\bullet}\elem$ refers to element level variables. Upon linearizing with respect to the unknown displacement and stress parameters, respectively $\displprm$ and $\stsprm\elem, \,\, e =1, \dots{}, N\elem$, and denoting by $j$ the iteration index, it results that:
\begin{align}
	\begin{aligned}
		&\res_{\displ}^{j+1} \approx \res_{\displ}^{j} + \cA_{e=1}^{N\elem} \sbr{\compmtx\elem\txtsup{T} \rbr{\stsprm\elem^{j+1} - \stsprm\elem^j}} = \bzero, \\[1ex]
		&\rbr{\res_{\sts}}\elem^{j+1} \approx \rbr{\res_{\sts}}\elem^{j} + 
			\compmtx\elem\rbr{\displprm\elem^{j+1} - \displprm\elem^j} \\[1ex]
			&\hspace{2.55cm}-\cmplmtx\ep\elem\!\rbr{\stsprm\elem^j}\rbr{\stsprm\elem^{j+1}-\stsprm\elem^j} = \bzero.
	\end{aligned}
\end{align}
In particular, exploiting the local character of the stress interpolation, the element compatibility conditions can be solved for the stress parameters:
\begin{equation}
	\stsprm\elem^{j+1} - \stsprm\elem^j = \sbr{\cmplmtx\ep\elem\!\rbr{\stsprm\elem^j}}^{-1} 
		\sbr{\rbr{\res_{\sts}}\elem^{j} + \compmtx\elem\rbr{\displprm\elem^{j+1} - \displprm\elem^j}},
\end{equation}
and the substitution in the structural equilibrium equations yields:
\begin{equation}
	\displprm^{j+1} =  \displprm^j - 
		\sbr{\cA_{e=1}^{N\elem} \stiff\elem^{j}}^{-1} \tilde{\res}_{\displ}^{j},
\end{equation}
where the consistent element stiffness matrix~$\stiff\elem^{j}$ has the same expression given in equation~\eqref{eq:HR_el_stiff} and a modified global residuum is introduced by: 
\begin{equation}
	\tilde{\res}_{\displ}^{j} = \res_{\displ}^{j} + \cA_{e=1}^{N\elem} \,\,
	\compmtx\elem\txtsup{T} 	\sbr{\cmplmtx\ep\elem\!\rbr{\stsprm\elem^j}}^{-1}\!\rbr{\res_{\sts}}\elem^{j}.
\end{equation}
It is worth noting that the modified global residuum~$\tilde{\res}_{\displ}$ accounts for both the contributions arising from the unbalanced stress state $\res_{\displ}$ and from the incompatible strains at the element level~$\rbr{\res_{\sts}}\elem$, which can be interpreted as inducing an unbalanced pre-stress. A possible drawback of such solution strategy is that, because the element equations are not exactly satisfied, substituting the relevant results in the global equilibrium equations make them deeply nonlinear in the unknown nodal DOFs~\cite{Bilotta_Garcea_CS_2012}. Consequently, the robustness of the algorithm might be worsened and a larger number of structural iterations might be necessary for achieving convergence.

\subsection{Element return mapping} \label{ss:element_return_mapping}
Within the context of a CM variational formulation (Section~\ref{ss:return_map}), the solution at element level of the compatibility equation~\eqref{eq:return_map_stat}$\txtsub{1}$ (involving a weak enforcement of the plastic flow law) and of the plastic admissibility conditions~\eqref{eq:return_map_stat}$\txtsub{2}$, can be performed by resorting to an elastic predictor/plastic corrector strategy, as inspired by the procedure customarily adopted for the integration of the elastoplastic evolution equations (for instance, see~\cite{Karaoulanis_ARCME_2013}). In~\cite{Simo_Taylor_CMAME_1989}, an algorithm assuming pointwise interpolation of the plastic multiplier field (Section~\ref{sss:plastic_mlt_pointwise}) and von Mises yield criterion is proposed. A natural generalization of that approach is here discussed.

\subsubsection{Elastic predictor/plastic corrector algorithm} \label{ss:element_return_mapping_sub}
The elastic predictor stage consists in the computation, from equation~\eqref{eq:return_map_stat}$\txtsub{1}$, of the trial stress corresponding to the assumption of no plastic evolution, i.e.~for $\pmltprm = \bzero$:
\begin{equation}
	\stsprm\trial = \rbr{\cmplmtx\e}^{-1}\rbr{\compmtx\displprm - \stnprm\p\n}.
\label{eq:trial_stsprm}
\end{equation}
In case the trial stress turns out to be plastically admissible, that is condition~\eqref{eq:return_map_stat}$\txtsub{2}$ holds true:
\begin{equation}
	\yfvct\!\rbr{\stsprm\trial} \leq \bzero,
\end{equation}
the loading is purely elastic and the solution is found:
\begin{equation}
	\stsprm = \stsprm\trial, \quad
	\pmltprm = \bzero.
\label{eq:return_map_elastic}	
\end{equation}
Contrarily, a plastic correction has to be computed on the basis of the trial stress~$\stsprm\trial$. Similarly to the problem of non-smooth multi-surface plasticity at material point level (e.g., see~\cite{Scalet_Auricchio_ARCME_2017}), the enforcement of the plastic admissibility conditions~\eqref{eq:return_map_stat}$\txtsub{2}$ requires explicitly introducing the set of active plastic multipliers:
\begin{equation}
	I\act = \cbr{i \in I \,\,\vert\,\, \sbr{\pmltprm}_i > 0}, 
\end{equation}
where $I = \cbr{1, \dots{}, N_{\pmlt}}$. On the other hand, as the set~$I\act$ is unknown and cannot be predicted by means of its trial counterpart:
\begin{equation}
	I\act\trial = \cbr{i \in I \,\,\vert\,\, \sbr{\yfvct\!\rbr{\stsprm\trial}}_i > 0},
\end{equation}
an iterative procedure has to be carried out~\cite{Simo_Govindjee_IJNME_1988}. Denoting by $k = 0, 1, \dots{}$ the last iteration index, equations~\eqref{eq:return_map_stat} are solved assuming the active set to be~$I\act^k$ (initialized by~$I\act^0 = I\act\trial$):
\begin{align}
	\begin{aligned}
		\bzero &= \cmplmtx\e\!\rbr{\stsprm^{k+1} - \stsprm\trial} + \flowmtx\!\rbr{\stsprm^{k+1}}\pmltprm^{k+1}, \\[1ex]
		0 &= \sbr{\yfvct\!\rbr{\stsprm^{k+1}}}_i,	\quad i \in I\act^{k}, \\[1ex]
		0 &= \sbr{\pmltprm^{k+1}}_i,	\quad i \in I \setminus I\act^{k},
	\end{aligned}
\label{eq:return_map_system}	
\end{align}
thus resulting in a nonlinear system of dimension~$N_{\sts}+N_{\pmlt}$ in the unknown stresses~$\stsprm^{k+1}$ and plastic multipliers~$\pmltprm^{k+1}$. Provided equations~\eqref{eq:return_map_system}$\txtsub{3}$ are substituted in the remaining ones, such nonlinear system reduces to dimension~$N_{\sts}+N^{k}_{\pmlt}$, with $N^{k}_{\pmlt} = \#{I\act^k}$, in the unknown stresses~$\stsprm^{k+1}$ and active plastic multipliers~$\pmltprm\act^{k+1}$:
\begin{align}
	\begin{aligned}
		\bzero &= \res_{\sts}^{k+1} &\hspace{-0.4cm}&= \cmplmtx\e\!\rbr{\stsprm^{k+1} - \stsprm\trial} + \flowmtx\act\!\rbr{\stsprm^{k+1}}\pmltprm\act^{k+1}, \\[1ex]
		\bzero &= \res_{\pmlt}^{k+1} &\hspace{-0.4cm}&= \yfvct\act\!\rbr{\stsprm^{k+1}}, \\[1ex]
	\end{aligned}
\label{eq:return_map_algorithm}
\end{align}
in which $\yfvct\act$ [resp., $\flowmtx\act$] is the subvector [resp., submatrix] of $\yfvct$ [resp., $\flowmtx$] identified by~$I\act^k$. Adopting Newton's method of solution, with $j$ as the iteration index, a sequence of approximations of the solution is derived by the linearized system:
\begin{align}
	\begin{aligned}
		&\res_{\sts}^{k+1, j+1} \approx \res_{\sts}^{k+1, j} \\[1ex] 
			&\hspace{1cm}+\cmplmtx\ep\!\rbr{\stsprm^{k+1, j}, \pmltprm^{k+1, j}}\rbr{\stsprm^{k+1, j+1} - \stsprm^{k+1, j}} \\[1ex]
			&\hspace{1cm}+\flowmtx\act\!\rbr{\stsprm^{k+1, j}}\rbr{\pmltprm\act^{k+1, j+1}-\pmltprm\act^{k+1, j}} = \bzero, \\[1ex]
		&\res_{\pmlt}^{k+1, j+1} \approx \res_{\pmlt}^{k+1, j} \\[1ex] 
			&\hspace{1cm}+\flowmtx\act\txtsup{T}\!\rbr{\stsprm^{k+1, j}}\rbr{\stsprm^{k+1, j+1} - \stsprm^{k+1, j}}  = \bzero.
	\end{aligned}
\label{eq:return_map_Newton}
\end{align}
Here, the elastoplastic compliance matrix has been introduced as:
\begin{equation}
	\cmplmtx\ep\!\rbr{\stsprm, \pmltprm} = \cmplmtx\e + \pdiff{\flowmtx}{\stsprm}\!\rbr{\stsprm}\!\sbr{\pmltprm}, 
\label{eq:return_map_complmtx}
\end{equation}
with the derivative of the element plastic flow operator contracted with the plastic multiplier (component-wise) given by:
\begin{multline}
	\sbr{ \pdiff{\flowmtx}{\stsprm}\!\rbr{\stsprm}\!\sbr{\pmltprm} }_{il} \\[1ex] = \rbr{
		\intbr{\sbr{\stsmtx}_{pi} 
		\sbr{\partial^{2}_{\sts\sts}\yf\!\rbr{\stsmtx \stsprm}}_{pq} 
		\sbr{\stsmtx}_{ql} \sbr{\pmltmtx}_{k} } } \sbr{\pmltprm}_k,
\end{multline}
summation over repeated indices understood. 
Once convergence has been reached within a prescribed tolerance in Newton's iterations, the relevant solution, i.e.~the stress and plastic multiplier parameters, $\stsprm^{k+1}$ and $\pmltprm^{k+1}$ respectively, are checked to fulfill the plastic admissibility conditions:
\begin{align}
	\begin{aligned}
		&\sbr{\pmltprm^{k+1}}_i \geq \bzero, \,\, i \in I\act^{k}, \\[1ex]
		&\sbr{\yfvct\!\rbr{\stsprm^{k+1}}}_i \leq \bzero, \,\, i \in \cbr{1, \dots{}, N_{\pmlt}} \setminus I\act^{k}.
	\end{aligned}
\label{eq:return_map_act_check}
\end{align}
If such is the case, the solution has been found. Contrarily, the set of active plastic multipliers is updated in:
\begin{align}
	\begin{aligned}
		&I\act^{k+1} = \cbr{i \in I \setminus I\act^{k} \,\,\vert\,\, \sbr{\yfvct\!\rbr{\stsprm^{k+1}}}_i > 0} \\[1ex]
		&\hspace{3cm}\bigcup
		\cbr{i \in I\act^{k} \,\,\vert\,\, \sbr{\pmltprm^{k+1}}_i > 0},
	\end{aligned}
\label{eq:return_map_update_act}
\end{align}
and next iteration is performed.

It is observed that the present element return mapping procedure couples the plastic behavior of all quadrature points to each other, irrespectively of the particular interpolation adopted for the plastic multiplier~\cite{Simo_Taylor_CMAME_1989}. In fact, even for a pointwise imposition of the plastic admissibility condition (Section~\ref{sss:plastic_mlt_pointwise}), if yielding occurs at any quadrature point, the entire element flows plastically to fulfill the element compatibility.

\begin{remark}
In~\cite{Simo_Govindjee_IJNME_1988} a condition different from equation \eqref{eq:return_map_update_act} is given for the update of the set of active plastic multipliers, namely:
\begin{equation}
	I\act^{k+1} = \cbr{i \in I\act^{k} \,\,\vert\,\, \sbr{\pmltprm^k}_i > 0}, 
\end{equation}
whence it follows that $I\act^{k} \subset I\act^{k-1} \subset \dots{} \subset I\act^{0} = I\act\trial$.
However, a simple counterexample, showing that the final set of active constraints may not be a subset of the active constraints determined by the predictor, is presented in~\cite{Bilotta_Garcea_FEAD_2011}. In fact, a component of the element yield function that was non-negative at the previous iteration might become strictly positive at the current one. Such a possibility is accounted for by the first set in equation~\eqref{eq:return_map_update_act}.
\end{remark}

\begin{remark}\label{r:constant_pmlt}
The present elastic predictor/plastic corrector algorithm results to be greatly simplified in the case of constant interpolation of the plastic multiplier field (Remark~\ref{r:pltmlt_constant}), as discussed in~\cite{Bilotta_Casciaro_CMAME_2007} for a material with von Mises yield criterion. In fact, referring to the discrete problem~\eqref{eq:plastic_admissibility_Bilotta_2}, the element yield function $\yfvct$ boils down to the scalar~$\overline\yf$ and the character of non-smooth multi-surface plasticity vanishes. Consequently, if plastic loading occurs, the algorithm reduces to Newton's method iterations corresponding to equations \eqref{eq:return_map_Newton}.
\end{remark}

\subsubsection{Consistent tangent stiffness matrix}\label{sss:return_map_tangent}
After the stress and plastic multiplier parameters, $\stsprm$ and $\pmltprm$ respectively, have been determined, the element forces~$\force\internal$ are computed from the equilibrium equation~\eqref{eq:return_map_q_int}. Accordingly, the consistent element stiffness operator~$\stiff$ is given by:
\begin{equation}
	\stiff 
	= \compmtx\txtsup{T}	\, \pdiff{\stsprm}{\displprm}, 
\label{eq:return_map_stiff}
\end{equation}
and requires the computation of the derivative of the stress parameters~$\stsprm$ with respect to the nodal DOFs~$\displprm$. In case of purely elastic loading, from equations~\eqref{eq:return_map_elastic}$\txtsub{1}$ and~\eqref{eq:trial_stsprm}, it results that:
\begin{equation}
	\stiff = \compmtx\txtsup{T} \, \rbr{\cmplmtx\e}^{-1} \, \compmtx.
\end{equation} 
Conversely, if element plastic flow occurs, differentiation of equations~\eqref{eq:return_map_algorithm} yields:
\begin{align}
	\begin{aligned}
		\bzero &= 
			\cmplmtx\ep\!\rbr{\stsprm, \pmltprm}\,\pdiff{\stsprm}{\displprm} - \compmtx+ 
			\flowmtx\act\!\rbr{\stsprm}\pdiff{\pmltprm\act}{\displprm}, \\[1ex]
		\bzero & = \flowmtx\act\txtsup{T}\!\rbr{\stsprm} \pdiff{\stsprm}{\displprm},
	\end{aligned}
\label{eq:return_map_Dini}
\end{align}
where equation~\eqref{eq:trial_stsprm} has been used. By solving equation~\eqref{eq:return_map_Dini}$\txtsub{1}$ for $\partial{\stsprm}/\partial{\displprm}$ and substituting the result in equation~\eqref{eq:return_map_Dini}$\txtsub{2}$, it follows that:
\begin{equation}
	\pdiff{\pmltprm\act}{\displprm} = \bW\act^{-1} \, \flowmtx\act\txtsup{T} \, \rbr{\cmplmtx\ep}^{-1} \, \compmtx,
\end{equation}
with the following position:
\begin{equation}
	\bW\act = \flowmtx\act\txtsup{T} \, \rbr{\cmplmtx\ep}^{-1} \, \flowmtx\act.
\end{equation}
Finally, from equation~\eqref{eq:return_map_Dini}$\txtsub{1}$ the derivative of the stress parameters~$\stsprm$ with respect to the nodal DOFs~$\displprm$ turns out to be:
\begin{equation}
	\pdiff{\stsprm}{\displprm} = \rbr{\cmplmtx\ep}^{-1} \sbr{ \bI_{N_{\sts}} -
		\flowmtx\act \, \bW\act^{-1} \, \flowmtx\act\txtsup{T} \, {\cmplmtx\ep}^{-1}}
		\compmtx,
\end{equation}
with $\bI_{N_{\sts}}$ as the identity matrix of dimension~$N_{\sts}$, whence the consistent element stiffness matrix can be derived from equation~\eqref{eq:return_map_stiff}:
\begin{equation}
	\stiff = \compmtx\txtsup{T}	\, \rbr{\cmplmtx\ep}^{-1} \sbr{ \bI_{N_{\sts}} - 
		\flowmtx\act \, \bW\act^{-1} \, \flowmtx\act\txtsup{T} \, {\cmplmtx\ep}^{-1} }
		\compmtx.
\end{equation}

\subsection{Mathematical programming}
An alternative strategy to the element return mapping algorithm for the solution of the CM formulation equations is derived by observing that compatibility and plastic admissibility conditions~\eqref{eq:return_map_stat} can be recast in the following nonlinear optimization problem:
\begin{equation}
	\begin{aligned}
		& \underset{\stsprm}{\text{maximize}}
		& & - \frac{1}{2} \rbr{\stsprm-\stsprm\trial}\txtsup{T} \cmplmtx\e \rbr{\stsprm-\stsprm\trial} \\[0.5ex]
		& \text{subject to}
		& & \yfvct\!\rbr{\stsprm} \leq \bzero.
	\end{aligned}
\label{eq:return_map_OPT}
\end{equation}
As corresponding to the convex projection of the trial stress~$\stsprm\trial$ onto the elastic domain given by~$\yfvct\!\rbr{\stsprm} = \bzero$, such formulation complies with standard closest-point-projection schemes adopted for the elastoplastic material state update. In the following two optimization techniques, i.e.~interior point algorithms and sequential quadratic programming, are discussed as solution method of problem~\eqref{eq:return_map_OPT}.

\subsubsection{Interior-point algorithms}\label{sss:IPM}
Interior-point (IP) methods tackle the nonlinear optimization problem~\eqref{eq:return_map_OPT} by (i) converting the inequality constraints into equalities, provided a slack variable has been added, and (ii) introducing a penalty term into the objective function for the slack variable to be strictly positively restricted~\cite{Krabbenhoft_Wriggers_IJNME_2007}. Accordingly, the modified problem turns out to be:
\begin{equation}
	\begin{aligned}
		& \underset{\cbr{\stsprm, \bs}}{\text{maximize}}
		& & \Big\{- \frac{1}{2} \rbr{\stsprm-\stsprm\trial}\txtsup{T} \cmplmtx\e \rbr{\stsprm-\stsprm\trial} \\
			&&&\hspace{3.2cm}+ \mu \textstyle\sum_{i=1}^{N_{\pmlt}} \log{s_i} \Big\} \\[1ex]
		& \text{subject to}
		& & \yfvct\!\rbr{\stsprm} + \bs = \bzero,
	\end{aligned}	
\end{equation}
where $\bs$ is the slack variable, $\mu$ is an arbitrarily small positive constant and $\mu \log{s_i}$ is the so-called logarithmic barrier function. The associated Lagrangian reads:
\begin{align}
	\begin{aligned}
		\cL\!\rbr{\stsprm, \pmltprm, \bs} &= - \frac{1}{2} \rbr{\stsprm-\stsprm\trial}\txtsup{T} \cmplmtx\e \rbr{\stsprm-\stsprm\trial} \\[1ex]
		&-\pmltprm\txtsup{T} \sbr{\yfvct\!\rbr{\stsprm} + \bs} +
		 \mu \textstyle\sum_{i=1}^{N_{\pmlt}} \log{s_i}, 
	\end{aligned}
\end{align}
and its optimality conditions are given by:
\begin{align}
	\begin{aligned}
		\bzero &= \res_{\sts} = \cmplmtx\e\rbr{\stsprm-\stsprm\trial} + \flowmtx\!\rbr{\stsprm}\pmltprm, \\[1ex]
		\bzero &= \res_{\pmlt} =  \yfvct\!\rbr{\stsprm} + \bs, \\[1ex]
		\bzero &= r_{\bs} = \mu \bS^{-1} \boldsymbol{1} - \pmltprm,
	\end{aligned}
\label{eq:return_map_IP}
\end{align}
where $\bS = \diag{\bs}$ and~$\boldsymbol{1}$ is a vector of ones. It is remarked that equations~\eqref{eq:return_map_IP}$\txtsub{2--3}$ imply the plastic admissibility conditions~\eqref{eq:return_map_stat}$\txtsub{2}$ up to small positive constants. In fact, as $\bs > \bzero$ for the presence of the logarithmic barrier and $\mu > 0$, also $\pmltprm > \bzero$; moreover, $\pmltprm\txtsup{T}\yfvct\!\rbr{\stsprm} = - N_{\pmlt} \mu$. Within the standard primal-dual IP method, Newton's iterations are applied to the optimality conditions~\eqref{eq:return_map_IP}. Specifically, at $j$-th iteration the following linearized system has to be solved:
\begin{align}
	\begin{aligned}
	&\rbr{\begin{array}{c} 
		\res_{\sts}^{j+1} \\[1ex]
		\res_{\pmlt}^{j+1} \\[1ex] 
		r_{\bs}^{j+1} \\[1ex]
	\end{array}} \approx
	\rbr{\begin{array}{c} 
		\res_{\sts}^{j} \\[1ex]
		\res_{\pmlt}^{j} \\[1ex] 
		r_{\bs}^{j} \\[1ex]
	\end{array}} + 	\\[1ex]
	&\sbr{\begin{array}{ccc} 
		\!\!\cmplmtx\ep\!\rbr{\stsprm^{j}, \pmltprm^{j}} & \flowmtx\!\rbr{\stsprm^j} & \bzero \\[1ex]
		\flowmtx\txtsup{T}\!\rbr{\stsprm^j} & \bzero & \bI_{N_{\Delta\lambda}} \\[1ex]
		\bzero & \bI_{N_{\Delta\lambda}} & \rbr{\bS^j}^{-2}\!\!
	\end{array}}\!\!\!
	\rbr{\begin{array}{c} 
		\!\!\stsprm^{j+1}-\stsprm^j \!\!\\[1ex]
		\!\!\pmltprm^{j+1}-\pmltprm^j \!\!\\[1ex] 
		\!\!\bs^{j+1} - \bs^j \!\!
	\end{array}} \!=\!
	\rbr{\begin{array}{c} 
		\bzero \\[1ex]
		\bzero \\[1ex] 
		\bzero \\[1ex]
	\end{array}}\!,
	\end{aligned}		
\end{align}
where $\bI_{N_{\Delta\lambda}}$ is the identity matrix of dimension~$N_{\Delta\lambda}$ and the elastoplastic compliance matrix $\cmplmtx\ep\!\rbr{\stsprm, \pmltprm}$ has been introduced in equation~\eqref{eq:return_map_complmtx}. IP methods take their name from the limitation imposed on the increments~$\rbr{\pmltprm^{j+1}-\pmltprm^j}$ and~$\rbr{\bs_i^{j+1} - \bs_i^j}$ to derive strictly positive updated approximations. Accordingly, a maximum allowable step length~$\omega$ is introduced by:
\begin{align}
	\begin{aligned}
		&\omega = \max\cbr{\omega_{\pmlt}, \, \omega_{\bs}}, \\[1ex]
		&\omega_{\pmlt} = \max_{\Delta l_i^{j+1} < \Delta l_i^{j}} \cbr{-\frac{\Delta l_i^{j}}{\Delta l_i^{j+1} - \Delta l_i^{j}}}, \\[1ex]
		&\omega_{\bs} = \max_{s_i^{j+1} < s_i^{j}} \cbr{-\frac{s_i^{j}}{s_i^{j+1} - s_i^{j}}},
	\end{aligned}
\end{align}
in conjunction with the update rules:
\begin{align}
	\begin{aligned}
		&\stsprm^{j+1} = \stsprm^j + \vartheta\omega\rbr{\stsprm^{j+1}-\stsprm^j}, \\[1ex]
		&\pmltprm^{j+1} = \pmltprm^j + \vartheta\omega\rbr{\pmltprm^{j+1}-\pmltprm^j}, \\[1ex]
		&\bs^{j+1} = \bs^j + \vartheta\omega\rbr{\bs^{j+1}-\bs^j},
	\end{aligned}
\end{align}
where $\vartheta$ is a relaxation parameter usually chosen to be close to unity, e.g.~$\vartheta\approx 0.95$--$0.9995$~\cite{Krabbenhoft_Wriggers_IJNME_2007}. An important aspect for the performances of IP methods is the application, as the algorithm proceeds, of a convenient reduction rule to the barrier parameter~$\mu$. In fact, fast reduction may induce a loss of convergence, whereas slow reduction may induce a large number of iterations. In~\cite{Krabbenhoft_Wriggers_IJNME_2007}, the following reduction rule is considered:
\begin{equation}
	\mu^{j+1} = \eta \frac{\delta^j}{N_{\pmlt}}, \quad
	\delta^j = \rbr{\bs^j}\txtsup{T} \pmlt^j,
\end{equation}
with $\eta \approx 0.2$--$0.7$. Nevertheless, the authors suggest a simple modification of the IP methods based on the assumptions that $\mu = 0$, $\vartheta = \omega = 1$ and after each iteration $\pmltprm^{j+1}$ and $\bs^{j+1}$ are corrected as:
\begin{align}
	\begin{aligned}
		&\Delta l_i^{j+1} = \max\cbr{\Delta l_i^{j+1}, \varepsilon_{\pmlt}}, \\[1ex]
		&s_i^{j+1} = \max\cbr{s_i^{j+1}, \varepsilon_{\bs}},
	\end{aligned}
\end{align}
with $\varepsilon_{\pmlt} \approx 10^{-9}$--$10^{-12}$ and $\varepsilon_{\bs} \approx 10^{3} \varepsilon_{\pmlt}$. The consistent tangent stiffness matrix can be finally derived as discussed in Section~\ref{sss:return_map_tangent}.

\begin{remark}\label{r:Krabbenhoft_Wriggers}
As proposed in~\cite{Krabbenhoft_Wriggers_IJNME_2007}, the IP algorithm at hand can be directly applied to the solution of structural equilibrium equations, i.e.~of the system comprising the local nonlinear optimization problems~\eqref{eq:return_map_OPT} and the global equilibrium equations (obtained from the element internal forces given in equation~\eqref{eq:return_map_q_int} by standard assembly procedure). Accordingly, the following nonlinear optimization problem is considered at structural level:
\begin{equation}
	\begin{aligned}
		& \min_{\displprm} \max_{\stsprm\elem}
		& & \Big\{\displprm\txtsup{T}\sbr{\cA_{e=1}^{N\elem} \,\compmtx\elem\txtsup{T}\stsprm\elem - \force\external} \\
			&&&-\textstyle\sum_{e=1}^{N\elem} 
			\sbr{\frac{1}{2} \, \stsprm\elem\txtsup{T} \, \cmplmtx\e\elem \, \stsprm\elem +
			\stsprm\elem\txtsup{T} \rbr{\stnprm\p\n}\elem} \Big\} \\[0.5ex]
		& \text{subject to}
		& & \yfvct\!\rbr{\stsprm\elem} \leq \bzero, \quad e = 1, \dots{}, N\elem,
	\end{aligned}	
\label{eq:return_map_OPT_structural}
\end{equation}
and its IP-modified version results to be:
\begin{equation}
	\begin{aligned}
		&\min_{\displprm} \max_{\cbr{\stsprm\elem, \bs\elem}}
		& & \Big\{\displprm\txtsup{T}\sbr{\cA_{e=1}^{N\elem} \, \compmtx\elem\txtsup{T}\stsprm\elem - \force\external} \\
			&&&- \textstyle\sum_{e=1}^{N\elem} 
			\sbr{\frac{1}{2} \, \stsprm\elem\txtsup{T} \, \cmplmtx\e\elem \, \stsprm\elem +\stsprm\elem\txtsup{T} \rbr{\stnprm\p\n}\elem} \\
			&&& + \mu \textstyle\sum_{i=1}^{N_{\pmlt}} \log{\rbr{\bs\elem}_i} \Big\}  \\[0.5ex]
		&\text{subject to}
		& & \yfvct\!\rbr{\stsprm\elem} + \bs\elem = \bzero, \quad e = 1, \dots{}, N\elem.
	\end{aligned}	
\end{equation}
As claimed by the authors~\cite{Krabbenhoft_Wriggers_IJNME_2007}, the solution of such system of equations by Newton's method yields a novel optimization point of view in elastoplasticity. In fact, the tradition of satisfying exactly (up to prescribed tolerance) the local constitutive equations at the end of each structural iteration is broken in favor of the rule-of-thumb, well estabilished in optimization community, that is desirable to have all residuals converging at comparable rates.
\end{remark}

\subsubsection{Sequential quadratic programming}
The sequential quadratic programming (SQP) approach aims at constructing a sequence of approximations of the solution of the nonlinear optimization problem~\eqref{eq:return_map_OPT}, denoted~$\rbr{\stsprm^{j}, \pmltprm^{j}}$, by iterative linearization of its constraints. Apart from inessential constants in the objective functions, the following quadratic programming (QP) problem is considered for a given iterative solution~$\stsprm^{j}$:
\begin{equation}
	\begin{aligned}
		& \underset{\Delta\stsprm}{\text{maximize}}
		& & - \frac{1}{2} \Delta\stsprm\txtsup{T} \cmplmtx\e \Delta\stsprm - \Delta\stsprm\txtsup{T} \cmplmtx\e\rbr{\stsprm^{j}-\stsprm\trial} \\[0.5ex]
		& \text{subject to}
		& & \yfvct\!\rbr{\stsprm^{j}} + \flowmtx\txtsup{T}\!\rbr{\stsprm^j}\Delta\stsprm \leq \bzero,
	\end{aligned}
\label{eq:return_map_SQP}	
\end{equation}
whence the updated solution follows as~$\stsprm^{j+1} = \stsprm^{j} + \Delta\stsprm$. In~\cite{Bilotta_Garcea_FEAD_2011}, it is suggested to solve such QP problem by using the Goldfarb-Idnani dual active set method~\cite{Goldfarb_Idnani_MP_1983}, which is highly efficient and, at no extra cost, furnishes the algorithmic tangent operator relative to each iteration. Moreover, a line search algorithm is employed to enforce global convergence of the algorithm \cite{Nocedal_Wright_2006}. 

\begin{remark}
With a point of view similar to the one in Remark~\ref{r:Krabbenhoft_Wriggers}, the solution of the structural nonlinear optimization problem~\eqref{eq:return_map_OPT_structural} through the so-called equality-constraint sequential quadratic programming (EC-SQP) is proposed in~\cite{Bilotta_Garcea_CS_2012}. The algorithm consists in two stages: (i) to estimate the active inequality constraints and (ii) to solve the equality constrained quadra\-tic program resulting by considering only the predicted active constraints. In particular, the estimation of the set of active constraints is derived by a single iteration of the SQP in the local optimization problems~\eqref{eq:return_map_SQP}. Consequently, a single Newton's method iteration is performed at structural level to obtain the updated approximation of the solution. The convergence is achieved in case relevant residuals do not exceed prescribed tolerances and conditions~\eqref{eq:return_map_act_check}, concerning the fulfillment of constraint equations, are met.

Some observations are in order for such an algorithm. First, a noteworthy aspect of the procedure is that, because of their linearity, the equilibrium equations at element level are exactly satisfied at each iteration. Accordingly, at each iteration the algorithm tries to satisfy the plastic admissibility condition by searching the solution among the equilibrated ones. Second, it may be noted that performing the stage (i) of the algorithm aims at deriving an intermediate improved prediction of the set of active constraint between two consecutive Newton's iterations. In~\cite{Bilotta_Garcea_CS_2012}, the authors state that an improved estimate of the set of active constraints, as obtained by performing more than one iteration in stage~(i), is not significant for the robustness of the algorithm. On the other hand, an approach that can be regarded as a particular instance of the algorithm under consideration, though not presented in the framework of nonlinear optimization, is proposed in~\cite{Mendes_Castro_FEAD_2009} and prescribes the complete elimination of the stage (i).
\end{remark}

\subsection{Nodal-force-based algorithm}\label{ss:element_iem}
When considering a HW formulation, the need of an accurate description of the strain field might lead to a quite large dimension of the strain interpolation space. Consequently, a direct application of Newton's method might results not satisfactorily robust for the solution of the element state determination problem~\eqref{eq:HW_fun_stat}, in which quite many coupled nonlinear and nonsmooth (constitutive) equations are involved. In this section, the iterative procedure proposed in~\cite{Nodargi_Bisegna_C&S_2017} for the case of piecewise constant strain interpolation (Section~\ref{sss:piecewise_stn}) is reviewed.

\subsubsection{Filtering out rigid body motions}\label{ss:corot}
The following developments presuppose the mapping of nodal DOFs into internal force vector, $\displprm \mapsto \bq\txtsup{int}$,  to be invertible. Because rigid body motions do not contribute to incremental energy, that is not the case. Nevertheless, a homomorphism $\overline{\bPi}$, whose kernel is composed by rigid body motions, can be chosen into the space of nodal DOFs $\displprm$:
\begin{equation}
	\fltdisplprm = \overline{\bPi} \displprm.
\label{eq:overline_P_map}
\end{equation}
Such operator behaves as a filter, mapping the nodal DOFs $\displprm$ into their purely deformational counterpart $\fltdisplprm$. Accordingly, a congruence relation, whose equivalence classes are nodal DOFs~$\displprm$ differing by a rigid body motion, is introduced. The canonical map onto the induced quotient space and the typical equivalence class are respectively denoted by $\bPi$ and $\frdisplprm$:
\begin{equation}
	\frdisplprm = \bPi \displprm.
\label{eq:Pi_map}	
\end{equation}
By the first isomorphism theorem, the quotient space is isomorphic to the image of $\overline{\bPi}$, whence there exists an isomorphism $\bV$ such that:
\begin{equation}
	\fltdisplprm = \bV \frdisplprm.
\label{eq:V_map}	
\end{equation}
Denoting by~$N\txtsub{r}$ the number of element rigid body modes, it follows that $\fltdisplprm$ is a $N_{\displ}$-components vector, whereas $\frdisplprm$ can be regarded as a $\rbr{N_{\displ}-N\txtsub{r}}$-component vector of deformational parameters. 

Exploiting such a framework, the nodal DOFs are replaced by their filtered counterparts (correspondence law $\displprm \rightarrow \fltdisplprm$) in stationary conditions~\eqref{eq:HW_fun_stat_piecewise_stn}, and correspondingly equation~\eqref{eq:HW_q_int} is intended to yield the element forces work-conjugated to the filtered nodal DOFs~$\fltdisplprm$, to be denoted by $\overline\bq\txtsup{int}$. The element forces $\bq\txtsup{int}$ and $\bq\txtsup{int}_{\blambda}$, respectively work-conjugated to the nodal DOFs $\displprm$ and to the deformational parameters $\frdisplprm$, can be then derived by energy equivalence:
\begin{equation}
	\var \displprm\txtsup{T} \bq\txtsup{int} = \var \frdisplprm\txtsup{T} \bq\txtsup{int}_{\blambda} = \var \fltdisplprm\txtsup{T} \overline\bq\txtsup{int},
\label{eq:energy_equivalence}
\end{equation}
whence, exploiting equations~\eqref{eq:Pi_map} and \eqref{eq:V_map}, it follows that:
\begin{equation}
	\bq\txtsup{int} = \proj\txtsup{T} \bq\txtsup{int}_{\blambda},
    \quad
	\bq\txtsup{int}_{\blambda} = \lift\txtsup{T} \overline\bq\txtsup{int}.
\label{eq:filter_q}
\end{equation}
It is pointed out that, differently from the mappings $\fltdisplprm \mapsto \overline\bq\txtsup{int}$ and $\displprm \mapsto \bq\txtsup{int}$, the mapping $\frdisplprm \mapsto \bq\txtsup{int}_{\blambda}$ turns out to be invertible, and this will be instrumental in Section~\ref{ss:iterative_solution_algorithm}.

\subsubsection{Iterative solution algorithm}\label{ss:iterative_solution_algorithm}
Upon filtering out rigid body motions by the technique discussed in the previous section, the stationary conditions of HW formulation with piecewise-constant strain interpolation~\eqref{eq:HW_fun_stat_piecewise_stn} become:
\begin{align}
	\begin{aligned}
		\bzero &= {\pdiffarg{\stn}{\ie} \! \rbr{\stnprm\qp} -\avgstsmtx\qp\stsprm} \,, \\[1ex]
		\bzero &= \textstyle\sum_{\qppnt=1}^{N\qp}{\avgstsmtx\qp^{\text{T}} \stnprm\qp \abs{\dmn\qp}} - \compmtx \lift \frdisplprm.
	\end{aligned}
\label{eq:elem_res_fem}
\end{align}
In the standard displacement-driven architecture of finite element computer programs, the nodal DOFs $\displprm$, and hence the deformational parameters $\frdisplprm$ descending from equation~\eqref{eq:Pi_map}, are known at element level. Let $\frdisplprm\txtsub{s}$ denote that actual value, where the index ``$\text{s}$'' stands for ``solution''. Accordingly, the system of $\rbr{N_\sigma + N_{\varepsilon}}$ equations~\eqref{eq:elem_res_fem} can be solved with respect to stress and strain parameters, respectively $\stsprm$ and $\stnprm\qp, \, \qppnt = 1, \dots{}, N\qp$. Unfortunately, because compatibility condition~\eqref{eq:elem_res_fem}$\txtsub{2}$ couples the unknown strains involved in the nonlinear and nonsmooth constitutive equations~\eqref{eq:elem_res_fem}$\txtsub{1}$, a direct application of Newton's method cannot guarantee satisfactory robustness. 

Basic idea of the nodal-force-based algorithm is to interpret the deformational parameters $\frdisplprm$ as depending on their work-conjugated nodal forces $\bq\txtsup{int}_{\blambda}$, where the latter descend from equations~\eqref{eq:HW_q_int} and \eqref{eq:filter_q}$\txtsub{2}$:
\begin{equation}
	\bq\txtsup{int}_{\blambda} = \rbr{\compmtx \lift} \txtsup{T} \stsprm.
\label{eq:equilibrium}
\end{equation}
Consequently, the following nodal compatibility condition has to be imposed:
\begin{equation}
	\res\!\rbr{\bq\txtsup{int}_{\blambda}} = \frdisplprm\!\rbr{\bq\txtsup{int}_{\blambda}} - \frdisplprm\txtsub{s} = \bzero,
\label{eq:nd_comp}
\end{equation}
the mapping $\bq\txtsup{int}_{\blambda} \mapsto \frdisplprm$ being defined through equations~\eqref{eq:elem_res_fem} and \eqref{eq:equilibrium}. It is remarked that such a definition is well-posed because the $\rbr{N_{\displ}-N\txtsub{r}} \times \rbr{N_{\displ}-N\txtsub{r}}$ matrix $\compmtx \lift$ is invertible. Thus, by means of element equilibrium~\eqref{eq:equilibrium}, the stress parameters~$\stsprm$ can be computed as:
\begin{equation}
	\stsprm = \rbr{\compmtx \lift}^{-\text{T}} \bq\txtsup{int}_{\blambda}.
\label{eq:el_eq}
\end{equation}
Then, the element subdomain constitutive equations \eqref{eq:elem_res_fem}$\txtsub{1}$ can be solved for the unknown strains $\stnprm\qp$ by inversion of material constitutive law. Here, the assumption of strictly convex stress potential $\ie$ is exploited. Recalling equation~\eqref{eq:characterization}, the inverse constitutive law is expressed in terms of the incremental complementary energy~$\ie\conj$, whence:
\begin{equation}
	\stnprm\qp = \partial_{\sts}\ie\conj\!\rbr{\avgstsmtx\qp\stsprm} \,,
\label{eq:el_cl}
\end{equation}
Finally, internal compatibility equation~\eqref{eq:elem_res_fem}$\txtsub{2}$ yields the deformational parameters~$\frdisplprm$:
\begin{equation}
	\frdisplprm = \rbr{\compmtx \lift}^{-1} \textstyle\sum_{\qppnt=1}^{N\qp}{\avgstsmtx\qp^{\text{T}} \stnprm\qp \abs{\dmn\qp}}.
\label{eq:el_comp}
\end{equation}
The key advantage of the present algorithm is the possibility to solve the constitutive equations at any element subdomain independently from each other, thus mitigating convergence difficulties and reducing computational cost. On the other hand, the same numerical techniques used for the enforcement of direct material model can be used for its inversion, at no further computational cost (see Remark~\ref{r:complementary_ie} and~\cite{Nodargi_Bisegna_C&S_2017}).

Equation~\eqref{eq:nd_comp} is finally solved by adopting Newton's method. Denoting by $j$ the iteration index and exploiting the implicit function theorem, the following sequence of approximations for the solution $\bq\txtsup{int}_{\blambda}$ is derived:
\begin{equation}
	\rbr{\bq\txtsup{int}_{\blambda}}^{j+1} = \rbr{\bq\txtsup{int}_{\blambda}}^{j} - \sbr{\pdiff{\bq\txtsup{int}_{\blambda}}{\frdisplprm}\rbr{\frdisplprm^{j}}} \rbr{\frdisplprm^{j} - \frdisplprm\txtsub{s}},
\label{eq:el_Newton}
\end{equation}
where $\rbr{\bq\txtsup{int}_{\blambda}}^{0}$ is assumed to be the relevant quantity at the end of the previous Newton structural iteration. After convergence has been achieved through such approximating sequence, the element forces $\bq\txtsup{int}$ are computed from equation~\eqref{eq:filter_q}$\txtsub{1}$.

\subsubsection{Consistent tangent stiffness matrix}
For the computation of the derivative of element nodal forces $\bq\txtsup{int}_{\blambda}$ with respect to deformational parameters $\frdisplprm$, element equilibrium~\eqref{eq:el_eq} yields:
\begin{equation}
	\pdiff{\bq\txtsup{int}_{\blambda}}{\frdisplprm} = \rbr{\compmtx \lift}\txtsup{T} \pdiff{\stsprm}{\frdisplprm},
\label{eq:el_eq_diff}
\end{equation}
and requires the computation of the derivative of stress interpolation parameters $\stsprm$. To this end, element subdomain constitutive equations~\eqref{eq:el_cl} give:
\begin{equation}
	\pdiff{\stnprm\qp}{\frdisplprm} = \partial^2_{\sts\sts}\ie\conj\!\rbr{\avgstsmtx\qp\stsprm} \avgstsmtx\qp \pdiff{\stsprm}{\frdisplprm}.
\end{equation}
On the other hand, internal compatibility equations~\eqref{eq:el_comp} imply the following linear relationships between the derivatives of strains $\stnprm\qp$ at element subdomains:
\begin{equation}
	\bI_{\rbr{N_{\displ}-N\txtsub{r}}} = \rbr{\compmtx \lift}^{-1} \sum_{\qppnt=1}^{N\qp}{\avgstsmtx\qp^{\text{T}} \pdiff{\stnprm\qp}{\frdisplprm} \abs{\dmn\qp}},
\label{eq:el_comp_diff}
\end{equation}
where $\bI_{\rbr{N_{\displ}-N\txtsub{r}}}$ is the identity matrix of dimension $\rbr{N_{\displ}-N\txtsub{r}}$. Using equations~\eqref{eq:el_eq_diff}--\eqref{eq:el_comp_diff}, it is a simple matter to check that the derivative of element nodal forces $\bq\txtsup{int}_{\blambda}$ turns out to be:
\begin{equation}
	\pdiff{\bq\txtsup{int}_{\blambda}}{\frdisplprm} = \rbr{\compmtx \lift}\txtsup{T} \elflex^{-1} \rbr{\compmtx \lift},	
\label{eq:q_lambda_diff}
\end{equation}
where $\elflex$ is the element compliance operator:
\begin{equation}
	\cmplmtx = \textstyle\sum_{\qppnt=1}^{N\qp}{\avgstsmtx\qp^{\text{T}} \, \partial^2_{\sts\sts}\ie\conj\!\rbr{\avgstsmtx\qp\stsprm} \avgstsmtx\qp \abs{\dmn\qp}}.
\end{equation}

Finally, the consistent element stiffness operator $\elstiff$, i.e. the derivative of the element forces $\bq\txtsup{int}$ with respect to nodal DOFs~$\displprm$ follows from conditions~\eqref{eq:Pi_map}, \eqref{eq:filter_q} and \eqref{eq:q_lambda_diff}:
\begin{equation}
	\elstiff = \proj\txtsup{T} \, \pdiff{\bq\txtsup{int}_{\blambda}}{\frdisplprm} \, \proj.
\end{equation}

\section{Numerical stability}\label{s:stability}
For a mixed formulation to be numerically stable, interpolation spaces of displacement and additional fields need to enjoy the compatibility requirements expressed by the ellipticity and the inf--sup conditions of Babu\v{s}ka and Brezzi \cite{Bathe_1996, Boffi_Brezzi_Fortin_2013}. 

Especially the strictness of the latter has motivated the twofold focus of researchers to modify mixed formulations in order to circumvent stability issues or to explore the stability of newly proposed elements. Among stabilization techniques proposed to achieve stable elements irrespective of the inf--sup condition, it is worth mentioning Galerkin least square stabilization \cite{Hughes_Franca_CMAME_1987}, bubble enrichment of displacement field \cite{Arnold_Fortin_Calcolo_1984, Mahnken_Laschet_CMAME_2008, Caylak_Mahnken_CS_2014} and variational multi-scale formulations \cite{Cervera_Codina_I_CMAME_2010, Cervera_Codina_III_CMAME_2015, Cervera_Chiumenti_CM_2016}. On the other hand, an analytical proof of mixed element stability would be an indisputable achievement. However, apart from few cases (for instance, see \cite{Chapelle_Bathe_CS_1993}), such a result is usually difficult to be accomplished. Accordingly, interest has also been devoted to the derivation of a reliable numerical test to draw a prediction of element stability, aiming to be a mixed-method counterpart of the patch test used in standard displacement based formulations (for instance, see \cite{Taylor_Chan_IJNME_1986}). Bathe and co-authors have proposed a test which amounts at the solution of a generalized eigenvalue problem involving the discrete operators relevant to the variational formulation under investigation~\cite{Chapelle_Bathe_CS_1993, Bathe_CS_2001}. Exploiting the usual displacement-driven structure of the numerical solution strategy for the structural discrete problem (Section~\ref{s:numerical_strategy}), a similar numerical test has been presented and validated in~\cite{Nodargi_Bisegna_IJNME_2016, Nodargi_Bisegna_C&S_2017}, exhibiting reliability, simplicity of implementation and limited computational cost as its main advantages.

In this section, the general algebraic structure of mixed formulation problems is first discussed in an abstract setting. Particular implications are then derived for the functionals presented in Section~\ref{s:gsm_variational}. Specifically, as the numerical stability problem is addressed referring to a linear elastic material, the analysis is restricted to HW, ES and HR functionals. The numerical test for checking mixed element stability proposed in~\cite{Nodargi_Bisegna_IJNME_2016, Nodargi_Bisegna_C&S_2017} is finally reviewed.

\subsection{Stability conditions} \label{ss:abstract_setting}
The stationary conditions of the mixed functionals investigated in~Section~\ref{s:gsm_variational} can be set within a unitary analytical setting, which is here presented. Let $X$ and $Y$ denote two Hilbert spaces. Upon introducing the two continuous bilinear forms:
\begin{equation}
	a: X \times X \rightarrow \mathbb{R}, \quad 
	b: X \times Y \rightarrow \mathbb{R}, \quad
\end{equation}
and given two continuous linear functionals $\bbf \in X\conj$ and $\bg \in Y\conj$, with $X\conj$ [resp., $Y\conj$] as the (topological) dual space of $X$ [resp., $Y$], the following problem is considered:
\begin{align}
	\begin{aligned}		 
		& a\!\rbr{\bbx, \var\bbx} + b\!\rbr{\var\bbx, \by} = \duality{\bbf}{\var\bbx}{X\conj}{X}, \quad &&\forall \,\var\bbx \in X, \\[0.5ex]
		& b\!\rbr{\bbx, \var\by} = \duality{\bg}{\var\by}{Y\conj}{Y}, \quad &&\forall \, \var\by \in Y.
	\end{aligned}
\label{eq:continuum1}
\end{align}
in the unknown $\rbr{\bbx, \by} \in X\times Y$. In particular, by the Riesz representation theorem, the bilinear forms $a$ and $b$ can be also represented in terms of two continuous linear operators:
\begin{align}
	\begin{aligned}
		&\mathbb{A} : X \rightarrow X\conj, && 
		\bduality{\mathbb{A}\bbx_1}{ \bbx_2}_X = a \rbr{\bbx_1, \bbx_2},\,\, \forall \, \bbx_1, \bbx_2 \in X, \\[0.5ex]
		&\mathbb{B} : X \rightarrow Y\conj, &&
		\bduality{\mathbb{B \bbx}}{\by}_{Y} = b\!\rbr{\bbx, \by},\,\, \forall \,\bbx \in X, \,\, \forall \,\by \in Y,
	\end{aligned}
\end{align}
whence problem~\eqref{eq:continuum1} is equivalent to:
\begin{align}
	\begin{aligned}
		&\mathbb{A}\bbx + \mathbb{B}\conj\by = \bbf, \\[0.5ex]
		&\mathbb{B}\bbx = \bg.
	\end{aligned}
\label{eq:continuum2}	
\end{align} 
Here $\mathbb{B}\conj:Y\rightarrow X\conj$ denotes the adjoint operator of $\mathbb{B}$, i.e.~the (unique) operator such that $\bduality{\mathbb{B \bbx}}{\by}_{Y} = \bduality{\bbx}{\mathbb{B}\conj \by}_{X}$, $\forall \,\bbx \in X, \,\,\forall \,\by \in Y$. Henceforth, the two formulations \eqref{eq:continuum1} and \eqref{eq:continuum2} are considered the same and one or the other is used according to the convenience of the moment.

With the position $K = \ker{\mathbb{B}}$, sufficient conditions for problem~\eqref{eq:continuum1} to have a unique solution for any $\bbf \in X\conj$ and $\bg \in Y\conj$ are: 
\begin{align}
	\begin{aligned}
		&\exists \,\, \alpha > 0 \bdot a\!\rbr{\bbx, \bbx} \geq \alpha \norm{\bbx}_{X}^2, \quad \forall \, \bbx \in K, \\[1ex]
		&\exists \,\, \beta > 0 \bdot \adjustlimits \inf_{\by \in Y} \sup_{\bbx \in X}
			\frac{b\!\rbr{\bbx, \by}}{\norm{\bbx}_{X}\norm{\by}_Y } \geq \beta.
	\end{aligned}
\label{eq:existence_cond}	
\end{align}
Condition~\eqref{eq:existence_cond}$\txtsub{1}$ requires the bilinear form $a$ to be coercive on $K$. Conversely, condition~\eqref{eq:existence_cond}$\txtsub{2}$ implies the operator $\mathbb{B}\conj$ to be bounding, and hence injective and with closed range. A detailed proof of this result can be found in~\cite{Boffi_Brezzi_Fortin_2013}. 

To derive a mixed finite element approximation of problem~\eqref{eq:continuum1}, an increasing sequence of (finite-dimen\-sional) interpolation spaces $X\interp \subset X$ and $Y\interp \subset Y$, with~$h \rightarrow 0$ as a mesh parameter, is introduced for the unknown fields. Accordingly, $\rbr{\bbx\interp, \by\interp} \in X\interp \times Y\interp$ is the discrete solution of problem~\eqref{eq:continuum1} in case:
\begin{align}
	\begin{aligned}
		& a\!\rbr{\bbx\interp, \var\bbx\interp} + b\!\rbr{\var\bbx\interp, \by\interp} = \bduality{\bbf}{\var\bbx\interp}, \quad &&\forall \, \var\bbx\interp \in X\interp, \\[1ex]
		& b\!\rbr{\bbx\interp, \var\by\interp} = \bduality{\bg}{\var\by\interp}, \quad &&\forall \, \var\by\interp \in Y\interp,
	\end{aligned}
\label{eq:discrete1}
\end{align}
or, equivalently, by using the same notation for the discrete operators~$\mathbb{A}\interp$ and $\mathbb{B}\interp$ and their matrix representations:
\begin{equation}
	\sbr{\begin{array}{cc} 
		\mathbb{A}\interp & \rbr{\mathbb{B}\interp}\txtsup{\!T} \\[1ex]
		\mathbb{B}\interp & \bzero \\[1ex]
	\end{array}}
	\rbr{\begin{array}{c} 
		\hat{\bbx}\interp \\[1ex]
		\hat{\by}\interp
	\end{array}} =
	\rbr{\begin{array}{c} 
		\bbf\interp \\[1ex]
		\bg\interp \\[1ex] 
	\end{array}},
\label{eq:discrete2}	
\end{equation}
where $\hat{\bbx}\interp$ [resp., $\hat{\by}\interp$] denotes interpolation parameters in $X\interp$ [resp., $Y\interp$], and $\bbf\interp$ [resp., $\bg\interp$] is the finite-dimen\-sional counterpart of~$\bbf$ [resp., $\bg$].

For the finite dimensional problem~\eqref{eq:discrete1}, conditions analogous to equations~\eqref{eq:existence_cond} guarantee the existence of a unique solution:
\begin{align}
	\begin{aligned}
		&\exists \,\, \alpha\interp > 0 \bdot a\!\rbr{\bbx\interp, \bbx\interp} \geq \alpha\interp \norm{\bbx\interp}_{X}^2, \quad \forall \bbx\interp \in K\interp, \\[1ex]
		&\exists \,\, \beta\interp > 0 \bdot \adjustlimits \inf_{\by\interp \in Y\interp} \sup_{\bbx \in X\interp}
			\frac{b\!\rbr{\bbx\interp, \by\interp}}{\norm{\bbx\interp}_{X}\norm{\by\interp}_Y} \geq \beta\interp,
	\end{aligned}
\label{eq:existence_discrete}	
\end{align}
where $K\interp = \ker{\rbr{\mathbb{B}\interp}}$ and constants $\alpha\interp$ and $\beta\interp$ in general depend on $h$. Condition~\eqref{eq:existence_discrete}$\txtsub{1}$, implying the operator $\mathbb{A}\interp$ to be coercive on $K\interp$, guarantees that the sub-matrix containing~$\mathbb{A}\interp$ and~$\mathbb{B}\interp$ is injective. Similarly, condition~\eqref{eq:existence_discrete}$\txtsub{2}$ implies that $\rbr{\mathbb{B}\interp}\txtsup{\!T}$ is injective. Consequently, the overall discrete operator is injective as well, and hence invertible, i.e.~system~\eqref{eq:discrete2} is solvable.

Finally, the sequence of problems~\eqref{eq:discrete1} is said numerically stable for $h \rightarrow 0$ if the two constants $\alpha\interp$ and $\beta\interp$ involved in conditions~\eqref{eq:existence_discrete} are bounded from below by strictly positive constants $\alpha$ and $\beta$:
\begin{align}
	\begin{aligned}
		&\exists \,\, \alpha > 0 \bdot a\!\rbr{\bbx\interp, \bbx\interp} \geq \alpha \norm{\bbx\interp}_{X}^2, \quad \forall \bbx\interp \in K\interp, \\[1ex]
		&\exists \,\, \beta > 0 \bdot \adjustlimits \inf_{\by\interp \in Y\interp} \sup_{\bbx\interp \in X\interp}
			\frac{b\!\rbr{\bbx\interp, \by\interp}}{\norm{\bbx\interp}_{X}\norm{\by\interp}_Y} \geq \beta.
	\end{aligned}
\label{eq:stability_discrete}	
\end{align}
In other words, the problem is numerically stable if it is uniformly solvable with respect to the mesh parameter~$h$.

\subsubsection{Hu--Washizu functional}
Referring to the HW formulation discussed in Section~\ref{ss:HW}, the functional spaces to be considered are $X = E \times V$ and $Y = S$, with the bilinear forms $a$ and $b$ defined by: 
\begin{align}
	\begin{aligned}
		&a\!\rbr{\rbr{\stn, \displ}\!, \rbr{\bbeta, \bv}} = \intbr{\fC\e \stn \cdot \bbeta \,}, \\[1ex]
		&b\!\rbr{\rbr{\stn, \displ}\!, \sts} = \intbr{-\sts \cdot \rbr{\stn - \sym\nabla\displ}}.
	\end{aligned}
\end{align}
Upon introducing interpolation spaces $X\interp = E\interp \times V\interp$ and $Y\interp = S\interp$, their discrete representation reads:
\begin{equation}
	\mathbb{A}\interp = \sbr{\begin{array}{cc} 
		\rbr{\bC\e}\interp & \bzero \\[1ex]
		\bzero & \bzero  
	\end{array}}, \quad
	\mathbb{B}\interp = \sbr{\begin{array}{cc} 
		-\rbr{\stsstnmtx\interp}\txtsup{\!T} & \compmtx\interp
	\end{array}},
\end{equation}
with $\compmtx\interp$ and~$\stsstnmtx\interp$ introduced in equations~\eqref{eq:comp_stsstn_operator} (superscript~$h$ is here added to remark their dependence on the mesh parameter) and with~$\rbr{\bC\e}\interp$ as the elastic energy matrix given by:
\begin{equation}
	\rbr{\bC\e}\interp = \intbr{\stnmtx\tp \fC\e \stnmtx \,},
\end{equation}
whence the element discrete problem turns out to be:
\begin{equation}
	\sbr{\begin{array}{ccc} 
		\rbr{\bC\e}\interp & \bzero & -\stsstnmtx\interp \\[1ex]
		\bzero & \bzero & \rbr{\compmtx\interp\large}\txtsup{\!T} \\[1ex]
		-\rbr{\stsstnmtx\interp}\txtsup{\!T} & \compmtx\interp & \bzero 
	\end{array}}
	\rbr{\begin{array}{c} 
		\stnprm\interp \\[1ex]
		\displprm\interp \\[1ex] 
		\stsprm\interp 
	\end{array}} =
	\rbr{\begin{array}{c} 
		\bzero \\[1ex]
		\rbr{\force\external}\interp \\[1ex] 
		\bzero
	\end{array}}.
\label{eq:HW_discrete}	
\end{equation}
Accordingly, the discrete solvability conditions~\eqref{eq:existence_discrete} result to be:
\begin{align}
	\begin{aligned}
		&\exists \, \alpha\interp \!>\! 0 \bdot a ((\stn\interp, \displ\interp), (\stn\interp, \displ\interp)) \geq \alpha\interp \norm{(\stn\interp, \displ\interp)}_{E\times V}^2, \\
		&\hspace{6cm} \forall (\stn\interp, \displ\interp) \in K\interp, \\[1ex]
		&\exists \,  \beta\interp \!>\! 0 \bdot \!\!\adjustlimits \inf_{\sts\interp \in S\interp} \sup_{\rbr{\stn\interp, \displ\interp} \in E\interp \times V\interp}
			\!\frac{b\!\rbr{\rbr{\stn\interp, \displ\interp}\!, \sts\interp}}{\norm{\sts\interp}_S \norm{\rbr{\stn\interp, \displ\interp}}_{E\times V}} \geq \beta\interp\!,
	\end{aligned}
\label{eq:HW_existence_discrete}	
\end{align}
where:
\begin{equation} 
	K\interp = \cbr{(\stn\interp,\displ\interp) \in E\interp \times V\interp \,\,\vert\, -\rbr{\stsstnmtx\interp}\txtsup{\!T}\stnprm\interp + \compmtx\interp\displprm\interp = \bzero}
\end{equation}
represents the subset of compatible interpolated strains and displacements. Provided the elastic energy matrix $\rbr{\bC\e}\interp$ is positive definite, that is a stable material is considered, necessary and sufficient condition for the operator $\mathbb{A}\interp$ to be coercive on $K\interp$, required by condition~\eqref{eq:existence_discrete}$\txtsub{1}$, is the injectivity of $\compmtx\interp$. Accordingly, there cannot exist any spurious mode, i.e.~non-trivial DOFs~$\displprm\interp$ such that~$\compmtx\interp\displprm\interp = \bzero$. On the other hand, the injectivity of $\rbr{\mathbb{B}\interp}\txtsup{\!T}$, required by condition~\eqref{eq:existence_discrete}$\txtsub{2}$, rules out the possibility that interpolated stresses $\stsstnmtx\interp\stsprm\interp$ and nodal internal forces~$\rbr{\compmtx\interp\large}\txtsup{\!T}\stsprm\interp$ simultaneously vanish for non-vanishing stress parameters~$\stsprm\interp$. Such conditions can be then recast in:
\begin{align}
	\begin{aligned}
		&\ker{\compmtx\interp} = \cbr{\bzero}, \\[0.5ex]
		&\ker{\rbr{\mathbb{B}\interp}\txtsup{\!T}} = \ker{\stsstnmtx\interp} \cap \ker{\rbr{\compmtx\interp}\txtsup{\!T}} = \cbr{\bzero}.
	\end{aligned}
\label{eq:HW_discrete_2}
\end{align}
Uniform bounds ensuring the validity of such conditions with respect to the mesh parameter $h\rightarrow 0$ yield the numerical stability of the formulation.

\begin{remark}
Necessary conditions for equations~\eqref{eq:HW_discrete_2} to hold are:
\begin{equation}
	N_{\displ} - N_{r} \leq N_{\sts} \leq N_{\stn}+N_{\displ} - N_{r},
\end{equation}
where it is recalled that $N_{r}$ denotes the number of element rigid body modes.
\end{remark}

\subsubsection{Enhanced strain functional}
The ES functional introduced in Section~\ref{ss:enhanced_strain} can be investigated in the light of the abstract setting discussed in Section~\ref{ss:abstract_setting} by assuming the spaces $X = V \times \tilde E$ and $Y = S$, and the bilinear forms:
\begin{align}
	\begin{aligned}
		&a\!\rbr{\rbr{\displ, \tilde\stn}\!, \rbr{\bv,\tilde\bbeta}} = \intbr{\fC\e \rbr{\sym\nabla\displ + \tilde\stn} \cdot \rbr{\sym\nabla\bv + \tilde\bbeta} \,}, \\[1ex]
		&b\!\rbr{\rbr{\displ,\tilde\stn}\!, \sts} = \intbr{-\sts \cdot \tilde\stn\,}.
	\end{aligned}
\end{align}
Selecting interpolation spaces~$X\interp = V\interp \times \tilde E\interp$ and~$Y\interp = S\interp$, the discrete representations of $a$ and $b$ are:
\begin{equation}
	\mathbb{A}\interp \!=\! \sbr{\begin{array}{cc} 
		\stiff_{\displ}\interp & \rbr{\stiff_{\tilde\stn,\displ}\interp}\txtsup{\!T}\! \\[1ex]
		\!\stiff_{\tilde\stn,\displ}\interp & \stiff_{\tilde\stn}\interp  
	\end{array}}\!\!, \quad
	\mathbb{B}\interp \!=\! \sbr{\begin{array}{cc} 
		\bzero & -\!\rbr{\tilde\stsstnmtx\interp}\txtsup{\!T}\!
	\end{array}}\!,
\end{equation}
where~$\stiff_{\displ}\interp$, $\stiff_{\tilde\stn,\displ}\interp$ and~$\stiff_{\tilde\stn}\interp$ are the linear elastic counterpart of the matrices introduced in equations~\eqref{eq:enh_K_u_K_eps_u} and~\eqref{eq:enh_K_eps}, respectively, and define the following discrete element problem:
\begin{equation}
	\sbr{\begin{array}{ccc} 
		\stiff_{\displ}\interp & \rbr{\stiff_{\tilde\stn,\displ}\interp}\txtsup{\!T} & \bzero \\[1ex]
		\stiff_{\tilde\stn,\displ}\interp & \stiff_{\tilde\stn}\interp & -\tilde\stsstnmtx\interp \\[1ex]
		\bzero & -\rbr{\tilde\stsstnmtx\interp}\txtsup{\!T} & \bzero 
	\end{array}}
	\rbr{\begin{array}{c} 
		\displprm\interp \\[1ex] 
		\tilde\stnprm\interp \\[1ex]		
		\stsprm\interp 
	\end{array}} =
	\rbr{\begin{array}{c}
		\rbr{\force\external}\interp \\[1ex] 	 
		\bzero \\[1ex]
		\bzero
	\end{array}}.
\label{eq:enh_discrete}	
\end{equation}
The relevant solvability conditions result to be:
\begin{align}
	\begin{aligned}
		&\exists \, \alpha\interp \!>\! 0 \bdot a ((\displ\interp,\tilde\stn\interp), (\displ\interp,\tilde\stn\interp)) \geq \alpha\interp \norm{(\displ\interp,\tilde\stn\interp)}_{V\times \tilde E}^2, \\
		&\hspace{6cm} \forall (\displ\interp,\tilde\stn\interp) \in K\interp, \\[1ex]
		&\exists \, \beta\interp \!>\! 0 \bdot \!\!\adjustlimits \inf_{\sts\interp \in S\interp} \sup_{\rbr{\stn\interp, \displ\interp} \in E\interp \times V\interp}
			\!\frac{b\,((\displ\interp,\tilde\stn\interp), \sts\interp)}{\norm{\sts\interp}_S \norm{(\displ\interp,\tilde\stn\interp)}_{V\times \tilde E}} \!\geq\!\beta\interp,
	\end{aligned}
\label{eq:enh_existence_discrete}	
\end{align}
where:
\begin{equation}
	K\interp = \cbr{(\displprm\interp,\tilde\stnprm\interp) \in V\interp \times \tilde E\interp \,\,\vert\, \rbr{\tilde\stsstnmtx\interp}\txtsup{\!T}\tilde\stnprm\interp = \bzero}
\end{equation}
represents the set of admissible enhanced strains. The coercivity condition of~$\mathbb{A}\interp$ on $K\interp$, equation~\eqref{eq:enh_existence_discrete}$\txtsub{1}$, is satisfied if a stable material is considered. For the operator~$\rbr{\mathbb{B}\interp}\txtsup{\!T}$ to be injective, equation~\eqref{eq:enh_existence_discrete}$\txtsub{2}$, it is required that~$\tilde\stsstnmtx\interp$ is injective, hence precluding the existence of spurious enhanced strain modes. The discrete solvability conditions can be recast as:
\begin{equation}
	\ker{\tilde\stsstnmtx\interp} = \cbr{\bzero},
\label{eq:enh_existence_discrete_2}		
\end{equation}
and numerical stability of the formulation follows by its uniform fulfillment for the mesh parameter~$h\rightarrow 0$.

\begin{remark}
Necessary condition for satisfying equations \eqref{eq:enh_existence_discrete_2} is that:
\begin{equation}
	N_{\sts} \leq N_{\tilde\stn}.
\end{equation}
\end{remark}

\begin{remark}
In case stress elimination is performed in the ES formulation by selecting $L^{2}$-ortho\-gonal interpolations of enhanced strain and stress fields, condition~\eqref{eq:enh_existence_discrete_2} becomes immaterial and the formulation numerical stability is guaranteed by energy coercivity as in displacement-based formulations. However, the issue of stress recovery requires to be solved. 
\end{remark}

\subsubsection{Hallinger--Reissner functional}
The HR formulation investigated in Section~\ref{ss:HR} complies with the abstract setting of Section~\ref{ss:abstract_setting} by assuming~$X = S$, $Y = V$ and the bilinear forms~$a$ and~$b$ to be given by:
\begin{align}
	\begin{aligned}
		&a\!\rbr{\sts, \btau} = -\intbr{{\fC\e}^{-1} \sts \cdot \btau \,}, \\[1ex]
		&b\!\rbr{\sts, \displ} = \intbr{\sts \cdot \sym\nabla\displ \,}.
	\end{aligned}
\end{align}
Consequently, in the discrete framework corresponding to selection of interpolation spaces~$X\interp = S\interp$ and $Y\interp = V\interp$, the following discrete operators represent $a$ and $b$:
\begin{equation}
	\mathbb{A}\interp = -\rbr{\cmplmtx\e}\interp, \quad
	\mathbb{B}\interp = \rbr{\compmtx\interp\large}\txtsup{\!T},
\end{equation}
where~$\rbr{\cmplmtx\e}\interp$ and~$\compmtx\interp$ have been introduced in equations~\eqref{eq:return_map_H} and~\eqref{eq:comp_stsstn_operator}$\txtsub{1}$, respectively, and the discrete element problem results in:
\begin{equation}
	\sbr{\begin{array}{cc} 
		-\rbr{\cmplmtx\e}\interp & \compmtx\interp \\[1ex]
		\rbr{\compmtx\interp\large}\txtsup{\!T} & \bzero\\[1ex]
	\end{array}}
	\rbr{\begin{array}{c} 
		\stsprm\interp \\[1ex]
		\displprm\interp
	\end{array}} =
	\rbr{\begin{array}{c} 
		\bzero \\[1ex]
		\rbr{\force\external}\interp \\[1ex] 
	\end{array}}.
\label{eq:HR_discrete}	
\end{equation}
The discrete solvability conditions turn out to be:
\begin{align}
	\begin{aligned}
		&\exists \,\, \alpha\interp > 0 \bdot a\!\rbr{\sts\interp, \sts\interp} \geq \alpha\interp \norm{\sts\interp}_{S}^2, \quad \forall \sts\interp \in K\interp, \\[1ex]
		&\exists \,\, \beta\interp > 0 \bdot \adjustlimits \inf_{\displ\interp \in V\interp} \sup_{\sts\interp \in S\interp}
			\frac{b\!\rbr{\sts\interp, \displ\interp}}{\norm{\sts\interp}_{S}\norm{\displ\interp}_V} \geq \beta\interp,
	\end{aligned}
\label{eq:HR_existence_discrete}	
\end{align}
where:
\begin{equation}
	K\interp = \cbr{\stsprm\interp \in S\interp \,\,\vert\, \rbr{\compmtx\interp\large}\txtsup{\!T}\stsprm\interp = \bzero}
\end{equation}
is the set of self-equilibrated stress parameters. In case a stable material is considered, the coercivity condition of~$\mathbb{A}\interp$ on $K\interp$, equation~\eqref{eq:HR_existence_discrete}$\txtsub{1}$, is automatically satisfied. Contrarily, the operator~$\rbr{\mathbb{B}\interp}\txtsup{\!T}$ to be injective, equation~\eqref{eq:HR_existence_discrete}$\txtsub{2}$, requires $\compmtx\interp$ to be injective, and hence that no spurious mode can exist. Accordingly, the discrete solvability conditions can be rephrased in:
\begin{equation}
	\ker{\compmtx\interp} = \cbr{\bzero},
\label{eq:HR_existence_discrete_2}	
\end{equation}
and its validity uniformly with respect to the mesh parameter~$h\rightarrow 0$ ensures numerical stability of the formulation.

\begin{remark}
Necessary condition for equations~\eqref{eq:HR_existence_discrete_2} to hold is:
\begin{equation}
	N_{\displ} - N_{r} \leq N_{\sts},
\end{equation}
with $N_{r}$ as the number of element rigid body modes.
\end{remark}

\subsection{Numerical tests}
In this section, the numerical test proposed in~\cite{Nodargi_Bisegna_IJNME_2016, Nodargi_Bisegna_C&S_2017} for checking the stability of mixed formulations is reviewed. Basic assumption is that a displacement-driven architecture is adopted in the numerical solution strategy of the discrete structural problem. Accordingly, it can be assumed that all additional fields with respect to displacements have been statically condensed out at element level in the element state determination procedure and the only equilibrium equations have to be solved at structural level. In particular, those equations take the form:
\begin{equation}
	\bK\interp \bu\interp = \rbr{\forces\external}\interp \,,
	\label{eq:FEM_problem}
\end{equation}
where $\bK\interp$ is the structural stiffness matrix, obtained by assembling element contributions derived in Section~\ref{s:numerical_strategy} for the variational formulations under consideration, and $\rbr{\forces\external}\interp$ is the external force vector. 

The structural stiffness matrix~$\bK\interp$ can be regarded as the representation of a linear operator from the displacement interpolation space~$V\interp$, which is a subset of the displacement ambient space~$V$, into itself. For a fixed mesh size $h$, necessary and sufficient conditions for it to be an isomorphism are (e.g., see \cite{Boffi_Brezzi_Fortin_2013}):
\begin{align}
	\begin{aligned}
		&\exists \,\, \beta\interp > 0 \,\, \bdot \,\, \norm{\bK\interp\bv\interp}_{V} \geq \beta\interp \norm{\bv\interp}_{V}\,, \quad \forall \, \bv\interp \in V\,, \\[1ex]
		&\norm{{\bK\interp}\txtsup{T} \bv\interp}_{V} > 0\,, \quad \forall \bv\interp \in V \backslash \cbr{\bzero}\,.
	\end{aligned}
\label{eq:K_bounding_Kt_injective}
\end{align}

The bounding condition \eqref{eq:K_bounding_Kt_injective}$\txtsub{1}$ is customarily recast in the following inf-sup form:
\begin{equation}
	\exists \,\, \beta\interp > 0 \,\, : \,\, \adjustlimits \inf_{\bv\interp \in V\interp} \sup_{\bw\interp \in V\interp}
	\frac{\bduality{\bK\interp\bv\interp}{\bw\interp}_{V}}{\norm{\bv\interp}_{V} \norm{\bw\interp}_{V}} \geq \beta\interp \,,
\label{eq:K_infsup}
\end{equation}
which yields a useful characterization of the constant $\beta\interp$. Indeed, introducing the $V$-norm operator:
\begin{equation}
	\bT : V \to {V}\dual\,, \quad \duality{\bT \bv}{\bv}{{V}\dual}{V}= \norm{\bv}_{V}^2\,,
\end{equation}
the inf-sup constant $\beta\interp$ is proven to be the minimum eigenvalue of the generalized eigenvalue problem \cite{Bathe_1996, Boffi_Brezzi_Fortin_2013}:
\begin{equation}
	\bK\interp \bv\interp = \lambda\interp \bT\bv\interp\,.
\label{eq:eigenvalue}
\end{equation}

\begin{remark}
A simple proof can be obtained arguing as in~\cite{Bathe_1996}. Let $\lambda_1\interp\leq \dots{} \leq \lambda_N\interp$ be the eigenvalues of problem \eqref{eq:eigenvalue} and $\cbr{\bphi_1\interp, \dots{}, \bphi_N\interp}$ the corresponding eigenvectors, supposed to form an orthonormal basis of $V\interp$. Exploiting the representation formulas $\bv\interp = v_i\interp \bphi_i\interp$ and $\bw\interp = w_i\interp \bphi_i\interp$, where summation is intended over repeated indices, the inf--sup in \eqref{eq:K_infsup} turns out to be:
\begin{align}
	\begin{aligned}
		&\adjustlimits 	\inf_{\bv\interp \in V\interp} \sup_{\bw\interp \in V\interp} \frac{\duality{\bK\interp\bv\interp}{\bw\interp}{{V\interp}\dual}{V\interp}}{\norm{\bv\interp}_{V\interp} \norm{\bw\interp}_{V\interp}} \\[1ex]
		&\hspace{1cm}= \inf_{\cbr{v_i\interp}} \frac{1}{\rbr{\sum_i {v_i\interp}^2}^{1/2}} \sup_{\cbr{w_i\interp}} \frac{\lambda_i\interp v_i\interp w_i\interp}{\rbr{\sum_i {w_i\interp}^2}^{1/2}}\,.
	\end{aligned}
\end{align}
Using the Cauchy--Schwarz inequality on the argument of the supremum and observing that the equality is attained for $\bv\interp = \bw\interp$, it follows that:
\begin{equation}
	\adjustlimits	\inf_{\bv\interp \in V\interp} \sup_{\bw\interp \in V\interp} \frac{\duality{\bK\interp\bv\interp}{\bw\interp}{{V\interp}\dual}{V\interp}}{\norm{\bv\interp}_{V\interp} \norm{\bw\interp}_{V\interp}} =
				\inf_{\cbr{v_i\interp}} \frac{{\lambda_i\interp} {v_i\interp}^2}{\sum_i {v_i\interp}^2}\,.
\end{equation}
Finally, because the infimum is attained for $\bv = \bphi_1$, it results that:
\begin{equation}
	\adjustlimits \inf_{\bv\interp \in V\interp} \sup_{\bw\interp \in V\interp} 	\frac{\duality{\bK\interp\bv\interp}{\bw\interp}{{V\interp}\dual}{V\interp}}{\norm{\bv\interp}_{V\interp} \norm{\bw\interp}_{V\interp}} = \lambda_1\interp\,,
\label{eq:end_proof}
\end{equation}
which concludes the proof.
\end{remark}

Such a result ensures that the bounding condition \eqref{eq:K_bounding_Kt_injective}$\txtsub{1}$ is satisfied if and only if the minimum eigenvalue of problem \eqref{eq:eigenvalue}, i.e.~$\lambda_1\interp$, is strictly positive. Actually, it is a simple matter to check that the same requirement also implies the adjoint injectivity condition \eqref{eq:K_bounding_Kt_injective}$\txtsub{2}$. \\
Although the strictly positiveness of $\lambda_1\interp$ ensures the solvability of problem \eqref{eq:FEM_problem} for a fixed mesh size $h$, the finite element formulation is numerically stable provided that such solvability condition holds true uniformly for $h \to 0$. Hence, the stability condition can be stated as follows (e.g., see \cite{Boffi_Brezzi_Fortin_2013}):
\begin{equation}
	\exists \,\, \beta > 0 \,\, : \,\, \beta \norm{\bv\interp}_{V} \leq \norm{\bK\interp\bv\interp}_{V\dual}\,, \quad \forall \, \bv\interp \in V\interp\,, \quad \forall \, h\,,
\end{equation}
where the constant $\beta$ is independent from the mesh size $h$. Finally, the stability of a finite element can be explored by proving that, for $h \to 0$, the minimum eigenvalue of problem \eqref{eq:eigenvalue} is bounded from below by a strictly positive constant.

A noticeable advantage of the present numerical test is the easiness of implementation in a displacement-driven finite element computer code, as it does not require any additional computational effort but the assemblage of the matrix $\bT$ and the solution of the eigenvalue problem~\eqref{eq:eigenvalue}. It is remarked that the prediction depends on the choice of problem domain, mesh and boundary conditions, whence, by its own nature, the numerical test cannot replace analytical investigations on element stability. However, it might be used as a reliable preliminary tool in assessing the stability of a mixed formulation.

\section{Numerical simulations}\label{s:numerical_simulations}
In this section, the performances of several mixed finite element formulations available in the literature are investigated in numerical simulations. As intrinsic differences in the formulations preclude uniform comparison conditions, the present exposition aims at giving just a clue on the accuracy, robustness and efficiency of the mixed elements at hand. 

Specifically, one mixed element is considered for each of the variational frameworks discussed in Section~\ref{s:gsm_variational}. With ``Q4'' and ``Q8'' respectively referring to four-node and eight-node quadrilateral elements, i.e.~to bilinear and serendipity isoparametric interpolations of the displacement field (Sections~\ref{sss:bilinear} and~\ref{sss:serendipity}), the following elements are selected:
\begin{itemize}
	\item{HW-Q8-D~\cite{Nodargi_Bisegna_C&S_2017}, based on HW formulation (Section \ref{ss:HW}). A self-equilibrated stress interpolation is derived by Airy's function approach (Section \ref{sss:Airy}), comprising a complete basis of second order enriched by two cubic modes . Piecewise-constant strain interpolation is assumed at element level (Section \ref{sss:piecewise_stn}). The element state determination is performed by the nodal-force-based algorithm discussed in Section~\ref{ss:element_iem}};
	\item{ES-Q4~\cite{Simo_Rifai_IJNME_1990}, based on ES formulation (Section \ref{ss:enhanced_strain}). The enhanced strain field is approximated as discussed in Section~\ref{ss:enh_stn_interp}, whereas the stress field is eliminated from the discrete problem by $L^2$-orthogonal interpolation with respect to the enhanced strain. Newton's method is implemented for the element state determination (Section~\ref{sss:Newton_enh});}
	\item{CM-Q4~\cite{Simo_Taylor_CMAME_1989}, based on CM formulation (Section \ref{ss:return_map}). A stress approximation similar to the Pian-Sumihara 5-$\beta$ interpolation is assumed (Section \ref{sss:Pian_Sumihara}). The plastic multiplier is pointwise interpolated (Section \ref{sss:plastic_mlt_pointwise}). The element state determination is carried out by the element return mapping procedure reviewed in Section~\ref{ss:element_return_mapping};}		
	\item{HR-Q4~\cite{Schroder_Starke_CMAME_2017}, based on HR formulation (Section~\ref{ss:HR}). The Pian-Sumihara 5-$\beta$ stress interpolation is adopted (Section \ref{sss:Pian_Sumihara}). The element state determination is tackled by Newton's method (Section~\ref{sss:Newton_HR}).}\end{itemize}
For the sake of completeness, Q4 and Q8 displacement-based formulations are also analyzed. 

Two classical benchmarks in elastoplasticity, i.e.~the bending of a tapered beam referred to as Cook's membrane~\cite{Cook_Witt_1989} and the extension of a perforated plate (e.g., see~\cite{deSouza_Peric_2008}), are investigated. In both simulations, the material is assumed to be linearly elastic, with von Mises yield criterion, associative flow law and isotropic linear hardening (Section~\ref{s:elastoplastic_materials}). For convenience, it is recalled that the von Mises yield function is given by (for instance, see~\cite{Armero_ECM_2004}):
\begin{equation}
	\yf\!\rbr{\sts, q\iso} = \norm{\deviatoric\sts} - c\!\rbr{\sigma\y - q\iso},
\end{equation}
where~$\deviatoric$ denotes the operator extracting the deviatoric part of the argument, and $c$ and $\sigma\y$ are positive constants (assuming~$c=\sqrt{2/3}$, $\sigma\y$ can be interpreted as the initial yield limit in tension). Moreover, both the analyses are performed under plane stress assumption.

\subsection{Cook's membrane}
\begin{figure}
	\centering
	\includegraphics[trim = 0cm 0cm 0cm 0cm, clip = true, scale = 0.7]{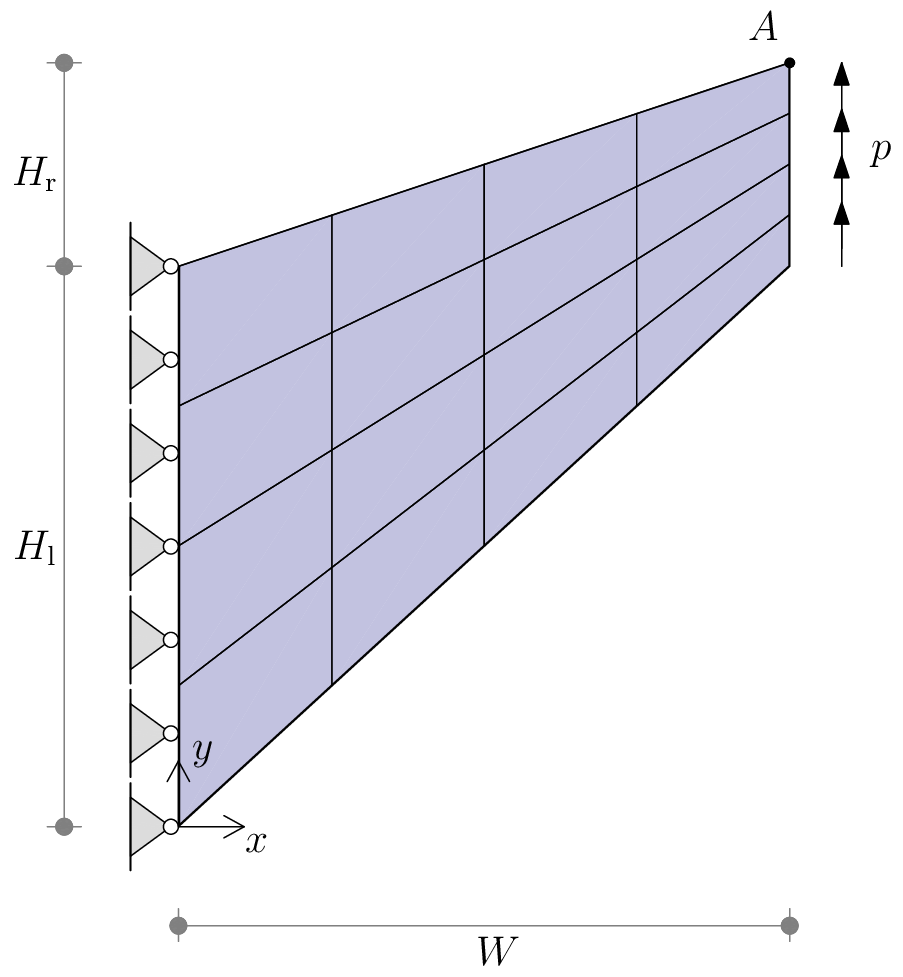}
	\caption{Cook's membrane: geometry, boundary conditions and finite element mesh with $n=4$ subdivisions.}
	\label{fig:cook_geometry}
\end{figure}
\begin{figure}
	\centering
	\includegraphics[trim = 0cm 0cm 0cm 0cm, clip = true, scale = 0.45]{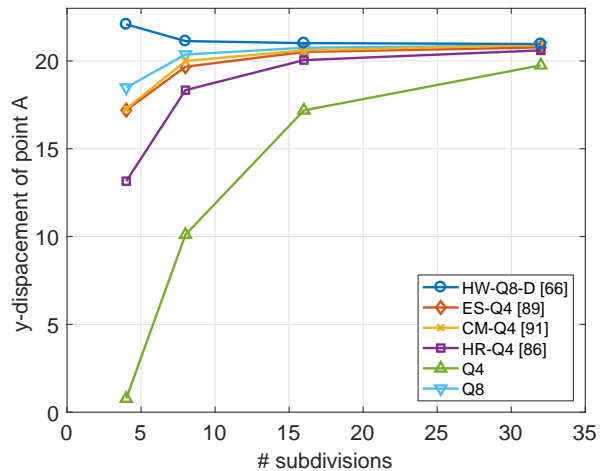}
	\caption{Cook's membrane: $h$-convergence analysis in terms of vertical displacement of point A.}
	\label{fig:Cook_conv}
\end{figure}
Cook's membrane problem is a popular benchmark to test element performances in bending dominated applications~\cite{Cook_Witt_1989}. Such membrane is characterized by the skewed shape shown in Figure~\ref{fig:cook_geometry}, with dimensions $W = 48$, $H\txtsub{l} = 44$, $H\txtsub{r} = 16$ and thickness $t = 1$. The left-end section is completely constrained, and a uniform tangential load is applied at free right edge. The applied uniform load magnitude is increased up to a resultant value $F\txtsub{max} = 1.8$. The material parameters adopted in the analysis are: Young's modulus $E = 70$, Poisson's ratio $\nu = 1/3$, yield stress $\sigma\y = 0.243$ and linear isotropic hardening modulus $k\iso = 0.2$. In order to explore $h$-convergence properties of the formulations under investigation, the problem is solved adopting finite element meshes corresponding to different level of refinement. In particular, each mesh is defined by the number~$n$ of subdivisions of the membrane sides. The values $n = 4,8,16, 32$ are taken into account (in Figure~\ref{fig:cook_geometry}, the finite element mesh with $n = 4$ is depicted). 

In Figure~\ref{fig:Cook_conv} the $h$-convergence analysis for the vertical displacement of point~$A$ is reported. Poor performances of Q4 displacement-based formulation can be highlighted. Among mixed four-node quadrilaterals, the ES and CM formulations respectively benefit from an improved description of the strain field and of the plastic strain field (via the plastic multiplier), and
converge faster than the HR formulation. HW quadrilateral is practically at convergence already for the mesh with~$n = 8$ subdivisions, in that taking advantage of the very accurate and flexible strain field approximation within each finite element.

\subsection{Extension of a perforated plate}
\begin{figure}
	\centering
	\includegraphics[trim = 0cm 0cm 0cm 0cm, clip = true, scale = 0.84]{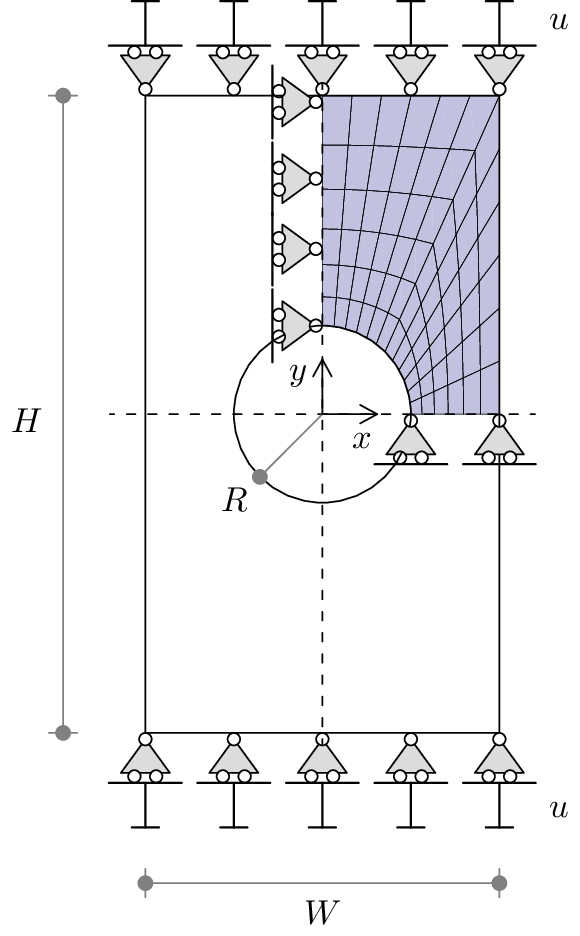}
	\caption{Extension of a perforated plate: geometry, boundary conditions and finite element mesh $6 \times 12$.}
	\label{fig:plate_geometry}
\end{figure}
\begin{figure}
	\centering
	\includegraphics[clip = true, scale = 0.55]{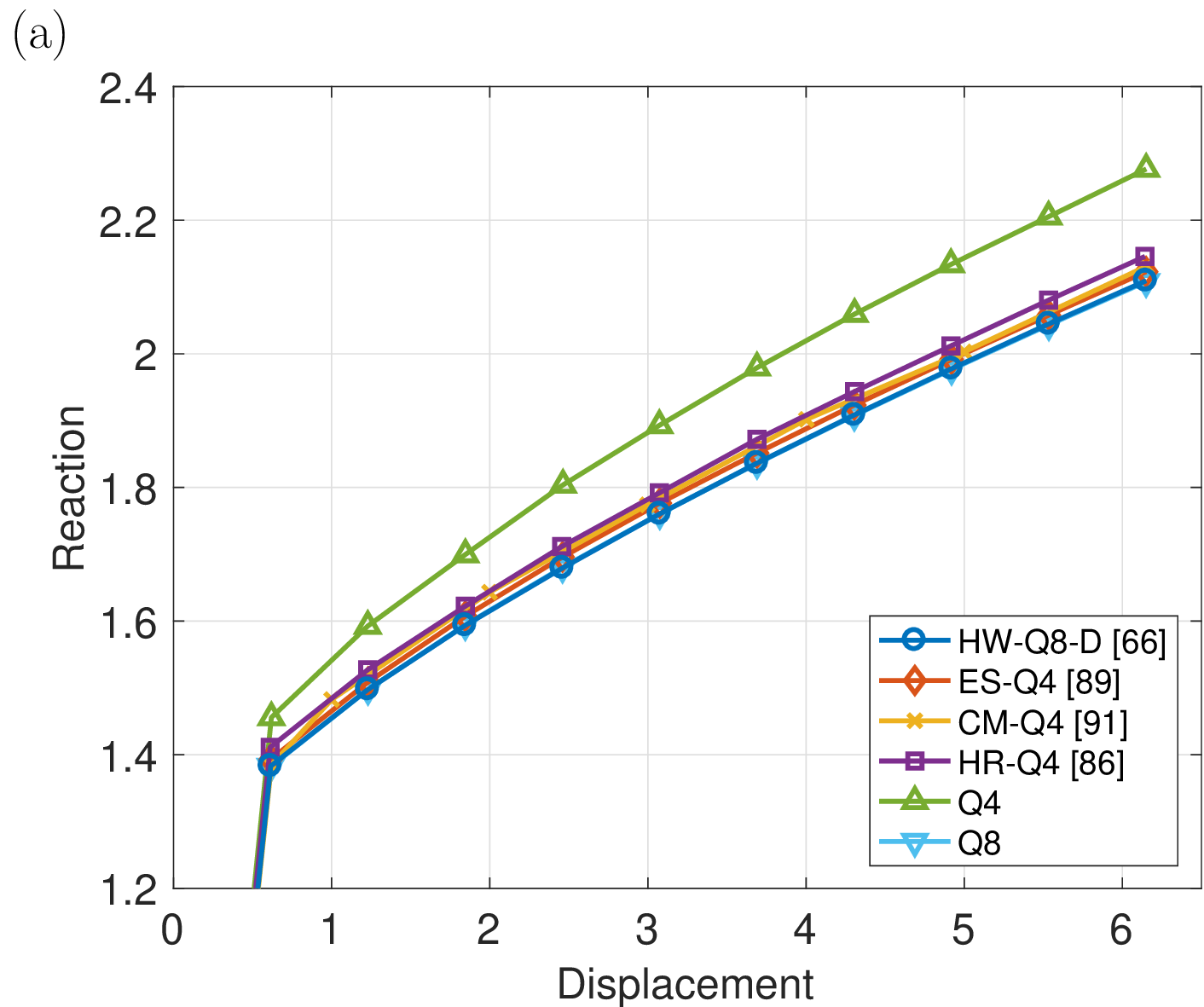}
	\includegraphics[clip = true, scale = 0.55]{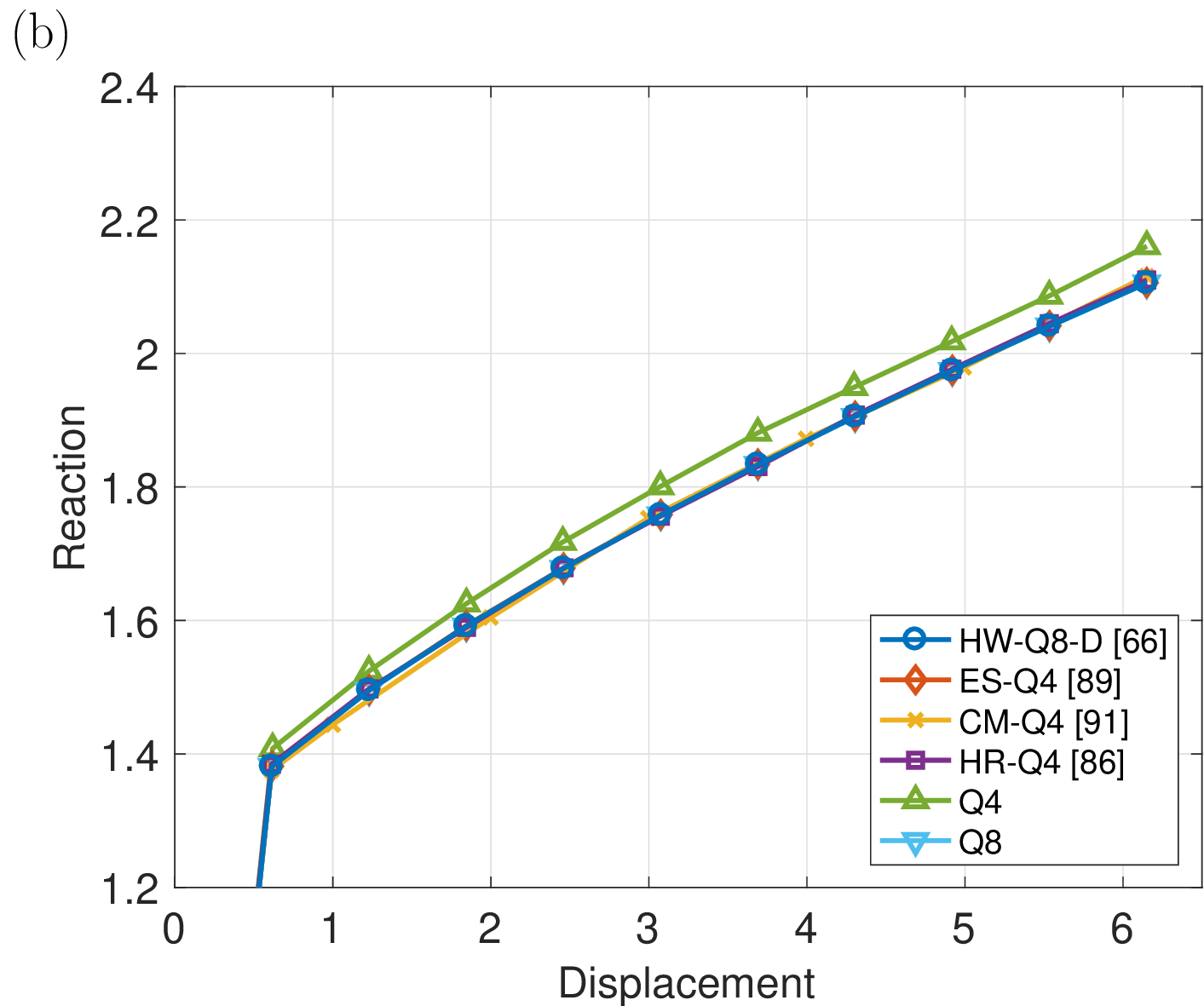}
	\caption{Extension of a perforated plate: load-displacement curve for (a) coarser mesh, comprising $6\times 12$ finite elements, and (b) finer mesh, comprising $19\times 38$ finite elements.}
	\label{fig:plate_response}
\end{figure}
The second benchmark of interest deals with a thin perforated plate, stretched along its longitudinal axis (for instance, see \cite{deSouza_Peric_2008}). The membrane has width $W = 20$, height $H = 36$, and thickness $t = 1$, and presents a circular hole of radius $R = 5$ at its center. The analysis is conducted under displacement control, prescribing a vertical displacement on the top and bottom edges. A displacement increasing up to maximum value $u\txtsub{max} = 6.15$ is applied. The material constituting the membrane is characterized by the following material parameters: Young's modulus $E = 70$, Poisson's ratio $\nu = 0.2$, yield stress $\sigma\y = 0.243$, and linear isotropic hardening modulus $k\iso = 0.2$. Because of problem symmetry, only one quarter of the plate is modeled by applying the appropriate boundary condition on the symmetry edges. The problem is solved considering two finite element meshes, defined in terms of the number of radial and circumferential subdivisions, $n\txtsub{r}$ and $n\txtsub{c}$ respectively, and labelled as $n\txtsub{r} \times n\txtsub{c}$. In particular, a coarser mesh~$6 \times 12$ and a finer mesh~$19 \times 38$ are analyzed. Geometry, boundary conditions and the former finite element mesh are illustrated in Figure~\ref{fig:plate_geometry}.

Figures~\ref{fig:plate_response}(a) and \ref{fig:plate_response}(b) report the vertical reaction (on the computational domain) versus the applied displacement for the coarser and the finer meshes, respectively. Apart for the displacement-based linear quadrilateral, all the curves are very close to each other for the analyses on the coarser mesh. Similarly to Cook's membrane problem, HW formulation is practically at convergence, whereas ES and CM formulations exhibit slightly superior results than HR formulation. When considering the finer mesh, all the mixed quadrilaterals are converged and the relevant curves coincide.

\section{Conclusions and perspectives}\label{s:conclusions}
The present work has proposed a review of mixed finite element methods for the analysis of bidimensional inelastic structures. Motivated by the poor performances of standard displacement-based formulations when considering inelastic materials, mixed formulations represent an opportunity to design high-performance finite elements by the introduction of independent additional fields, such as stresses, strains and enhanced strains. Though the main subject covered in the exposition deals with elastoplasticity, the more general framework of generalized standard materials has been adopted. Accordingly, exploiting the (incremental) energetic structure the material is equipped with, a unified mathematical setting has been delivered for investigating several potential variational formulations, as a starting point for the derivation of the relevant discrete problems. Some attention has been devoted to the selection of suitable interpolations for the unknown field with the aim of a general overview. Deeper analysis has been devoted to the numerical solution strategy adopted for the nonlinear discrete problem resulting from the variational formulation under investigation, a key aspect for the accuracy, robustness and efficiency of whatever computational method for inelasticity. Necessary and sufficient conditions for the numerical stability of the formulation in linear problems have been explored, in order to make clear the compatibility constraints to be fulfilled by the interpolation spaces of the unknown fields. Finally, numerical performances of several mixed formulations available in the literature have been explored in a few classical benchmark problems.

A remarkable literature deals with the treated topic, both considering its methodological aspects and the design of finite element formulations. Despite the consolidated variational foundations that can be relied upon, such continuous and ongoing research effort proves the need of further developments, especially (i) in exploiting the interpolation flexibility of mixed methods for capturing peculiar features of the behavior of inelastic structures and (ii) in exploring new numerical solution strategies for the nonlinear structural discrete problem. As shown in recent contributions, the two aspects are closely related to each other, because of the increasing number of element interpolation parameters required for achieving accuracy even on coarse meshes. In the former aspect, a promising direction is represented by interpolations that are discontinuous within each finite element, so to mimic the introduction of a \emph{macro-element} concept. In the latter aspect, mathematical programming strategies have nowadays reached firm maturity to be reliably exploited in (nonlinear) structural mechanical applications. In addition, numerical algorithms might benefit from the incremental energy minimization format that still appears not conventional in elastoplasticity. Attention should also be focused on the numerical stability analysis, which is well established from both theoretical and numerical standpoints, but is so far limited to linear elastic behavior. In closing, apart for the extension to the finite strain regime, it is worth recalling the essential assumption, here considered, of hardening behavior (i.e., of strictly convex incremental energy). The extension of such results to softening behavior represents a central topic of research in computational mechanics.

\vspace{2ex}
\footnotesize{\noindent \textbf{Acknowledgements} The author expresses his sincere gratitude to Professor Paolo Bisegna for valuable comments and stimulating discussions on the present work.}

\bibliographystyle{plainnat} 
\bibliography{quad}

\end{document}